\documentclass[a4paper,12pt]{amsart}

\usepackage{amsthm,amsmath,amsfonts,amssymb}
\usepackage{color}
\newcommand{\Red}{}
\renewcommand{\Red}[1]{{\color{red} #1}}

\newlength{\hchng}
\newlength{\vchng}
\newtheorem{thm}{Theorem}[section]

\newtheorem{defn}[thm]{Definition}

\newtheorem{preremark}[thm]{Remark}
\newenvironment{remark}{\begin{preremark}\rm}{\medskip \end{preremark}}
\numberwithin{equation}{section}

\DeclareMathOperator{\dv}{div}

\begin{document}

\title{On the stationary Navier-Stokes flow with isotropic streamlines in all latitudes on a sphere or a 2D hyperbolic space.}
\author{Chi Hin Chan and Tsuyoshi Yoneda}
\maketitle
\begin{center}
Department of Applied Mathematics, National Chiao Tung University, Hsinchu, Taiwan 30010, ROC \\
\end{center}
\begin{center}
Department of Mathematics, Hokkaido University, Sapporo 060-0810, Japan
\end{center}

\bibliographystyle{plain}
\noindent{\bf Abstract:} In this paper, we show the existence of real-analytic stationary Navier-Stokes flows with isotropic streamlines in all latitudes in some simply-connected flow region on a rotating round sphere. We also exclude the possibility of having a Poiseuille's flow profile to be one of these stationary Navier-Stokes flows with isotropic streamlines.
When the sphere is replaced by a $2$-dimensional hyperbolic space, we also give the analog existence result for stationary parallel laminar Navier-Stokes flows along a circular-arc boundary portion of some compact obstacle in the $2$-D hyperbolic space. The existence of stationary parallel laminar Navier-Stokes flows along a straight boundary of some obstacle in the $2$-D hyperbolic space is also studied. In any one of these cases, we show that a parallel laminar flow with a Poiseuille's flow profile ceases to be a stationary Navier-Stokes flow, due to the curvature of the background manifold.

\vskip0.3cm \noindent {\bf Keywords:}
Navier-Stokes equation, Riemannian manifold, streamlines

\vskip0.3cm \noindent {\bf Mathematics Subject Classification:}
76D03, 76D05, 53Z05
\section{Introduction: main results of this paper}\label{IntroductionSEC}

\noindent{\bf Part I : About the study of stationary Navier-Stokes flow with isotropic streamlines on a round sphere}\\

Yoden and Yamada \cite{YY} studied a two dimensional flow on a rotating sphere without any external force.
They investigated the morphology of the stream-function and the vorticity field at several rotation rates.
As the rotation rate increases, the temporal evolution of the flow field changes drastically. In particular, they observed in \cite{YY} that an easterly circumpolar vortex starts to appear in high latitudes, and that the flow field becomes anisotropic in all the latitudes.
However, to the best of our knowledge, it appears that these observations made in \cite{YY} have not been investigated by means of a pure mathematical approach in the research literature (see however a recent article \cite{Cheng} for a analytical study of time-dependent solutions to the Euler equation on a rotating round sphere).
In this paper we investigate existence of solutions to the following $2$-dimensional Stationary Navier-Stokes equation on a rotating sphere, which is written in the language of differential $1$-forms on the sphere.
\begin{equation}\label{NavierStokesintroduction}
\begin{split}
\nu ((-\triangle)u^* -2 Ric (u^*) ) +\beta \cos (ar) * u^* +[\overline{\nabla}_{u}u]^* + dP & = 0, \\
d^* u^* & = 0 .
\end{split}
\end{equation}
In \eqref{NavierStokesintroduction},$*$ is the Hodge-star operator sending the space of differential $1$-forms into itself. Intuitively, the action $u^* \to * u^*$ should be interpreted as the rotational action $u \to u^{\perp}$ by the angle $\frac{\pi}{2}$ in the anti-clockwise direction. So, it turns out that the term $\beta \cos (ar) * u^*$ which appears in \eqref{NavierStokesintroduction} represents the effect upon the velocity field $u$ due to the rotation of the sphere with some \emph{constant} rotational speed $\beta \geqslant 0$. The operator $d^*$ acting on the space of smooth differential $1$-forms on the sphere should be interpreted as the operator $- div$ acting on the space of smooth vector fields on the sphere. Notice also that the viscosity term of \eqref{NavierStokesintroduction} is represented as the linear combination of the standard Hodge Laplacian $(-\triangle ) = dd^* + d^*d$ acting on the space of $1$-forms on the sphere and $-2Ric (u^*)$, with $Ric$ to be the standard Ricci tensor on the sphere in differential geometry. Notice that the multiplicative \emph{constant} $\nu > 0$ stands for the viscosity coefficient for the stationary Navier-Stokes flows governed by equation \eqref{NavierStokesintroduction}. Moreover, the symbol $\overline{\nabla}$ as appears in the nonlinear convection term of \eqref{NavierStokesintroduction} is the Levi-Civita connection (which operates on the space of smooth vector fields) induced by the intrinsic Riemannian geometry of the sphere. \\

\begin{remark}
For the precise definitions of the operators $(-\triangle)$, $d^*$, $*$, and the Levi-Civita connection $\overline{\nabla}$ in the general Riemannian manifold setting, we refer our readers to \textbf{Definitions \ref{VolumeformDef}, \ref{HodgeDefinition} , \ref{LeviCivitadefinition}} in \textbf{Section \ref{BasicsphereSEC}} and \textbf{Definition \ref{dstarDef}} in \textbf{Section \ref{Spheremainsection}} of this paper. For the intuitive meaning of these operators, we refer our readers to \textbf{Remarks \ref{remarkgradient}, \ref{curlremark}, and \ref{divergenceremark}} in \textbf{Section \ref{BasicsphereSEC}} of this paper.
\end{remark}

\begin{remark}\label{1formremark}
Here, we would like to explain the meaning of the notation $u^*$ which appears in \eqref{NavierStokesintroduction}. In the case of a general $N$-dimensional Riemannian manifold $M$ equipped with a Riemannian metric $g(\cdot , \cdot )$, for any given smooth vector field $u$ on $M$, we can always construct its associated $1$-form, namely $u^* \in C^{\infty}(T^*M)$, in accordance with the following relation.
\begin{equation}
u^* = g(u , \cdot ).
\end{equation}
Since by definition, the Riemannian metric $g(\cdot,\cdot )$ on $M$ is actually a smoothly varying way of assigning to each $p\in M$ the positive definite inner product $g_p(\cdot , \cdot )$ on the tangent space $T_pM$ of $M$ at $p$, the above construction of the associated $1$-form $u^*$ on $M$ for each smooth vector field $u$ on $M$ actually provides a one-to-one correspondence between the space of smooth vector fields and the space of smooth $1$-forms on $M$.
\end{remark}

\begin{remark}\label{NSremark}
The structure of \eqref{NavierStokesintroduction} is based on the standard stationary Navier-Stokes equation as given in \eqref{NavierStokeshyperbolic} with an extra term $\beta \cos (ar) * u^*$ being included in order to take the effect due to the rotational action of the sphere into our account. The structure of the Navier-Stokes equation as specified in \eqref{NavierStokeshyperbolic} first appeared in the Historical work \cite{ EbinMarsden} by D. Ebin and J. Marsden. Since then, equation \eqref{NavierStokeshyperbolic} (with the extra term $\partial_{t}u^*$ in the time-dependent case) had been accepted as the standard form of the Navier-Stokes equation being written on a general finite dimensional Riemannian manifold $M$. Having a discussion about the research literature of the Navier-Stokes equation on a general Riemannian manifold   is out of the scope of this paper. However, for further historical remarks, we refer our readers to the textbook \cite{MichaelTaylor} by M. Taylor, and also some recent works on this subject matter such as the work \cite{ DindosMitrea} by M. Dindos and M. Mitrea, the work \cite{nonuniquenesshyperbolic} by the first author and M. Czubak, and work \cite{Khesin} by B. Khesin and G. Misiolek for instance.
\end{remark}

\noindent
The first goal of our paper is to show that there exists a stationary Navier-Stokes flow solving equation \eqref{NavierStokesintroduction} (the same stationary flow for any rotation rates) with isotropic streamlines in all the latitudes on a (rotating) sphere, and with given boundary values near the north polar of the sphere.
Moreover, we show that there is \emph{no} stationary flow composed by a quadratic profile (Poiseuille flow profile).


Here, we just mention that in the case when the background manifold is just the standard Euclidean $2$-dimensional space $\mathbb{R}^2$, the analog investigation of the behavior of parallel laminar Navier-Stokes flows around a cicuclar arc boundary portion of an obstacle in $\mathbb{R}^2$ has been carried out in a recent paper \cite{Y} by the second author by means of  elementary method. \\

\noindent
Let us to be more precise. Consider the $2$-dimensional space form $S^2(a^2)$ of positive sectional curvature $a^2$, which can be realized as the standard sphere $\{ (x_1 ,x_2 , x_3 ) : x_{1}^2 + x_{2}^2 + x_{3}^2 = \frac{1}{a^2}       \}$ of radius $\frac{1}{a}$ in $\mathbb{R}^3$. Consider the selected point $O = (0,0,1) \in  S^2(a^2)$, which can be regarded to be the North Pole of the sphere $S^2(a^2)$, and introduce the standard \emph{normal polar coordinate system} $(r,\theta )$ about the based point $O$ on $S^2(a^2)$ via the exponential map $\exp_{O} : T_{O}S^2(a^2) \rightarrow S^2(a^2)$ (See \textbf{Definition \ref{expsphereDef}} for the precise meaning of $\exp_O$, and \textbf{Definition \ref{PolarDefnSphere}} for the precise meaning of $(r,\theta)$ on $S^2(a^2)$ ).\\

\noindent
The first goal of this paper is to study the existence and non-existence of locally defined parallel laminar Stationary Navier-Stokes flow on some local exterior region near the circular-arc \emph{boundary portion} of some compact obstacle $K$ in $S^2(a^2)$.
As a preparation for the statement of Theorem \ref{Noparallel}, let $K$ to be a compact set in $S^2(a^2)$ whose boundary $\partial K$ has a \emph{circular-arc boundary portion}. $K$ will represent an obstacle around which stationary Navier-Stokes flows occur. We will study stationary Navier-Stokes flows in some \emph{simply-connected} exterior region whose boundary shares the same circular-arc boundary portion with $\partial K$. More precisely, suppose that, for some positive numbers $\delta \in (0, \frac{\pi}{a})$ and $\tau \in (0, 2\pi)$, $\partial K$ contains the following circular-arc portion
\begin{equation}
C_{\delta , \tau} = \{p\in S^2(a^2) : d (p, O) = \delta , 0 < \theta (p) < \tau \}.
\end{equation}
That is, we have $C_{\delta , \tau } \subset \partial K$. Here, the symbol $d(p,O)$ stands for the geodesic distance between $p$ and $O$ on the sphere $S^2(-a^2)$. For convenience, we will assume also that $K$ is confined within the closed geodesic ball $\overline{B_{O}(\delta)} = \{ p \in S^2(a^2) : d(p,O) \leqslant \delta\}$. That is, we have $K \subset \overline{B_O(\delta )}$. Then, Consider the following \emph{sector-shaped} open region $R_{\delta , \tau , \epsilon_{0}}$ in $S^2(a^2)-K$
\begin{equation}\label{exteriorregion}
R_{\delta , \tau , \epsilon_0 } = \{p\in S^2(a^2): \delta<d(p,O)< \delta + \epsilon_{0} , 0 < \theta (p) < \tau \},
\end{equation}
with $\epsilon_0$ to be any small positive number which satisfies $\delta + \epsilon_{0} < \frac{\pi}{a}$.
Let $\{\frac{\partial }{\partial r} , \frac{\partial }{\partial \theta }\}$ to be the natural coordinate frame induced by the normal polar coordinate system $(r,\theta)$ on $S^2(a^2)$ (see \textbf{Definition \ref{coordFRAMEsphere}} for the precise meaning of $\{\frac{\partial }{\partial r} , \frac{\partial }{\partial \theta }\}$ on $S^2(a^2)$ ). Then, the following assertion holds:
\begin{thm}\label{Noparallel}
Consider the simply-connected exterior region $R_{\delta , \tau , \epsilon }$ as defined in \eqref{exteriorregion}.
For the quadratic profile (Poiseuille flow profile) $h(\lambda ) = \alpha_1 \lambda - \frac{\alpha_2}{2} \lambda^2$, with both $\alpha_1 > 0$ and $\alpha_2 > 0$, there does not exist any parallel laminar flow in the form of $u = -h(r-\delta ) \frac{a}{\sin (ar)}\frac{\partial}{\partial \theta }$ which solves the stationary Navier-Stokes equation \eqref{NavierStokesintroduction} in the region $R_{\delta , \tau , \epsilon_0 }$ as specified in \eqref{exteriorregion}.
\end{thm}

\begin{remark}
In the statement of Theorem \ref{Noparallel}, a velocity field in the form of $u = -h(r-\delta ) \frac{a}{\sin (ar)}\frac{\partial}{\partial \theta }$ is called a parallel laminar flow along the circular-arc boundary portion $C_{\delta , \tau}$ of the given obstacle $K$ in $S^2(a^2)$, exactly because each streamline of such a velocity field $u$ is by itself a circular arc which keeps a constant geodesic distance from the boundary portion $C_{\delta , \tau}$. In general, a smooth velocity field $u$ as specified on some open flow region near some smooth boundary portion $\Gamma$ of some obstacle $K$ in a $2$-D manifold $M$ is called a parallel laminar flow along the smooth boundary portion $\Gamma$, if every single streamline of $u$ keeps a constant geodesic distance from $\Gamma$. This notion of parallel laminar flows as employed here is consistence with the definition of parallel laminar flow as given in the recent work \cite{Y} by the second author.
\end{remark}

\noindent
Here, let us say something about the idea of the proof of Theorem \ref{Noparallel}, which will be given in \textbf{Section \ref{Spheremainsection}} in details. The first step of the argument is to obtain an explicit formula of the following expression, for $u = -h(r-\delta ) \frac{a}{\sin (ar)}\frac{\partial}{\partial \theta }$.
\begin{equation}\label{dclosedtest}
d \{  \nu ((-\triangle)u^* -2 Ric (u^*) ) + \overline{\nabla}_{u}u^*  \},
\end{equation}
In \eqref{dclosedtest}, the symbol $u^*$ stands for the associated $1$-form of the vector field $u$, which is defined by the relation $u^* = g(u,\cdot )$, with $g(\cdot, \cdot )$ to be the Riemannian metric of the sphere $S^2(a^2)$.
Also, the symbol $d$ stands for the exterior differential operator $d : C^{\infty}(T^*M) \rightarrow C^{\infty}(\wedge^2 T^*M)$ sending smooth $1$-forms to smooth $2$-forms on a Riemannian mainfold $M$, which in case of Theorem \ref{Noparallel}, is taken to be $M = S^2(a^2)$.\\
So, \textbf{Step 1} and \textbf{Step 2} in \textbf{Section \ref{Spheremainsection}} are carried out in order to compute the terms $(-\triangle )u^*$, and the nonlinear convection term $\overline{\nabla}_{u}u^*$ involved in equation \eqref{NavierStokesintroduction}. These efforts eventually lead to the following representation formula of expression \eqref{dclosedtest} in \textbf{step 3} of \textbf{Section \ref{Spheremainsection}}.

\begin{equation}\label{ONE}
\begin{split}
d \big \{ \nu (-\triangle u^* -2 a^2 u^*) + \overline{\nabla}_{u}u^*     \big \}
& = \nu \bigg \{ h'''(r-\delta ) \frac{\sin (ar)}{a}  + 2 h''(r-\delta ) \cos (ar) \\
& + a h'(r-\delta ) (\sin (ar) -\frac{1}{\sin (ar)}) \\
& + a^2 h(r-\delta ) \cos (ar)(2+ \frac{1}{\sin^2(ar)}) \bigg \} dr\wedge d\theta .
\end{split}
\end{equation}
Moreover, since $\beta\cos(a r) *u^*=\beta h(r-\delta)\cos (ar)\frac{\sin(ar)}{a}dr$, we can easily get
\begin{equation}\label{TWO}
d\{\beta \cos r *u^*\}=0.
\end{equation}
Then, expressions \eqref{ONE} and \eqref{TWO} together will imply the following representation formula for $d \{ \nu (-\triangle u^* -2 a^2 u^*) + \beta\cos(a r) *u^* + \overline{\nabla}_{u}u^*     \}$.
\begin{equation}\label{generalexpINTRO}
\begin{split}
& d \big \{ \nu (-\triangle u^* -2 a^2 u^*) + \beta\cos(a r) *u^* + \overline{\nabla}_{u}u^*   \big \} \\
& = \nu \bigg \{h'''(r-\delta ) \frac{\sin (ar)}{a}  + 2 h''(r-\delta ) \cos (ar) \\
& + a h'(r-\delta ) (\sin (ar) -\frac{1}{\sin (ar)}) \\
& + a^2 h(r-\delta ) \cos (ar)(2+ \frac{1}{\sin^2(ar)}) \bigg \} dr\wedge d\theta .
\end{split}
\end{equation}
The representation formula \eqref{generalexpINTRO} is valid for \emph{any smooth functions} $h$ defined on any open interval around the point $\delta$.
So, the representation formula \eqref{generalexpINTRO} is the key tool which allows us to decide whether $u = -h(r-\delta ) \frac{a}{\sin (ar)}\frac{\partial}{\partial \theta }$ is a solution to \eqref{NavierStokesintroduction} or not. This is because, in accordance with the basic knowledge in differential topology, the \emph{vanishing} of
$d \{ \nu (-\triangle u^* -2 a^2 u^*) + \beta\cos(a r) *u^* + \overline{\nabla}_{u}u^*     \}$ over some simply-connected open region (say for instance $R_{\delta , \tau ,\epsilon_0}$ as specified in \eqref{exteriorregion}) near the circular-arc portion of the boundary of the obstacle $K$ will immediately imply the existence of some locally defined smooth pressure $P$ which solves \eqref{NavierStokesintroduction} on the same simply-connected open region. In other words, this basic idea of \textbf{$d$-closed implies $d$-exact on simply-connected region} in differential topology allows us to reduce our problem to the one of testing whether
or not $d \{ \nu (-\triangle u^* -2 a^2 u^*) + \beta\cos(a r) *u^* + \overline{\nabla}_{u}u^*     \}$ totally vanishes on some prescribed simply-connected open region near the circlar-arc portion of $\partial K$. So, in \textbf{Step 3} of \textbf{Section \ref{Spheremainsection}}, we finish the proof of Theorem \ref{Noparallel} by showing, case by case, that for the quadratic function $h(\lambda ) = \alpha_1 \lambda - \frac{\alpha_2}{2} \lambda^2$, the expression $d \{ \nu (-\triangle u^* -2 a^2 u^*) + \beta\cos(a r) *u^* + \overline{\nabla}_{u}u^*     \}$ would never vanish identically on the simply-connected region $R_{\delta , \tau , \epsilon_0}$ as specified in \eqref{exteriorregion}, \emph{no matter how small the positive number $\epsilon_0$ would be}. In this way, we get a neat and clean argument showing that $u = -h(r-\delta ) \frac{a}{\sin (ar)}\frac{\partial}{\partial \theta }$ is not a solution to \eqref{NavierStokesintroduction}, for any quadratic profile $h(\lambda ) =  \alpha_1 \lambda - \frac{\alpha_2}{2} \lambda^2$, with $\alpha_1 >0$, and $\alpha_2 > 0$.\\

\noindent
As a by product of this approach which we used in the proof of Theorem \ref{Noparallel}, we also make another important observation that the function as appears in the right hand side of \eqref{generalexpINTRO} is a third order linear differential operator $L$ acting on the unknown function $Y (r) = h (r-\delta )$. In other words, for a general parallel laminar flow in the form of $u = -Y(r) \frac{a}{\sin (ar)}\frac{\partial}{\partial \theta }$,  the vanishing of $d \{ \nu (-\triangle u^* -2 a^2 u^*) + \beta\cos(a r) *u^* + \overline{\nabla}_{u}u^*     \}$ over a certain simply-connected region near the circular-arc portion of $\partial K$ is equivalent to saying that the unknown function $Y(r) = h(r-\delta )$ will solve the following $3$-order linear ODE over some interval $[\delta , \delta + \epsilon )$.

\begin{equation}\label{thirdorderODE}
\begin{split}
0 & = Y'''(r) \frac{\sin (ar)}{a}  + 2 Y''(r) \cos (ar) \\
& + a Y'(r) (\sin (ar) -\frac{1}{\sin (ar)}) + a^2 Y(r) \cos (ar)(2+ \frac{1}{\sin^2(ar)}) .
\end{split}
\end{equation}
Since all the coefficient functions involved in the above $3$-order ODE are all real-analytic over the interval $(0, \frac{\pi}{a})$, the basic existence theorem (Theorem \ref{ODEtheorem} in \textbf{Section \ref{Spheremainsection}}) in the theory of linear ODE will ensure, for the prescribed initial data $Y(\delta ) = 0$, $Y'(\delta) = \alpha_{1}$, and $Y''(\delta ) = - \alpha_{2}$,  the existence of a unique smooth function $Y (r) = h(r-\delta) $ over $[\delta , \frac{\pi}{a})$ which solves equation \eqref{thirdorderODE} on $(0, \frac{\pi}{a})$, and which at the same time turns out to be real-analytic on $[\delta , \frac{\pi}{a} )$. This basically leads to our second basic theorem (Theorem \ref{existenceEasy}), which says that, for any prescribed constants $\alpha_1>0 $, $\alpha_2>0 $, there exists a unique smooth function $Y \in C^{\infty}([\delta , \frac{\pi}{a} ))$ satisfying $Y(\delta ) =0$, $Y'(\delta ) = \alpha_1$, and $Y''(\delta ) = -\alpha_2$, which turns out to be real-analytic over $[\delta , \frac{\pi}{a})$, such that the associated parallel laminar flow $u = -Y(r) \frac{a}{\sin (ar)}\frac{\partial}{\partial \theta }$ will be a solution to \eqref{NavierStokesintroduction} on some simply-connected open region near the circular-arc portion of the boundary of the compact obstacle $K$.

\begin{thm}\label{existenceEasy}
Consider the space form $S^2(a^2) = \{ (x_1 ,x_2 , x_3 ) : x_{1}^2 + x_{2}^2 + x_{3}^2 = \frac{1}{a^2}       \}$, with $a > 0$. Let $O \in S^2(a^2)$ to be a selected based point, and let $(r, \theta )$ to be the normal polar coordinate system on $S^2(a^2)$ about the based point $O$, which is introduced through the standard exponential map $\exp_{O} : \{ v \in T_OS^2(a^2) : \|v\| < \frac{\pi}{a} \} \rightarrow S^2(a^2)$.\\
Consider a fixed positive number $\delta \in (0 ,\frac{\pi}{a})$, and let $K$ to be some \emph{compact} region which is a subset of $\{p\in S^2(a^2): d(p,O)\leq \delta \}$, and which plays the role of an obstacle in $S^2(a^2)$. Suppose further that for some positive number $\tau \in (0, 2\pi )$, the circular arc $C_{\delta , \tau} = \{p\in S^2(a^2) : d (p, O) = \delta , 0 < \theta (p) < \tau \}$ constitutes a smooth boundary portion of $\partial K$ in that
\begin{equation}
\{p\in S^2(a^2) : d (p, O) = \delta , 0 < \theta (p) < \tau \} \subset \partial K .
\end{equation}
Then, it follows that, for any prescribed positive numbers $\alpha_1 > 0$, $\alpha_2 > 0$, there exists some unique smooth function $Y : [\delta , \frac{\pi}{a}) \rightarrow \mathbb{R}$ satisfying $Y(\delta ) = 0$, $Y'(\delta  ) = \delta_1$, $Y''(\delta ) = -\alpha_2$, such that the associated parallel laminar flow $u = -Y(r) \frac{a}{\sin (ar)}\frac{\partial}{\partial \theta }$ will solve \eqref{NavierStokesintroduction} in the following simply connected exterior region
\begin{equation}\label{OMEGASPHERE}
\Omega_{\delta, \tau} = \{p\in S^2(a^2) : \delta < d (p, O) < \frac{\pi}{a} , 0 < \theta (p) < \tau  \} ,
\end{equation}
which shares the same circular-arc boundary portion $C_{\delta , \tau }$ with the boundary of the compact obstacle $K$. Moreover, such a unique smooth function $Y : [\delta , \frac{\pi}{a}) \rightarrow \mathbb{R}$ which we desire turns out to be actually \textbf{real-analytic} on $[\delta , \frac{\pi}{a})$.
\end{thm}

\begin{remark}
Indeed, according to the basic existence theorem (Theorem \ref{ODEtheorem}) in the O.D.E. theory, the conclusion of Theorem \ref{existenceEasy} will remain valid even if we generalize $h$ to $h(\lambda)=\alpha_0+\alpha_1\lambda-\frac{\alpha_2}{2}\lambda^2$, with three prescribed parameters $\alpha_0 \geq 0$, $\alpha_1 > 0$, and $\alpha_2 > 0$.
\end{remark}

\noindent
Indeed, Theorem \ref{existenceEasy} is a by-product which follows from an application of representation formula \eqref{generalexpINTRO} together with the use of the basic existence and uniqueness Theorem (Theorem \ref{ODEtheorem} in \textbf{Section \ref{Spheremainsection}}) for $3$-order ODE with real-analytic coefficients. The very brief and easy proof of Theorem \ref{existenceEasy} will be given in \textbf{Step 4} of \textbf{Section \ref{Spheremainsection}}.\\

\noindent
At a first glance, the conclusion in Theorem \ref{existenceEasy} seems to state a conclusion which is opposite to that of Theorem \ref{Noparallel}. However, this is \textbf{not} true at all, since Theorem \ref{Noparallel} only rules out the existence of stationary Navier-Stokes flow in the form of $u = -Y(r) \frac{a}{\sin (ar)}\frac{\partial}{\partial \theta }$, with the \emph{quadratic profile (Poiseuille flow profile)} $Y (r) =  \alpha_{1} (r-\delta ) - \frac{\alpha_2}{2} (r-\delta )^2$.
However, Theorem \ref{existenceEasy} says that, as long as \emph{higher order terms} beyond the quadratic power $(r-\delta )^2$ are allowed in the Taylor series expansion
of the unknown function $Y(r)$, the basic existence and uniqueness theory for $3$-order linear ODE will ensure the existence of a unique real-analytic function $Y$ for which $u = -Y(r) \frac{a}{\sin (ar)}\frac{\partial}{\partial \theta }$ will solve \eqref{NavierStokesintroduction} in some simply-connected open region near the circular-arc portion of $\partial K$. Before we leave Part I of the introduction, let us make another interesting remark here which helps us to relate the result in Theorem \ref{existenceEasy} with the observations made in \cite{YY} by Yoden and Yamada.

\begin{remark}
It is worthwhile to notice that, in Theorem \ref{existenceEasy}, the unique real-analytic function $Y$ on $[\delta , \frac{\pi}{a} )$ which makes  $u = -Y(r) \frac{a}{\sin (ar)}\frac{\partial}{\partial \theta }$ become a solution to \eqref{NavierStokesintroduction} is \emph{independent} of the constant rotational speed $\beta > 0$ of the round sphere $S^2(a^2)$. In other words, for the \emph{same} real analytic function $Y$ on $[\delta , \frac{\pi}{a} )$ satisfying equation \eqref{thirdorderODE} and the initial values $Y(\delta )=0$, $Y'(\delta )= \alpha_1$, and $Y''(\delta ) = -\alpha_2$, the velocity field $u = -Y(r) \frac{a}{\sin (ar)}\frac{\partial}{\partial \theta }$ can be realized as a stationary Navier-Stokes flow on the rotating sphere $S^2(a^2)$ with any prescribed rotational speed $\beta > 0$. This remark seems to be even more interesting when one compares the existence result as given in Theorem \ref{existenceEasy} with the Numerical experiment as carried out in the work \cite{YY} by Yoden and Yamada, according to which easterly circumpolar vortex starts to appear in high latitudes when the rotational speed $\beta > 0$ increases, and that the flow field becomes anisotropic in all the latitudes. So, by combining the result as given in Theorem \ref{existenceEasy} and the observation made by Yoden and Yamada in \cite{YY}, it is very tempting for one to speculate that a large rotational speed $\beta$ of the rotation of the sphere $S^2(a^2)$ about the North-South poles axis may have some considerable effect in \emph{stabilizing} the flow pattern of the (real analytic) Stationary Navier-Stokes flow $u = -Y(r) \frac{a}{\sin (ar)}\frac{\partial}{\partial \theta }$ as ensured by Theorem \ref{existenceEasy}.
\end{remark}

\noindent{\bf Part II : About the study of stationary parallel laminar Navier-Stokes flows around an obstacle in a hyperbolic manifold with constant negative sectional curvature.}\\

\noindent
In the second part of the introduction, we will study the existence and non-existence of Stationary Navier-Stokes flows with circular-arc streamlines around some compact obstacle $K$ in a $2$-dimensional space-form $\mathbb{H}^2(-a^2)$ of constant negative sectional curvature $-a^2$. On such a $2$-dimensional space-form $\mathbb{H}^2(-a^2)$ with constant negative sectional curvature, which will also be called the $2$-dimensional hyperbolic space with constant sectional curvature $-a^2 < 0$, we will study the following stationary Navier-Stokes equation on $\mathbb{H}^2(-a^2)$, which again is formulated in terms of the language of differential $1$-forms on $\mathbb{H}^{2}(-a^2)$ instead of that of smooth vector fields (Again, see \textbf{Remark \ref{1formremark}}).

\begin{equation}\label{NavierStokeshyperbolic}
\begin{split}
\nu ((-\triangle)u^* -2 Ric (u^*) ) +[\overline{\nabla}_{u}u]^* + dP & = 0, \\
d^* u^* & = 0 .
\end{split}
\end{equation}
In equation \eqref{NavierStokeshyperbolic}, $\overline{\nabla}$ stands for the standard Levi-Civita connection (covariant derivative) acting on the space of smooth vector fields on $\mathbb{H}^2(-a^2)$. Again, the operator $d^*$ sending smooth $1$-forms into the space of smooth functions on $\mathbb{H}^2(-a^2)$ is interpreted as $-div$. The viscosity term in \eqref{NavierStokeshyperbolic} consists of two terms, namely $(-\triangle )u^*$ and $-2 Ric (u^*)$, where $(-\triangle ) = dd^* + d^*d$ is the standard Hodge Laplacian acting on the space of $1$-forms on $\mathbb{H}^2(-a^2)$ and $Ric$ is the standard Ricci tensor with respect to the Riemannian metric of the hyperbolic manifold $\mathbb{H}^2(-a^2)$.\\

\noindent
Here, let us explain a little bit about why we would like to extend our study of stationary parallel laminar flows around circular-arc boundary portion of some obstacle to the case in which the background space is a hyperbolic manifold $\mathbb{H}^2(-a^2)$ of constant negative sectional curvature $-a^2$. Indeed, it is understandable that a P.D.E. specialist with a more practical mind (or a scientist working in the area of fluid dynamics in general) may find it difficult to comprehend or appreciate the meaning or significance of studying Navier-Stokes flow on such a hyperbolic manifold. This kind of negative attitude towards the study of Navier-stokes flows on hyperbolic manifolds is understandable because a classical theorem due to Hilbert \cite{Hilbert} and Efimov \cite{Efimov} states that a complete $2$-dimensional Riemannian manifold with negative sectional curvature to be bounded above by a negative constant cannot be isometrically embedded  into the standard Euclidean space $\mathbb{R}^3$ equipped with the standard Euclidean metric. However, if one looks at this issue from the view point of pure mathematics, a hyperbolic space $\mathbb{H}^2(-a^2)$ with constant negative sectional curvature $-a^2 < 0$ is exactly the ''negative counterpart'' of the sphere $S^2(a^2)$ with radius $\frac{1}{a}$, which is the space-form of constant positive sectional curvature $a^2$. So, methodically speaking, if a set of methods, such as those we used in \textbf{Section \ref{Spheremainsection}}, works well in the study of parallel laminar flows around an obstacle in the round sphere $S^2(a^2)$, one expects that the same set of methods, once being adopted to the setting of a hyperbolic space $\mathbb{H}^2(-a^2)$, should work equally well and should yield equally interesting results analogical to those being obtained in the spherical case $S^2(a^2)$. Indeed, if one tries to compare the mathematical content of \textbf{Section \ref{Spheremainsection}}, which contains the proofs of Theorems \ref{Noparallel} and \ref{existenceEasy}, with \textbf{Section \ref{hyperMAINSECTION}}, which contains the proof of Theorem \ref{existenceHyperbolic}, the similarities between the spherical case and the hyperbolic counterpart are striking. In order words, from the mathematical view-point, it is completely natural to state and prove Theorem \ref{existenceHyperbolic}, which is analogical to the results of Theorems \ref{Noparallel} and Theorem \ref{existenceEasy} in the spherical case.

\begin{thm}\label{existenceHyperbolic}
Consider the $2$-dimensional space form $\mathbb{H}^2(-a^2)$ of constant negative sectional curvature $-a^2$, with $a > 0$ to be given. Let $O \in \mathbb{H}^2(-a^2)$ to be a selected based point, and let $(r, \theta )$ to be the normal polar coordinate system on $\mathbb{H}^2(a^2)$ about the based point $O$, which is introduced through the standard exponential map $\exp_{O} : T_{O}\mathbb{H}^2(-a^2) \rightarrow \mathbb{H}^2(-a^2)$ (see \textbf{Definition \ref{exphyper}} for the precise meaning of $\exp_O$, and also \textbf{Definition \ref{polarDefn}} for the precise meaning of $(r,\theta )$ on $\mathbb{H}^2(-a^2)$ ).\\
Consider a fixed choice of positive number $\delta \in (0 , \infty )$, and let $K$ to be some \emph{compact} region which is a subset of $\{p\in \mathbb{H}^2(a^2): d(p,O)\leq \delta \}$, and which plays the role of an obstacle in $\mathbb{H}^2(-a^2)$. Here, $d(p,O)$ stands for the geodesic distance between $O$ and $p$ in $\mathbb{H}^2(a^2)$. Suppose further that for some positive number $\tau \in (0, 2\pi )$, the circular arc $C_{\delta , \tau} = \{p\in \mathbb{H}^2(-a^2) : d (p, O) = \delta , 0 < \theta (p) < \tau \}$ constitutes a smooth boundary portion of $\partial K$ in that
\begin{equation}
\{p\in \mathbb{H}^2(-a^2) : d (p, O) = \delta , 0 < \theta (p) < \tau \} \subset \partial K .
\end{equation}
Moreover, let $\frac{\partial}{\partial r}$, and $\frac{\partial}{\partial \theta }$ to be the two natural vector fields induced by the normal polar system $(r,\theta )$ on $\mathbb{H}^2(a^2)$ about the based point $O$, which are defined in \textbf{Definition \ref{FRAMEHyperDef}} of \textbf{Section \ref{HYPERBOLIDSECTION}}. Then, it follows that the following two assertions are valid.\\

\noindent
\textbf{Assertion I} For the quadratic profile $h(\lambda ) = \alpha_1 \lambda -\frac{\alpha_2}{2} \lambda^2$ with any prescribed
 constants $\alpha_1 > 0$, and $\alpha_2 > 0$, the velocity field $u = -h(r-\delta ) \frac{a}{\sinh (ar)} \frac{\partial}{\partial \theta }$ does not satisfy equation \eqref{NavierStokeshyperbolic} on any sector-shaped region $R_{\delta , \tau , \epsilon_0 }$ of $\mathbb{H}^2(-a^2)$ specified as follow, regardless of how small the positive number $\epsilon_0 > 0$ may be:
\begin{equation}\label{RegionHYPER}
R_{\delta , \tau , \epsilon_0 } = \{ p \in \mathbb{H}^2(-a^2) : \delta < d(p,O) < \delta + \epsilon_0, 0 < \theta (p) < \tau   \}.
\end{equation}

\noindent
\textbf{Assertion II} For any prescribed positive numbers $\alpha_1 > 0$, $\alpha_2 > 0$, there exists some unique smooth function $Y : [\delta , \infty ) \rightarrow \mathbb{R}$ satisfying $Y(\delta ) = 0$, $Y'(\delta  ) = \alpha_1$, $Y''(\delta ) = -\alpha_2$, such that the associated parallel laminar flow $u = -Y(r) \frac{a}{\sinh (ar)}\frac{\partial}{\partial \theta }$ will satisfy equation \eqref{NavierStokeshyperbolic} on the following simply connected exterior region:
\begin{equation}\label{OMEGAhyper}
\Omega_{\delta, \tau} = \{p\in \mathbb{H}^2(a^2) :  d (p, O) > \delta , 0 < \theta (p) < \tau  \} ,
\end{equation}
which shares the same circular-arc boundary portion $C_{\delta , \tau }$ with the boundary of the compact obstacle $K$. Moreover, such a unique smooth function $Y : [\delta , \infty ) \rightarrow \mathbb{R}$ which we desire turns out to be actually \textbf{real-analytic} on $[\delta , \infty )$.
\end{thm}

We just remark that the proof of Theorem \ref{existenceHyperbolic} will be given in \textbf{Section \ref{hyperMAINSECTION}}.

\noindent

The last mathematical result which we are going to give concerns the existence or non-existence of stationary parallel laminar Navier-Stokes flows along a geodesic representing the ''straight edge''-boundary of some obstacle in $\mathbb{H}^2(-a^2)$. In order to state this last mathematical result (Theorem \ref{parallelflowhyperbolicThm}) we need to consider another pairs of vector fields $\frac{\partial}{\partial \tau}$ and $\frac{\partial}{\partial s}$ on $\mathbb{H}^2(-a^2)$ as defined in expression \eqref{intuitivehyp} of \textbf{Section \ref{Cartesianhyperbolicsection}}.

\begin{thm}\label{parallelflowhyperbolicThm}
Let $\mathbb{H}^2(-a^2)$ to be the $2$-dimensional space form with constant negative sectional curvature $-a^2$. Consider the ''Cartesian coordinate system'' $\Phi : \mathbb{R}^2 \rightarrow \mathbb{H}^2(-a^2)$ on $\mathbb{H}^2(-a^2)$, which is constructed in \textbf{Section \ref{Cartesianhyperbolicsection}} of this paper, with the two geodesics $\tau \to c(\tau )$ and $s \to \gamma (s)$ playing the roles of the ''$x$-axis'' and the ''$y$-axis'' on $\mathbb{H}^2(-a^2)$ respectively (For the precise definition of $\Phi (\tau , s)$ and the roles played by the geodesics $c$ and $\gamma$, see expression \eqref{coordinatesystemhyper} and the discussion in \textbf{Section \ref{Cartesianhyperbolicsection}}). Here, the geodesic $\gamma$ playing the role the ''$y$-axis'' of the coordinate system $\Phi (\tau, s)$ will also represent the ''straight-edge'' boundary of the an obstacle $K$ occupying the infinite region lying on the ''left-hand side'' of $\gamma$ $($see \eqref{K} in \textbf{Section \ref{Cartesianhyperbolicsection}} for a precise definition of $K$$)$. Let $\frac{\partial}{\partial \tau }$ and $\frac{\partial}{\partial s } $ to be the natural vector fields on $\mathbb{H}^2(-a^2)$, which we define in expression \eqref{intuitivehyp} of \textbf{Section \ref{Cartesianhyperbolicsection}}. Then, we claim that the following two assertions hold.\\

\noindent
\textbf{Assertion I} For the quadratic function $h(\tau ) = \alpha_1 \tau - \frac{\alpha_2}{2}\tau^2$, with any given constants $\alpha_1 > 0$ and $\alpha_2 >0$, the parallel laminar flow given by $u = -h(\tau ) \frac{1}{\cosh (a\tau )} \frac{\partial}{\partial s }$ does not satisfy equation \eqref{NavierStokeshyperbolic} on any simply connected region $\Omega_{\tau_{0}} \subset \mathbb{H}^2(-a^2)-K$ in the following form.
\begin{equation}
\Omega_{\tau_{0}} = \{p \in \mathbb{H}^2(-a^2) : 0 < \tau (p) < \tau_0 \} = \Phi ((0 ,\tau_0 )\times \mathbb{R}),
\end{equation}
with $\tau_0 > 0$ to be an arbitrary constant, and the coordinate system $\Phi (\tau , s)$ on $\mathbb{H}^2(-a^2)$ to be the one given in \eqref{coordinatesystemhyper}. \\

\noindent
\textbf{Assertion II} For any prescribed positive numbers $\alpha_1 > 0$, and $\alpha_2 > 0$, there exists a unique smooth function $Y \in C^{\infty} \big ( [0, \infty )\big ) $ satisfying $Y(0)=0$, $Y'(0)= \alpha_1$, and $Y''(0) = -\alpha_2$ such that the associated parallel laminar flow
$u = -Y(\tau ) \frac{1}{\cosh (a\tau )} \frac{\partial}{\partial s}$ will solve equation \eqref{NavierStokeshyperbolic} in the whole exterior region $\mathbb{H}^2(-a^2)- K = \Phi \big ( (0,\infty )\times \mathbb{R} \big )$. Moreover, such a unique smooth function $Y : [0,\infty ) \rightarrow \mathbb{R}$ which we desire turns out to be actually \textbf{real-analytic} on $[0,\infty )$.
\end{thm}

The proof of Theorem \ref{parallelflowhyperbolicThm} will be given in \textbf{Section \ref{ProofofHypstraightSec}}. We also point out that, in the case when the background manifold is the round sphere $S^2(a^2)$, the study of Stationary parallel laminar Navier-Stokes flows along a \emph{geodesic} representing the boundary of some obstacle in $S^2(a^2)$ is essentially covered by Theorem \ref{Noparallel} and Theorem \ref{existenceEasy}, simply due to the fact that the great circle $\{p\in S^2(a^2) : d(p,O) = \frac{\pi}{2a}\}$ is a geodesic on $S^2(a^2)$ (and is the only possible type of geodesic on $S^2(a^2)$ up to the group of rotational isometries).

\section{Some Basic Geometry and Geometric Construction on the Space Form $S^{2}(a^{2})$ with Positive Constant Sectional Curvature $a^{2}$ }\label{BasicsphereSEC}

\noindent
As a preparation for the proofs of Theorems \ref{Noparallel} and \ref{existenceEasy} in \textbf{Section \ref{Spheremainsection}}, we will discuss first the \emph{normal polar coordinate system} and then the \emph{Hodge Laplacian} on the space form $S^{2}(a^{2})$ with positive sectional curvature $a^{2}$, with $a > 0$ to be a given constant. We would like to stress that all the material as presented in this section, as well as those in \textbf{Sections \ref{HYPERBOLIDSECTION}} and \textbf{\ref{Cartesianhyperbolicsection}}, pertains to the standard, fundamental working knowledge which is basic and well-known to researchers working in differential geometry, geometric analysis, or differential topology. However, it seems to us that many of our potential readers, including P.D.E. specialists working in the areas of Navier-Stokes equations, may not be familiar with those basic geometric language which we employ in this paper. Based on this consideration, we include some necessary geometric background here as a preparation for the differential geometric computation being done in the proofs of Theorems \ref{Noparallel} and \ref{existenceEasy} in the \textbf{Section \ref{Spheremainsection}}.  \\

\noindent{\bf The normal polar coordinate system $(r, \theta )$ on $S^2(a^2)$, with the induced natural moving frame $\{\frac{\partial}{\partial r} , \frac{\partial}{\partial \theta }\}$ on $S^2(a^2)$ .}\\

\noindent
Geometrically, the space form $S^{2}(a^{2})$ can be \emph{realized} as the sphere in $\mathbb{R}^{3}$ with radius $\frac{1}{a}$ and centered at the origin $(0,0,0)$ of $\mathbb{R}^{3}$. That is, we have the following identification

\begin{equation}
S^{2}(a^{2}) = \left\{ (x_1 ,x_2 , x_3 ) : x_{1}^2 + x_{2}^2 + x_{3}^2 = \frac{1}{a^2}      \right \} ,
\end{equation}
where the sphere $\{ (x_1 ,x_2 , x_3 ) : x_{1}^2 + x_{2}^2 + x_{3}^2 = \frac{1}{a^2}   \}$ is understood to be a submanifold of $\mathbb{R}^{3}$, equipped with the standard induced Riemannian metric $g(\cdot ,\cdot )$ inherted from the background Euclidean space $\mathbb{R}^3$. Before we introduce the normal polar coordinate system on $S^{2}(a^{2})$, we choose the point $O = (0,0, \frac{1}{a} ) \in S^{2}(a^{2})$, which is the North Pole of the sphere $S^{2}(a^{2})$, to be the reference point at which our normal polar coordinate system will be based. Let $T_{O}S^{2}(a^{2})$ to be the tangent space of $S^{2}(a^{2})$ at $O \in S^{2}(a^{2})$.
Intuitively, $T_{O}S^{2}(a^2)$ can be realized as the two dimensional plane in $\mathbb{R}^{3}$ that passes through $(0,0, \frac{1}{a})$ and that is parallel to $\{ (x,y,0) : x , y \in \mathbb{R}  \}$. So, it is natural to identify $T_{O}S^{2}(a^2)$ with $\mathbb{R}^{2} = \{ (x,y) : x,y \in \mathbb{R}   \}$. Recall that
the tangent space $T_{p}M$ of a manifold $M$ at $p \in M$ is, by definition, a vector space of the same dimension as that of the manifold $M$. \\

\noindent
Here, we need to introduce the concept of \emph{exponential map} $\exp_{O} : T_{O} S^2 (a^2) \rightarrow  S^2 (a^2)$.

\begin{defn}\label{expsphereDef}
\textbf{The exponential map $\exp_O$ about the reference point $O\in S^2(a^2)$:}
For any vector $v \in T_{O} S^2 (a^2)$, $\exp_{O}(v)$ is defined as
\begin{equation}\label{exp}
\exp_{O}(v) = \gamma_{v}(1),
\end{equation}
where $\gamma_{v} : [0,\infty ) \rightarrow S^2 (a^2)$ is the unique geodesic on $S^2 (a^2)$ which satisfies $\gamma (0) = O$ and $\frac{d \gamma }{dt}|_{t=0} = v$.
\end{defn}

\noindent
Indeed, in accordance with \eqref{exp}, it is plain to see that the following relation holds for any vector $v \in T_{O} S^2 (a^2)$ with $|v| = 1$ and any $t > 0$.

\begin{equation}
\exp_{O}(t v) = \gamma_{v}(t).
\end{equation}
It is a well known fact that the exponential map $\exp_{O}$, once restricted on the open disc
$D_{0} (\frac{\pi }{a}) = \{ v \in T_{O} S^2 (a^2) : |v| < \frac{\pi}{a}   \} $, becomes a \emph{diffeomorphism} of $D_{0} (\frac{\pi }{a})$ onto $S^{2}(a^{2})-\{(0,0,-\frac{1}{a})\}$. This means that  $\exp_{O} : D_{0} (\frac{\pi }{a}) \rightarrow S^{2}(a^{2})-\{(0,0,-\frac{1}{a})\}$ is a smooth bijective map with a smooth inverse map $\exp_{O}^{-1} : S^{2}(a^{2})-\{(0,0,-\frac{1}{a})\} \rightarrow  D_{0} (\frac{\pi }{a}) $. Since we have the natural identification of the space form $S^{2}(a^{2})$ with $\{ (x_1 ,x_2 , x_3 ) : x_{1}^2 + x_{2}^2 + x_{3}^2 = \frac{1}{a^2}   \}$, $\exp_{O} : T_{O} S^2 (a^2) \rightarrow  S^2 (a^2)$ can be given explicitly as follow. Since $T_{O} S^2 (a^2)$ is naturally identified with $\mathbb{R}^{2}$, any vector $v \in T_{O} S^2 (a^2)$ with $|v| = 1$ can be represented as $v = (\cos \lambda , \sin \lambda )$, for some $\lambda \in [0 , 2\pi )$. Then, the unique geodesic $\gamma_{v} : [0,\infty )  \rightarrow S^2 (a^2)$ satisfying $\gamma_{v}(0) = O$, and $\frac{d \gamma_{v}}{dt}|_{t=0} = v$ is given explicitly by

\begin{equation}
\gamma_{v}(t) = \frac{1}{a} (\sin (ta)  \cos \lambda , \sin (ta) \sin \lambda , \cos (ta)).
\end{equation}
Hence, $\exp_{O} : D_{0} (\frac{\pi }{a}) \rightarrow  S^{2}(a^{2})-\{(0,0,-\frac{1}{a})\}$ can be given explicitly as follow.
\begin{equation}
\exp_{O} (t v) = \frac{1}{a} (\sin (t a)  \cos \lambda , \sin (ta) \sin \lambda , \cos (ta)) ,
\end{equation}
where $v = (\cos \lambda , \sin \lambda )$ is a unit vector in $T_{O} S^2 (a^2)$, and $t \in [0 , \frac{\pi}{a})$.\\

\noindent
We can now introduce the \emph{normal polar coordinate system} on $S^2 (a^2)$, through the use of the exponential map $\exp_O$ as follow.

\begin{defn}\label{PolarDefnSphere}
\textbf{The normal polar coordinate system on the sphere $S^2(a^2)$ :} The normal polar coordinate system on $S^2 (a^2)$ about the based point $O \in S^2 (a^2)$ is the bijective smooth map $( r,\theta ) : S^2 (a^2)-\{(x_1 , 0 , x_3) : x_{1}^2 + x_{3}^2  = \frac{1}{a^2}, x_1 \geqslant 0 \} \rightarrow
(0 , \frac{\pi}{a})\times (0 , 2\pi )$ defined by
\begin{equation}
(r, \theta )  = (\overline{r} , \overline{\theta}) \circ (exp_{O}^{-1}) ,
\end{equation}
where $(\overline{r} , \overline{\theta})$ stands for the \emph{standard polar coordinate system} $(\overline{r} , \overline{\theta }) : \mathbb{R}^{2}-\{(x_1 , 0) : x_{1} \geqslant 0 \} \rightarrow (0, \infty )\times (0 , 2\pi)$ on $\mathbb{R}^{2}$. In order words, the normal polar coordinate system $(r,\theta )$ on $S^2(a^2)$ is \emph{defined} to be the composite of the standard polar coordinate system $(\overline{r} , \overline{\theta})$ on $\mathbb{R}^2$ with
$\exp_{O}^{-1} : S^2 (a^2)-\{(x_1 , 0 , x_3) : x_{1}^2 + x_{3}^2  = \frac{1}{a^2}, x_1 \geqslant 0 \} \rightarrow D_{0}(\frac{\pi}{a}) - \{ (x_{1} , 0) : 0 \leqslant x_1 < \frac{\pi}{a}  \}$.
\end{defn}

\noindent
With such a normal polar coordinate system $(r , \theta )$ on $S^2(a^2)$, we now \emph{define} two \emph{natural, everywhere linearly independent vector fields} $\frac{\partial}{\partial r}$, and $\frac{\partial}{\partial \theta}$ on $S^2(a^2)$ as follow.

\begin{defn}\label{coordFRAMEsphere}
\textbf{Natural coordinate frame $\big \{ \frac{\partial}{\partial r} , \frac{\partial}{\partial \theta}   \big \}$ on the sphere $S^2(a^2)$ : }For any point $p \in S^2 (a^2)-\{(x_1 , 0 , x_3) : x_{1}^2 + x_{3}^2  = \frac{1}{a^2}, x_1 \geqslant 0 \}$, the vectors $\frac{\partial}{\partial r}|_{p}$,
$\frac{\partial}{\partial \theta}|_{p}$ in the tangent space $T_{p}(S^2(a^2))$ can be regarded as \emph{linear functionals} on the space $C^{\infty}(S^{2}(a^{2}))$ of  smooth functions on $S^2 (a^2)$, and hence, $\frac{\partial}{\partial r}|_{p}$, $\frac{\partial}{\partial \theta}|_{p}$ in $T_{p}S^{2}(a^2)$ can be \emph{defined} through the following characterization, where $f \in C^{\infty}(S^{2}(a^2))$.
\begin{equation}\label{naturalvector}
\begin{split}
\frac{\partial}{\partial r}\bigg |_{p} f &= \frac{\partial}{\partial \overline{r}}\bigg |_{(r,\theta )(p)}[f \circ (r,\theta)^{-1}] =
     \frac{\partial}{\partial \overline{r}}\bigg |_{(r,\theta )(p)}[f \circ \exp_{O} \circ (\overline{r} , \overline{\theta})^{-1} ] \\
\frac{\partial}{\partial \theta }\bigg |_{p} f & =  \frac{\partial}{\partial \overline{\theta}}\bigg |_{(r,\theta )(p)}[f \circ (r,\theta)^{-1}] =
\frac{\partial}{\partial \overline{\theta}}\bigg |_{(r,\theta )(p)}[f \circ \exp_{O} \circ (\overline{r} , \overline{\theta})^{-1} ].
\end{split}
\end{equation}
\end{defn}

\noindent
The two vector fields $\frac{\partial}{\partial r}$, and $\frac{\partial}{\partial \theta}$ on $S^2(a^2)$, which are characterized by relation \eqref{naturalvector}, are everywhere linearly independent on $S^2(a^2)$, and they together constitutes \emph{the natural moving frame} induced by the normal polar coordinate system on $S^{2}(a^{2})$. Now, by means of the identification $S^{2}(a^2) = \{(x_1 , x_2 , x_3) \in \mathbb{R}^3 : x_1^2 + x_2^2 + x_3^2 = \frac{1}{a^2}  \}$, we can express the vector fields $\frac{\partial}{\partial r}$, and $\frac{\partial}{\partial \theta}$ concretely as follow. \\

\noindent
Since $\exp_{O} : \{v \in T_{O}S^2(a^2) : |v| < \frac{\pi}{a}\} \rightarrow S^2(a^{2}) -\{(0,0,-\frac{1}{a})\}$ is a diffeomorphism, any given $p \in S^{2} - \{(0,0, -\frac{1}{a})\}$ which is away from the base point $O$ can be represented as $p = \exp_{O}(rv)$, with some uniquely determined unit vector $v = (\cos \lambda , \sin \lambda ) \in T_{O}S^2(a^2)$, and $r \in (0, \frac{\pi}{a})$. It turns out that $\frac{\partial}{\partial r}|_p$, and $\frac{\partial}{\partial \theta}|_p$ can be expressed concretely as follow, with $\gamma_{v} ;[0, \infty ) \rightarrow S^2(a^2)$ to be the unique geodesic satisfying $\gamma_{v}(O) =0$ and $\frac{d \gamma_{v}}{dt}|_{t=0} =v$, and
$v^{\perp} = (-\sin \lambda , \cos \lambda )$.
\begin{equation}\label{vectorexplicit}
\begin{split}
\frac{\partial}{\partial r}\bigg |_p &= \frac{d \gamma_{v}}{dt} \bigg |_{t=r} = \cos (ar) (v,0) - \sin (at) (0,0,1) \in T_{p}S^2(a^2) \\
\frac{\partial}{\partial \theta }\bigg |_p &= \frac{1}{a}\sin (ar) (v^{\perp} , 0) = \frac{1}{a} \sin (ar) (-\sin \lambda , \cos \lambda ,0 ) \in T_{p}S^2(a^2).
\end{split}
\end{equation}
 From \eqref{vectorexplicit}, it follows at once that the vector fields $e_{1} = \frac{\partial}{\partial r}$, and $e_{2} = \frac{a}{\sin (ar)} \frac{\partial}{\partial \theta }$ together constitute a moving frame on $S^{2}(a^2)$ which is \emph{everywhere orthonormal} on $S^{2}(a^2)$. That is, we know that $|e_{1}(p)|=1$, $|e_{2}(p)|=1$, and $g(e_{1}(p), e_{2}(p)) = 0$ holds for all $p \in S^2(a^2)$. Here, of course the symbol $g(\cdot , \cdot )$ denotes the Riemannian metric on $S^2(a^2)$, with $|v| = [g(v,v)]^{\frac{1}{2}}$. We also observe that the vector field $e_{2}$ satisfies the following relation, where $v \in T_{O}S^2(a^2) $ is a unit vector, and $0 < r < \frac{\pi}{a}$.

\begin{equation}
e_{2} |_{\exp_{O}(rv)} = (v^{\perp} , 0),
\end{equation}
which indicates that the vector field $e_{2}$, once being restricted to the geodesic ray $\gamma_{v}$ starting from $O$, is a \emph{parallel along} $\gamma_{v}$. So, it follows that $\overline{\nabla}_{\frac{d \gamma }{dt}} (e_{2}) = 0$, with
$\overline{\nabla } : C^{\infty }(TS^{2}(a^2)) \rightarrow C^{\infty} (T^{*}S^2(a^2)\otimes TS^2(a^2)) $ to be the Levi-Civita connection induced by the Riemannian metric $g (\cdot , \cdot )$ of $S^{2}(a^2)$ (See \textbf{Definition \ref{LeviCivitadefinition}} in \textbf{Section \ref{Spheremainsection}} for a precise definition of the Levi-Civita connection as induced by the intrinsic geometry of a Riemannian manifold). Upon this setting, we can now discuss the Hodge-star operator and the Hodge Laplacian on $S^2(a^2)$ in the following subsection.\\

\noindent{\bf The Hodge-star operator and Hodge Laplacian on $S^2(a^2)$ in terms of the normal polar coordinate on $S^2(a^2)$.}

\noindent
The Hodge Laplacian $(-\triangle ) = d d^* + d^* d$ on a Riemannian manifold $M$ is an operator sending the space $C^{\infty}(T^{*}M)$ of all smooth $1$-forms into  $C^{\infty}(T^*M)$ itself. So, the Hodge Laplacian $(-\triangle )$ on $S^2(a^2)$ acts on differential $1$-forms \emph{instead of vector fields} on $S^2(a^2)$. This leads (\emph{actually forces}) us to consider the natural identification of the space $C^{\infty}(TS^2(a^2))$ of smooth vector fields on $S^2(a^2)$ with the space
$C^{\infty}(T^*S^2(a^2))$ of smooth $1$-forms, which identifies a smooth vector field $v \in C^{\infty}(TS^2(a^2))$ with the associated $1$ form $v^* = g(v, \cdot ) \in C^{\infty}(T^*M)$. Here, $g(\cdot , \cdot )$ is the Riemannian metric of the sphere $S^2(a^2)$.\\

\noindent
Based upon the polar normal coordinate system $(r, \theta )$ on $S^2(a^2)$ as given in \textbf{Definition \ref{PolarDefnSphere}}, we have two smooth vector fields $e_{1} = \frac{\partial}{\partial r}$, and $e_{2} = \frac{a}{\sin (ar)} \frac{\partial}{\partial \theta }$, which together constitute a positively oriented orthonormal moving frame on\\
 $S^2(a^2)-\{(0,0, \frac{1}{a}) , (0,0,\frac{-1}{a})\}$ for the tangent bundle $TS^2(a^2)$ on $S^2(a^2)$. Then, the associated $1$-forms $e_1^* = g (e_1, \cdot )$, and $e_2^* = g(e_2 , \cdot )$ together constitute a positively oriented orthonormal moving \emph{co-frame} on for the cotangent bundle $T^{*}S^2(a^2)$ on $S^2(a^2)$. Indeed, the $1$ forms $e_{1}^*$, and $e_2^*$ can be expressed by
\begin{equation}\label{dualcoframe}
\begin{split}
e_{1}^{*} &= dr ,\\
e_{2}^{*} & = \frac{\sin ar}{a} d\theta .
\end{split}
\end{equation}
Here, we first define the volume form on a general oriented $2$-dimensional Riemannian manifold $M$, which is a everywhere non-vanishing globally defined $2$-form on $M$ induced by the intrinsic Riemannian geometry of $M$.

\begin{defn}\label{VolumeformDef}
\textbf{Volume form on a $2$-dimensional oriented Riemannian manifold :}
Given $M$ to be an oriented $2$-dimensional Riemannian manifold equipped with a Riemannian metric $g(\cdot , \cdot )$. Consider two locally defined vector fields $e_{1}$, $e_2$ on some open region $U$ of $M$ which together constitute an positively oriented orthonormal moving frame $\{e_1 , e_2\}$ of the tangent bundle $TM$ of $M$ over $U$. Then, the (locally defined) associated $1$-forms $e_{1}^* = g(e_1 ,\cdot )$, $e_2^* = g(e_2 , \cdot )$ will constitute the so-called orthonormal co-frame for the cotangent bundle $T^*M$ over $M$. Then, the volume form $Vol_{M}$ can be locally defined through the relation
\begin{equation}
Vol_M = e_1^* \wedge e_2^*.
\end{equation}
It is an elementary fact in differential geometry that such a local definition of $Vol_{M}$ turns out to be \textbf{independent} of the choice of the locally defined positively oriented orthonormal frame $\{e_1 , e_2\}$. Hence, our local construction gives a globally-defined Volume form $Vol_M$ on the whole $2$-dimensional manifold $M$ (see Chapter 2 of \cite{Jost} for a more general discussion).
\end{defn}

\noindent
Hence, the \emph{globally defined} volume form $Vol_{S^2(a^2)}$ on $S^{2}(a^2)$ is the $2$-form which can \emph{locally be expressed} by $Vol_{S^{2}(a^2)} = e_{1}^{*}\wedge e_{2}^* = \frac{\sin ar}{a} dr\wedge d\theta \in C^{\infty}(\wedge^2 T^*S^2(a^2))$. Here, the symbol $C^{\infty}(\wedge^2 T^*S^2(a^2))$ stands for the space of all smooth $2$-forms on $T^*S^2(a^2)$.\\

\noindent
Before we talk about the Hodge-star operator on $S^{2}(a^2)$, we recall the general definition of the Hodge-Star operator on a given $N$-dimensional manifold (See Chapter $2$ of Jost \cite{Jost}).
\begin{defn}\label{HodgeDefinition}
\textbf{The Hodge-Star operator :}
In the case of a general oriented Riemannian manifold $M$ of dimension $N$,for each integer $0 \leqslant k \leqslant N$,  the Hodge-star operator $* : C^{\infty }(\wedge^k T^*M ) \rightarrow C^{\infty }(\wedge^{N-k} T^*M ) $ sending the space of $k$ forms on $M$ to the space of $N-k$ forms on $M$ is characterized by the following relation: For any $k$-forms $\alpha$, $\beta$, where $\overline{g} (\cdot , \cdot )$ stands for the metric on the vector bundle $\wedge^{k}T^*M$ induced by the Remannian metric $g(\cdot , \cdot )$ on $M$, and that $Vol_{M}$ is the volume form on $M$ as defined in \textbf{Definition \ref{VolumeformDef}}.
\begin{equation}
\alpha \wedge *\beta = \overline{g}(\alpha , \beta ) Vol_{M}.
\end{equation}
\end{defn}
\noindent
Since the dimension of $S^2(a^2)$ is just $2$, in the case of $M = S^{2}(a^2)$ the Hodge-star operator $*$ can easily be described as follow.
First, $* : C^{\infty}(T^{*}S^{2}(a^2)) \rightarrow C^{\infty}(T^{*}S^{2}(a^2)) $ is characterized by the following relation
\begin{equation}\label{Hodgestar}
\begin{split}
*e_{1}^* &= e_{2}^* \\
*e_{2}^* & = -e_{1}^{*}.
\end{split}
\end{equation}
Since $* : C^{\infty }(\wedge^k T^*M ) \rightarrow C^{\infty }(\wedge^{N-k} T^*M ) $ is \emph{tensorial} in that the relation
$*(f \alpha ) = f *(\alpha)$ holds for any smooth function $f \in C^{\infty}(M)$ and $k$-form $\alpha$ on a Riemannian manifold $M$. It follows that in the case of
$M = S^{2}(a^2)$, we have
\begin{equation}
\begin{split}
* dr &= e_2^* = \frac{\sin ar}{a} d\theta ,\\
* d\theta &= *[\frac{a}{\sin ar} e_{2}^*] = -\frac{a}{\sin ar} e_{1}^* = -\frac{a}{\sin ar} dr.
\end{split}
\end{equation}
In the structure of the Hodge Laplacian sending $C^{\infty}(T^*S^2(a^2))$ into itself, one also encounters the Hodge-star operator
$* : C^{\infty}(\wedge^{2} T^* S^2(a^2)) \rightarrow C^{\infty}(S^2(a^2))$ sending smooth $2$-forms to smooth functions on $S^2(a^2)$, which can be characterized via the following relation.
\begin{equation}
*Vol_{S^2(a^2)} = * (e_1^* \wedge e_2^*) = 1,
\end{equation}
where we note that the $2$-form $e_1^* \wedge e_2^*$ is the volume form $Vol_{S^2(a^2)}$ on $S^2(a^2)$. By using the tensorial property of $*$, it follows that $* (dr \wedge d\theta ) = \frac{a}{\sin ar}$.
Now, we can discuss the Hodge Laplacian $(-\triangle) = d d^* + d^* d$, where the exterior differential operators $d : C^{\infty}(\wedge^kT^*S^2(a^2)) \rightarrow C^{\infty}(\wedge^{k+1}T^*S^2(a^2))$ and their associated adjoint operators $d^* : C^{\infty}(\wedge^{k+1}T^*S^2(a^2)) \rightarrow  C^{\infty}(\wedge^kT^*S^2(a^2))$ are involved for $k =0,1$. The operator $d : C^{\infty}(S^2(a^2)) \rightarrow C^{\infty}(T^*S^2(a^2))$ sends a smooth function $f$ to the $1$-form $df$, which can locally be expressed, in terms of the natural co-frame $\{dr , d \theta \}$, as follow.
\begin{equation}
df = \frac{\partial f}{\partial r} dr + \frac{\partial f}{\partial \theta }  d \theta ,
\end{equation}
where the smooth functions $\frac{\partial f}{\partial r}$ and $\frac{\partial f}{\partial \theta } $ are defined via \eqref{naturalvector}. It is an elementary fact that $df$ is \emph{independent} of the choice of the coordinate frame being used in its characterization. Here, we express $df$ in terms of the normal polar coordinate system $(r, \theta )$, since this will give us a quick and easy way in computing $(-\triangle) u^*$, with $u^* = g (u, \cdot )$ to be the associated $1$-form of the vector field $u$ representing a parallel laminar flow near some arc-shaped boundary of some obstacle in the sphere $S^2(a^2)$.

\begin{remark}\label{remarkgradient}
On a general $N$-dimensional Riemannian manifold $M$, one should \textbf{interpret} the $1$-form $df$ as the gradient field $\nabla f$ in the sense that the gradient of $f$ is the unique vector field $\nabla f$ on $M$ characterized by the following relation, with $g(\cdot , \cdot )$ to be the Riemannian metric on $M$.
\begin{equation}
df = g(\nabla f , \cdot ).
\end{equation}
\end{remark}

\noindent
Next, one may characterize the operator $d : C^{\infty}(T^* S^2(a^2) ) \rightarrow C^{\infty}(\wedge^2 T^*S^2(a^2)) $ sending smooth $1$-forms to $2$-forms in terms of the normal polar coordinate $(r, \theta )$ as follow. Since every smooth $1$-form $\alpha$ on $S^2(a^2)$ can locally be expressed as
\begin{equation}
\alpha = \alpha_{r} dr + \alpha_{\theta} d\theta ,
\end{equation}
with $\alpha_{r}$, and $\alpha_{\theta}$ to be some locally defined smooth functions on $S^2(a^2)$. Then, the $2$-form $d\alpha$ is locally expressed by

\begin{equation}
\begin{split}
d \alpha &= d \alpha_{r} \wedge dr + d \alpha_{\theta} \wedge d\theta \\
& = \frac{\partial  \alpha_{r}}{\partial \theta } d\theta \wedge dr + \frac{\partial \alpha_{\theta}}{\partial r} dr \wedge d\theta \\
& = \big \{ \frac{\partial \alpha_{\theta}}{\partial r} - \frac{\partial  \alpha_{r}}{\partial \theta }    \big \}  dr \wedge d\theta \\
& = \frac{a}{\sin (ar)} \big \{ \frac{\partial \alpha_{\theta}}{\partial r} - \frac{\partial  \alpha_{r}}{\partial \theta }    \big \} Vol_{S^2(a^2)} .
\end{split}
\end{equation}
where in the above computation, we have implicity used the facts that $dr \wedge dr = 0$, $d\theta \wedge d\theta = 0$, and $d\theta \wedge dr = -dr \wedge d\theta$, each of these follows directly from the definition $\alpha \wedge \beta = \frac{1}{2} \{\alpha \otimes \beta - \beta \otimes \alpha \}$, for two given smooth $1$-forms $\alpha$, and $\beta$ on a smooth manifold $M$.\\

\begin{remark}\label{curlremark}
We should think of the operator $d$ sending the space of $1$-forms into the space of $2$-forms on an oriented $2$-dimensional Riemannian manifold $M$ to be the \textbf{curl-operator} which assigns sends each smooth vector field $u$ to its vorticity function $\omega$ on $M$. More precisely, consider $u$ to be a smooth vector field on an oriented $2$-dimensional Riemannian manifold $M$ with Riemannian metric $g(\cdot, \cdot )$. Then, by taking the operator $d$ on the associated $1$-form $u^* = g(u,\cdot )$ of the vector field $u$, we yield $du^* = \omega Vol_{M}$, with some uniquely determined smooth function $\omega$ on $M$, which is exactly the vorticity of the vector field $u$. So, we should really think of the operator $d$ sending $1$-forms into $2$-forms as the natural generation of the \textbf{curl-operator} in the more general sense.
\end{remark}

\noindent
We can now easily characterize the co-adjoint operators $d^* : C^{\infty}(\wedge^p T^*S^2(a^2)) \rightarrow C^{\infty}(\wedge^{p-1} T^*S^2(a^2))$, with $p = 1, 2$, by means of the following standard definition in differential geometry (see for instance, Chapter 2 of \cite{Jost}).

\begin{defn}\label{dstarDef}
\textbf{The co-adjoint operators $d^*$ :}
Given $M$ to be a $2$-dimensional oriented manifold equipped with a Riemannian metric $g(\cdot , \cdot )$. For each $p= 1, 2$, the coadjoint operator  $d^* : C^{\infty}(\wedge^p T^*M) \rightarrow C^{\infty}(\wedge^{p-1} T^*M)$ sending the space of $p$-forms on $M$ into the space of $p-1$-forms on $M$ is defined through the following relation.
\begin{equation}
d^{*} = (-1)^{2(p + 1) + 1}* d * = -*d* ,
\end{equation}
where the symbol $*$ stands for the Hodge-Star operators as defined in \textbf{Definition \ref{HodgeDefinition}}.
\end{defn}

\noindent
Two remarks about \textbf{Definition \ref{dstarDef}} are in order here.

\begin{remark}
We remark that, in the above formula for $d^*$, we have included the extra index $2(p +1)$ in the power of $(-1)$, simply because for a general $N$-dimensional Riemannian manifold $M$, the formula for $d^*$ acting on smooth $p$-forms is exactly $d^* = (-1)^{N(p+1)+1}*d*$.
\end{remark}

\begin{remark}\label{divergenceremark}
On an oriented $2$-dimensional Riemannian manifold $M$ equipped with Riemannian metric $g(\cdot, \cdot )$, the co-adjoint operator $d^* : C^{\infty}(T^*M) \rightarrow  C^{\infty}(M)$, which sends $1$-forms into smooth functions, should be interpreted as the \textbf{divergence operator $-div$} acting on smooth vector fields on $M$ in the following sense: If $u$ is a smooth vector field on $M$, with associated $1$-form $u^* = g(u ,\cdot )$. Then, it is a standard fact in Riemannian geometry that the following relation holds for any smooth test function $f \in C^{\infty}_c(M)$ (see Chapter $2$ of \cite{Jost}, for instance).
\begin{equation}\label{integrationbypartI}
\int_{M} g(u , \nabla f ) Vol_{M} =  \int_{M} f  d^*u^* Vol_{M} .
\end{equation}
However, we know that $-div (u)$ also satisfies the following integration by parts formula, for any test function $f \in C^{\infty}_c(M)$,
\begin{equation}\label{integrationbyparttwo}
\int_{M} g(u , \nabla f ) Vol_{M} = - \int_{M} f  div(u) Vol_{M} .
\end{equation}
So, by comparing \eqref{integrationbypartI} with \eqref{integrationbyparttwo}, we are forced to conclude that $d^*u^*= -div (u)$.
\end{remark}

\section{Proof of Theorem \ref{Noparallel}, and the proof of Theorem \ref{existenceEasy} }\label{Spheremainsection}

\noindent
To begin the proof of Theorem \ref{Noparallel}, let us consider the space form $S^2(a^2)$ of constant sectional curvature $a^2 > 0$. Let $O \in S^2(a^2)$ to be an selected reference point on $S^2(a^2)$, and let $(r,\theta )$ to be the normal polar coordinate on $S^2(a^2)$ about the base point $O$, which we introduce in the previous section.
Here, let $K$ to be some compact region in $S^2(a^2)$, which is contained in the closed geodesic ball
$\overline{B_{O}(\delta )} = \{p\in S^2(a^2) : d(p,O) \leq \delta \}$ for some positive radius $0 < \delta < \frac{\pi}{a}$.
As in the hypothesis of Theorem \ref{Noparallel}, we assume that $\partial K$ contains a circular-arc portion in that, for some positive angle
$\tau \in (0, 2\pi )$ the circular arc $C_{\delta , \tau }= \{p \in S^2(a^2) : r(p)= \delta ,   0< \theta (p)< \tau  \}$ is contained in $\partial K$.\\

\noindent
Now, let $u$ to be some smooth vector field defined over the simply-connected open region $R_{\delta , \tau , \epsilon_0}$ as specified in \eqref{exteriorregion} of \textbf{Section \ref{IntroductionSEC}}, with $\epsilon_0 > 0$ to be some positive number. Again, remember that the exterior region $R_{\delta , \tau , \epsilon_0}$ shares the same circular-arc boundary portion $C_{\delta , \tau}$ with $K$. Suppose that $u$ is a \emph{parallel laminar flow} on $R_{\delta , \tau , \epsilon_0}$ in the sense that $u$ can be expressed as
\begin{equation}\label{circularStreamlines}
u  = - h( r -\delta ) e_{2} = - h( r -\delta ) \frac{a}{\sin ar} \frac{\partial}{\partial \theta }.
\end{equation}
In the expression \eqref{circularStreamlines}, $r$, which is the first component of the normal polar coordinate $(r, \theta )$, measures the distance of a point $p \in S^{2}(a^2)$ from the base point $O$. That is, $r(p) = d(p,O)$, with $d(p,O)$ to be the Riemannian distance of $p$ from $O$ in the Riemannian manifold $S^2(a^2)$. Also, in \eqref{circularStreamlines}, $e_{2} = \frac{a}{\sin (ar)} \frac{\partial}{\partial \theta }$ is the second vector field in the orthonormal moving frame
$e_{1}= \frac{\partial}{\partial r}$, $e_{2} = \frac{a}{\sin (ar)} \frac{\partial}{\partial \theta }$. Recall that from \eqref{dualcoframe}, the natural dual coframe of the moving frame $\{e_{1}, e_{2}\}$ is given by $e_{1}^* = dr$, and $e_{2}^* = \frac{\sin (ar)}{a} d\theta$. Hence the volume form on $S^2(a^2)$ is $Vol_{s^2(a^2)} = e_{1}^*\wedge e_{2}^* = \frac{\sin ar}{a} dr \wedge d\theta$.\\

\noindent
Against such a setting, the first step which we take is to compute $(-\triangle )u^*$, with $u^* = g(u, \cdot )$ to be the associated $1$-form of the parallel laminar flow $u$.\\

\noindent{\bf Step 1 : The divergence free property $d^* u^* = 0$ and the computation of $(-\triangle )u^*$, for $u$ as given in \eqref{circularStreamlines}}

\noindent
Here, it follows from \eqref{circularStreamlines} that the associated $1$-form $u^* = g (u,  \cdot )$ is locally given by
\begin{equation}\label{associatedform}
u^* = -h(r-\delta ) e_{2}^* = - h(r-\delta )\frac{\sin (ar)}{a} d \theta .
\end{equation}
First, we point out that the \emph{divergence free condition} $d^*u^* = 0$ is just a direct consequence which follows from the following computation.
\begin{equation}\label{divfreesphere}
\begin{split}
d^* u^* &= (-1)^{2(1 +1)  +1} *d* \{ -h(r-\delta ) \frac{\sin (ar) }{a} d\theta  \} \\
& = *d*[h(r-\delta ) e_{2}^{*}] \\
& = -* d[h(r-\delta ) dr ]\\
& = -* \{ d[h(r-\delta )] \wedge dr   \}\\
& = -* \{ [\frac{\partial h(r - \delta  )}{\partial r} dr +  \frac{\partial h(r - \delta  )}{\partial \theta } d\theta] \wedge dr                 \} \\
& = -* \{  [\frac{\partial h(r - \delta  )}{\partial r} dr \wedge dr  \}\\
&=0.
\end{split}
\end{equation}
In the above calculation, the third equality follows from \eqref{Hodgestar}, the sixth equal sign follows from the fact that $\frac{\partial h(r - \delta  )}{\partial \theta } = 0$, and the last equality is due to the fact that $dr \wedge dr = 0$. Here, $d^* u^* = 0$ means that the vector field $u$ is divergence free, since one thinks of $d^*$ as $-\dv$ when the language of smooth $1$-form is \emph{translated} back to the language of smooth vector field via the correspondence $u^* = g(u , \cdot )$, just as one thinks of the 1-form $df$ as the gradient field $\nabla f$ in the language of vector fields.\\

\noindent
Now, the divergence free property $d^* u^* = 0$ implies that $(-\triangle )u^* = (dd^* + d^*d)u^* = d^*du^*$, since the term $dd^*u^*$ vanishes.
Then we have
\begin{equation}
\begin{split}
du^{*} &= -\frac{\partial}{\partial r} [h(r-\delta )  \frac{\sin (ar)}{a}] dr \wedge d\theta \\
& = -\frac{\partial}{\partial r} [h(r-\delta )  \frac{\sin (ar)}{a}] \cdot \frac{a}{\sin (ar)} Vol_{S^2(a^2)} .
\end{split}
\end{equation}
As a result, we have the following direct computation, in accordance with the definition of the operators $d$, and $d^*$ as discussed in the previous section.
\begin{equation}
\begin{split}
d^{*}du^* &= (-1)^{2(2 + 1 )  +1}*d*\{  -\frac{\partial}{\partial r} [h(r-\delta )  \frac{\sin (ar)}{a}]\cdot \frac{a}{\sin (ar)} Vol_{S^2(a^2)}  \}\\
& = * d \{  \frac{\partial}{\partial r} [h(r-\delta )  \frac{\sin (ar)}{a}]\cdot  \frac{a}{\sin (ar)}     \}\\
&= *  \frac{\partial}{\partial r} \{   \frac{\partial}{\partial r} [h(r-\delta )  \frac{\sin (ar)}{a}]\cdot  \frac{a}{\sin (ar)}       \}  dr \\
& =   \frac{\partial}{\partial r} \{   \frac{\partial}{\partial r} [h(r-\delta )  \frac{\sin (ar)}{a}]\cdot  \frac{a}{\sin (ar)}       \} *e_{1}^* \\
&=    \frac{\partial}{\partial r} \{   \frac{\partial}{\partial r} [h(r-\delta )  \frac{\sin (ar)}{a}]\cdot  \frac{a}{\sin (ar)}       \} e_{2}^* \\
& =  \frac{\partial}{\partial r} \{   \frac{\partial h(r-\delta )}{\partial r} + a h(r-\delta ) \frac{\cos (ar)}{\sin (ar)}    \} e_{2}^{*}.
\end{split}
\end{equation}
That is, we have the following local expression for $(-\triangle)u^{*}$
\begin{equation}\label{triangleu}
\begin{split}
(-\triangle )u^* & = \frac{\partial}{\partial r} \{   \frac{\partial h(r-\delta )}{\partial r} + a h(r-\delta ) \frac{\cos (ar)}{\sin (ar)}    \} e_{2}^{*}\\
& = \frac{\partial}{\partial r} \{   \frac{\partial h(r-\delta )}{\partial r} + a h(r-\delta ) \frac{\cos (ar)}{\sin (ar)}    \} \cdot \frac{\sin (ar)}{a} d\theta \\
& = \{ h''(r-\delta )  \frac{\sin (ar)}{a} + h'(r-\delta ) \cos (ar) - \frac{a}{\sin (ar)} h(r-\delta )   \} d\theta .
\end{split}
\end{equation}
Next, we need to compute the convection term $[\overline{\nabla }_{u}u]^{*}$, with $\overline{\nabla} : C^{\infty}(TS^2(a^2)) \rightarrow C^{\infty} (T^*S^2(a^2) \otimes TS^2(a^2) )$ to be the Levi-Civita connection (covariant derivative) acting on the space of smooth vector fields on $S^2(a^2)$
(See \textbf{Definition \ref{LeviCivitadefinition}} for the precise meaning of $\overline{\nabla}$).\\


\noindent{\bf Step 2 : The computation of $\overline{\nabla}_{u}u^*$, for $u$ as given in \eqref{circularStreamlines}.}

\noindent
For a $N$-dimensional Riemannian manifold $M$ equipped with a Riemannian metric $g (\cdot , \cdot )$ in general, the Levi-Civita connection
$\overline{\nabla} : C^{\infty}(TM) \rightarrow C^{\infty}(T^*M\otimes TM)$ is uniquely characterized as follows.

\begin{defn}\label{LeviCivitadefinition}
The Levi Civita connection $\overline{\nabla} : C^{\infty}(TM) \rightarrow C^{\infty}(T^*M\otimes TM)$ acting on the space of smooth vector fields on a $N$-dimensional Riemannian manifold $M$ $($ equipped with a Riemannian metric $g(\cdot , \cdot )$ $)$ is uniquely determined by the following characterizing properties:

\begin{itemize}
\item \textbf{(1)} $($compatibility condition with the Riemannian metric $g(\cdot , \cdot )$ on $M$$)$ The relation $X(g(Y,Z)) = g(\overline{\nabla}_{X}Y , Z) + g(Y, \overline{\nabla}_{X}Z)$, holds for any smooth vector fields $X$, $Y$, $Z$ on $M$.
(Here, the notation $X(g(Y,Z))$ stands for the directional derivative of the function $g(Y,Z)$ along the direction of $X$. See also \textbf{Remark \ref{LastRemark}}.)
 \item \textbf{(2)} $($torsion free property of the connection$)$ $\overline{\nabla}$ is \emph{torsion free} in that the relation $\overline{\nabla}_{X}Y - \overline{\nabla}_{Y}X - [X,Y] = 0$, holds for any smooth vector fields $X$, $Y$ on $M$. Here, $[X,Y] = XY-YX$ is the Lie Bracket of $X$ and $Y$.
\item \textbf{(3)} $($Leibniz rule of $\overline{\nabla}$ as a covariant derivative$)$ For any smooth function $f$ on $M$, and any smooth vector fields $X$, $Y$ on $M$, the relation $\overline{\nabla}_X\big ( f X\big ) = X\big ( f \big )\cdot Y + f \overline{\nabla}_XY $ holds $($here, the symbol $X\big ( f \big )$ stands for the derivative of $f$ along the direction of the vector field $X$$)$.
\item \textbf{(4)} $($Tensorial property of $\overline{\nabla}$$)$ For any smooth function $f$ on $M$, and any smooth vector fields $X$, $Y$ on $M$, we always have $\overline{\nabla}_{fX} Y = f \overline{\nabla}_{X}Y $.
\end{itemize}
\end{defn}

\begin{remark}\label{LastRemark}
In conditions \textbf{(1)}, and \textbf{(3)} as given in \textbf{Definition \ref{LeviCivitadefinition}}, we have employed, for a given smooth vector field $X$ and a smooth function $f$ on $M$, the notion $Xf$, which is the \emph{rate of change of $f$ in the direction of $X$}. More precisely, for each $p\in M$, $Xf|_{p}$ is \emph{defined as} $Xf|_{p} = \frac{d}{dt}(f\circ \gamma )|_{t=0}$, with $\gamma : [0, \epsilon ) \rightarrow M$ to be some smooth path with $\gamma (0) = p \in M$, and $\frac{d\gamma }{dt}|_{t=0} = X_{p}$.
\end{remark}

\noindent
Now, in the case of $M = S^2(a^2)$, we again consider the locally defined parallel laminar flow $u = -h(r-\delta )e_{2}$ near some circular arc portion of $\partial \Omega$, as defined in \eqref{circularStreamlines}, and compute $\overline{\nabla}_{u}u$, with $\overline{\nabla} : C^{\infty}(TS^2(a^2)) \rightarrow C^{\infty}(T^*S^2(a^2) \otimes TS^2(a^2))$ to be the Levi-Civita connection with respect to the standard metric $g(\cdot , \cdot )$ on $S^2(a^2)$. As a first step in computing  $\overline{\nabla}_{u}u$, we set
\begin{equation}
\overline{\nabla}_{u}u = A \frac{\partial}{\partial r} + B e_{2},
\end{equation}
where $A$ and $B$ are some smooth functions on the simply-connected open region $R_{\delta , \tau , \epsilon_0}$ as specified in \eqref{exteriorregion}. Notice that $\frac{\partial}{\partial r}$, and $e_{2} = \frac{a}{\sin (ar)} \frac{\partial}{\partial \theta }$ constitute an orthonormal moving frame on $S^2(a^2)$, it follows that
\begin{equation}
\begin{split}
A &= g(\overline{\nabla}_{u}u , \frac{\partial}{\partial r}  ) , \\
B & = g( \overline{\nabla}_{u}u , e_{2}    ).
\end{split}
\end{equation}
Now, since $|u|^2 = g(u, u) = [h(r-\delta )]^2$ is independent of the $\theta$ variable, it follows from condition \textbf{(1)} of
\textbf{Definition \ref{LeviCivitadefinition}} that we have
\begin{equation}
0 = e_{2}(g(u,u)) = 2 g (\overline{\nabla}_{e_{2}} u , u),
\end{equation}
from which it follows that
\begin{equation}
\begin{split}
B & = g( \overline{\nabla}_{u}u , e_{2}    ) \\
& = -h(r-\delta ) g ( \overline{\nabla}_{e_2}u , e_{2}) \\
& = g (\overline{\nabla}_{e_2}u ,   -h(r-\delta )  e_{2}) \\
& = g (\overline{\nabla}_{e_{2}} u , u)\\
& = 0 ,
\end{split}
\end{equation}
with the second equality follows from the tensorial property of $\overline{\nabla}$ as specified in condition \textbf{(4)} in \textbf{Definition \ref{LeviCivitadefinition}}.
In the same way, we have

\begin{equation}\label{computeA}
\begin{split}
A & = g (\overline{\nabla}_{u} u , \frac{\partial}{\partial r}) \\
& = -h(r-\delta ) \frac{a}{\sin (ar)} g ( \overline{\nabla}_{\frac{\partial}{\partial \theta}} u ,\frac{\partial}{\partial r}  ) \\
& = h(r-\delta ) \frac{a}{\sin (ar)} g (u , \overline{\nabla}_{\frac{\partial}{\partial \theta}} \frac{\partial}{\partial r} ) \\
& = h(r-\delta ) \frac{a}{\sin (ar)} g (u , \overline{\nabla}_{\frac{\partial}{\partial r}} \frac{\partial}{\partial \theta } ).
\end{split}
\end{equation}
In the above calculation, the fourth equal sign follows from the torsion-free property of $\overline{\nabla}$ $($condition \textbf{(2)} in \textbf{Definition \ref{LeviCivitadefinition}}$)$ which gives
$ \overline{\nabla}_{\frac{\partial}{\partial \theta}} \frac{\partial}{\partial r} =  \overline{\nabla}_{\frac{\partial}{\partial r}} \frac{\partial}{\partial \theta }$. Yet the validity of the third equal sign is based on the following observation, which on its own is a direct consequence of condition \textbf{(1)} of \textbf{Definition \ref{LeviCivitadefinition}}.
\begin{equation}
0 = \frac{\partial}{\partial \theta } g(u,\frac{\partial}{\partial r} ) = g ( \overline{\nabla}_{\frac{\partial}{\partial \theta}} u ,\frac{\partial}{\partial r}  ) +
g (u , \overline{\nabla}_{\frac{\partial}{\partial \theta}} \frac{\partial}{\partial r} ) .
\end{equation}
Now, to compute the term $\overline{\nabla}_{\frac{\partial}{\partial r}} \frac{\partial}{\partial \theta }$, one recall that $\overline{\nabla}_{\frac{\partial}{\partial r}} e_{2} = 0$ holds, since the restriction of $e_{2}$ along each geodesic ray $\gamma$ starting from $O \in S^2(a^2)$ must be parallel along $\gamma$. Hence, it follows from the Leibniz rule $($condition \textbf{(3)} in \textbf{Definition \ref{LeviCivitadefinition}} $)$ of the connection $\overline{\nabla}$ that the following relation holds.
\begin{equation}
\overline{\nabla}_{\frac{\partial}{\partial r}} \frac{\partial}{\partial \theta } = \overline{\nabla}_{\frac{\partial}{\partial r}} (\frac{\sin (ar)}{a} e_{2})
= \cos (ar) e_{2} .
\end{equation}
Hence, it follows from \eqref{computeA} that
\begin{equation}
\begin{split}
A & = h(r-\delta ) \frac{a}{\sin (ar)} g (u, \cos (ar) e_{2}) \\
& = h(r-\delta ) \frac{a}{\sin (ar)} \cos (ar) g (-h(r-\delta ) e_{2} , e_{2})\\
& = -[h(r-\delta )]^2 \frac{a \cos (ar)}{\sin (ar)}.
\end{split}
\end{equation}
So, finally, we have the following expression for $\overline{\nabla}_{u}u$
\begin{equation}
\overline{\nabla}_{u}u = -[h(r-\delta )]^2 \frac{a \cos (ar)}{\sin (ar)} \frac{\partial}{\partial r} ,
\end{equation}
which immediately gives the following expression for $[\overline{\nabla}_{u}u ]^*$,
\begin{equation}\label{connectionform}
[\overline{\nabla}_{u}u]^*  =  -[h(r-\delta )]^2 \frac{a \cos (ar)}{\sin (ar)} dr .
\end{equation}

\noindent{\bf Step 3 : The proof of Theorem \ref{Noparallel}. }

\noindent
With relations \eqref{triangleu}, and \eqref{connectionform}, we can now complete the proof of Theorem \ref{Noparallel} here. Assume towards contradiction that there exits some sufficiently small positive number $\epsilon_{0} < \frac{\pi}{a} - \delta $, such that there exists some smooth function $P$ defined on the sector-shaped region $R_{\delta , \tau , \epsilon_{0}}$ as specified in \eqref{exteriorregion} which solves the following stationary Navier-Stokes equation on $R_{\tau , \epsilon_{0}}$, with $\beta \in \mathbb{R}$ to be a given positive constant.
\begin{equation}\label{Navier1}
\nu (-\triangle u^* -2Ric (u^*)) + \beta \cos (ar) *u^* + \overline{\nabla}_{u}u^* + dP = 0.
\end{equation}
On the space form $S^2(a^2)$, we have $Ric (u^*) = (2-1)a^2 u^* = a^2 u^*$, and we can rephrase equation \eqref{Navier1} as follow.
\begin{equation}\label{Navier2}
\nu (-\triangle u^* -2 a^2 u^*) + \beta \cos (ar) *u^* + \overline{\nabla}_{u}u^* + dP = 0.
\end{equation}
According to the well-known fact that $d\circ d = 0$, the existence of a smooth function $P$ solving equation \eqref{Navier2} on the simply-connected open region $R_{\delta , \tau , \epsilon_{0}} $ as specified in \eqref{exteriorregion} immediately implies that the following relation holds everywhere in the region $R_{\delta , \tau , \epsilon_{0}}$,
\begin{equation}\label{wrong1}
d \big \{  \nu (-\triangle u^* -2 a^2 u^*) + \beta \cos (ar) *u^* + \overline{\nabla}_{u}u^*     \big   \} = 0.
\end{equation}
Recall that $u^* = - h(r-\delta ) \frac{\sin ar}{a} d\theta$. Hence, we have, through direct computation, that
\begin{equation}\label{curlu}
d u^* = -\{h'(r-\delta ) \frac{\sin (ar)}{a} + h(r-\delta ) \cos (ar)\} dr\wedge d\theta.
\end{equation}
Moreover, we deduce from \eqref{connectionform} that
\begin{equation}\label{trivialnew}
d (\overline{\nabla}_{u}u^*) = 0.
\end{equation}
In addition, it is also obvious that we always have the following relation
\begin{equation}\label{drotat}
d \big \{ \beta \cos (ar) *u^* \big \} = \beta \frac{\partial}{\partial \theta } \big \{ \cos (ar) h(r-\delta )\big \} d\theta \wedge dr = 0 .
\end{equation}
Now, in accordance with \eqref{triangleu},\eqref{trivialnew} , \eqref{curlu}, and \eqref{drotat}, we have
\begin{equation}\label{generalexp}
\begin{split}
&d \big \{ \nu (-\triangle u^* -2 a^2 u^*) + \beta \cos (ar) *u^* +  \overline{\nabla}_{u}u^*    \big \}\\
& = \nu \bigg\{h'''(r-\delta ) \frac{\sin (ar)}{a}  + 2 h''(r-\delta ) \cos (ar)  + a h'(r-\delta ) \left(\sin (ar) -\frac{1}{\sin (ar)}\right) \\
& + a^2 h(r-\delta ) \cos (ar)\left(2+ \frac{1}{\sin^2(ar)}\right)\bigg\} dr\wedge d\theta .
\end{split}
\end{equation}
Now, recall that $h(\lambda)$ is chosen to be $h(\lambda) = \alpha_{1} \lambda -\frac{\alpha_{2}}{2}\lambda^2$ in Theorem \ref{Noparallel}. So, we have $h'(\lambda ) = \alpha_{1} - \alpha_{2} \lambda$, $h''(\lambda ) = - \alpha_{2}$, and eventually $h'''(\lambda ) = 0$, for every $\lambda \geqslant 0$. In which case, \eqref{generalexp} reduces down to
\begin{equation}
\begin{split}
&d \big \{ \nu (-\triangle u^* -2 a^2 u^*) + \beta \cos (ar) *u^* + \overline{\nabla}_{u}u^*   \big \}\\
& = \nu \bigg\{ -2 \alpha_{2}  \cos (ar)  + a (\alpha_{1} -\alpha_{2}(r-\delta )) \left(\sin (ar) -\frac{1}{\sin (ar)}\right) \\
& + a^2 (r-\delta )\left[\alpha_{1} -\frac{\alpha_{2}(r-\delta )}{2}\right] \cos (ar)\left(2+ \frac{1}{\sin^2(ar)}\right)\bigg\} dr\wedge d\theta .
\end{split}
\end{equation}
For convenience, we will use the following abbreviation
\begin{equation}\label{expressionF}
\begin{split}
F_{\alpha_{1}, \alpha_{2}, \delta } (r)
& = \bigg\{ -2 \alpha_{2}  \cos (ar) \\
& + a (\alpha_{1} -\alpha_{2}(r-\delta )) \left(\sin (ar) -\frac{1}{\sin (ar)}\right) \\
& + a^2 (r-\delta )\left[\alpha_{1} -\frac{\alpha_{2}(r-\delta )}{2}\right] \cos (ar)\left(2+ \frac{1}{\sin^2(ar)}\right)\bigg\} ,
\end{split}
\end{equation}
so that $F_{\alpha_{1}, \alpha_{2}, \delta }$ is a smooth function defined on $(0,\frac{\pi}{a})$, and that we have
\begin{equation}\label{cheapexpression}
d \big \{ \nu (-\triangle u^* -2 a^2 u^*) + \beta \cos (ar) *u^* + \overline{\nabla}_{u}u^*    \big \} = \nu F_{\alpha_{1}, \alpha_{2}, \delta } (r) dr\wedge d\theta .
\end{equation}
Now, by insisting on the existence of a smooth $P$ which solves \eqref{Navier2} on the simply-connected region $R_{\delta , \tau , \epsilon_{0}}$ as specified in \eqref{exteriorregion}, the validity of \eqref{wrong1} on $R_{\delta , \tau , \epsilon_{0}}$ will follow as a by-product, which in turns forces us to admit that the following relation should hold for every $r \in  [\delta , \delta + \epsilon_{0} )$,
\begin{equation}\label{vanishingofF}
F_{\alpha_{1}, \alpha_{2}, \delta } (r) = 0.
\end{equation}
That is, $F_{\alpha_{1}, \alpha_{2}, \delta } $ should identically vanish on $r \in  [\delta , \delta + \epsilon_{0} )$. In that which follows, we will split our argument into three cases, namely the case of $0< \delta < \frac{\pi}{2a}$, then the case of $\delta = \frac{\pi}{2a}$, and finally the case of $\frac{\pi}{2a} < \delta < \frac{\pi}{a}$. In each of these cases, we will derive a contradiction towards the validity of the vanishing property of $F_{\alpha_{1}, \alpha_{2}, \delta } $ on $[\delta , \delta + \epsilon_{0} )$.\\

\noindent
\textbf{Case One.} We first discuss the case of $0 < \delta < \frac{\pi}{2a} $. In this case, we simply observe that we have $\cos (a\delta ) > 0$, and $\frac{1}{\sin (a\delta )} - \sin (a\delta ) > 0$. Based upon such an observation, we deduce at once that the following property holds as long as $\alpha_{1} > 0$, and $\alpha_{2} > 0$,
\begin{equation}
F_{\alpha_{1}, \alpha_{2}, \delta }(\delta ) = -2\alpha_{2} \cos (a \delta ) -a \alpha_{1} \left(\frac{1}{\sin (a\delta )} - \sin (a\delta )\right) < 0.
\end{equation}
The validity of the above relation at once implies that, for some sufficiently small $\epsilon_{1} \in (0,\epsilon_{0})$, we will have the following property
\begin{itemize}
\item $F_{\alpha_{1}, \alpha_{2}, \delta }(r) < 0$ holds, for all $r \in [\delta , \delta + \epsilon_{1})$,
\end{itemize}
which directly contradicts the \emph{everywhere vanishing property} of $F_{\alpha_{1}, \alpha_{2}, \delta }(r)$ on $[\delta , \delta + \epsilon_{0})$. So, in this case, a contradiction has been arrived, which ensures the non-existence of a smooth $P$ solving \eqref{Navier2} on the sector-shaped region $R_{\delta , \tau , \epsilon_{0}}$ as specified in \eqref{exteriorregion}, regardless of how small its angle $\tau$ or thickness $\epsilon_{0}$ would be.\\

\noindent
\textbf{Case Two.} We now deal with the case of $\delta = \frac{\pi}{2a}$. The problem involved here is that $F_{\alpha_{1}, \alpha_{2}, \frac{\pi}{2a} }(\frac{\pi}{2a}) =0 $. So, we look at the quantity $\frac{\partial}{\partial r}(\sin^2(ar)F_{\alpha_{1}, \alpha_{2}, \frac{\pi}{2a} })|_{r= \frac{\pi}{2a}}$ instead.

First, we have
\begin{equation}\label{expressionSinarF}
\begin{split}
\sin^2(ar) F_{\alpha_{1}, \alpha_{2}, \frac{\pi}{2a} } (r)
& =  -2 \alpha_{2}  \cos (ar) \sin^2(ar) \\
& + a \left(\alpha_{1} -\alpha_{2}\left(r- \frac{\pi}{2a}\right )\right) (\sin^3(ar) -\sin (ar)) \\
& + a^2\left (r-  \frac{\pi}{2a}\right )\left[\alpha_{1} -\frac{\alpha_{2}}{2}\left(r-  \frac{\pi}{2a} \right)\right] \cos (ar)(2 \sin^2(ar)+ 1 ) .
\end{split}
\end{equation}
Now, observe that
\begin{equation}\label{boring}
\begin{split}
\frac{\partial }{\partial r} \{ -2 \alpha_{2}  \cos (ar) \sin^2(ar)  \}\bigg |_{r= \frac{\pi}{2a}} &= 2a \alpha_{2} , \\
\frac{\partial }{\partial r} \left\{ a \left(\alpha_{1} -\alpha_{2}\left(r- \frac{\pi}{2a} \right)\right) (\sin^3(ar) -\sin (ar))    \right \} \bigg |_{r= \frac{\pi}{2a}} & = 0 ,\\
\frac{\partial }{\partial r} \left\{a^2 \left(r-  \frac{\pi}{2a} \right)\left[\alpha_{1} -\frac{\alpha_{2}}{2}\left(r-  \frac{\pi}{2a} \right)\right] \cos (ar)(2 \sin^2(ar)+ 1 )\right\} \bigg |_{r= \frac{\pi}{2a}} & = 0 .
\end{split}
\end{equation}
Hence, it follows from \eqref{expressionSinarF}, and \eqref{boring} that we have
\begin{equation}
\frac{\partial}{\partial r}(\sin^2(ar)F_{\alpha_{1}, \alpha_{2}, \frac{\pi}{2a} })|_{r= \frac{\pi}{2a}} = 2a\alpha_2 > 0 ,
\end{equation}
which, together with $F_{\alpha_{1}, \alpha_{2}, \frac{\pi}{2a} }(\frac{\pi}{2a}) =0 $,  imply that the following property holds for some sufficiently small $\epsilon_{1} \in (0, \epsilon_{0})$,
\begin{itemize}
\item $\sin^2(ar)F_{\alpha_{1}, \alpha_{2}, \frac{\pi}{2a} } > 0 $ holds for all $r \in (\frac{\pi}{2a} , \frac{\pi}{2a} + \epsilon_{1}) $.
\end{itemize}
However, the above property is again in direct conflict with the fact that $F_{\alpha_{1}, \alpha_{2}, \frac{\pi}{2a} }$ should identically vanish on $[\frac{\pi}{2a} , \frac{\pi}{2a} + \epsilon_{0})$, should there be a smooth $P$ solving equation \eqref{Navier2} on the sector-shaped region $R_{\delta , \tau , \epsilon_{0}}$ as specified in \eqref{exteriorregion}. Again, this contradiction ensures the non-existence of a smooth function $P$ solving \eqref{Navier2} on $R_{\delta , \tau , \epsilon_{0}}$ as specified in \eqref{exteriorregion}.\\

\noindent
\textbf{Case Three.} We now consider the case of $\frac{\pi}{2a} < \delta < \frac{\pi}{a}$, which is the most delicate one among the three cases. The delicate issue here is that $\cos (a \delta ) < 0$. Here, we will use the basic theory of second order linear ODE to treat this case. As in the previous cases, by insisting on the existence of a smooth $P$ solving \eqref{Navier2} on $R_{\delta , \tau , \epsilon_{0}}$ as specified in \eqref{exteriorregion}, we will arrive at the consequence that the function $F_{\alpha_1 , \alpha_2 , \delta }$ must identically vanish over $[\delta , \delta + \epsilon_{0})$. But this is the same as saying that the quadratic function $Y(r) = \alpha_1 (r-\delta) -\frac{\alpha_2}{2}(r-\delta )^2$
will be a \emph{local solution} to the following linear second order ODE on the interval $[\delta , \delta + \epsilon_{0}  )$.
\begin{equation}\label{linearODE}
y''(r) + Q_{1}(r) y'(r) + Q_{2}(r) y(r) = 0,
\end{equation}
where $Q_1$, $Q_2$ are the real-analytic functions on $(\frac{\pi}{2a} , \frac{\pi}{a})$ defined by
\begin{equation}
\begin{split}
Q_1(r) & = \frac{a}{2\cos (ar)} \left(\sin (ar) - \frac{1}{\sin (ar)}  \right ) ,\\
Q_2(r) & = \frac{a^2}{2} \left(2+ \frac{1}{\sin^2(ar)}\right).
\end{split}
\end{equation}
In accordance with the basic existence and uniqueness theory for second order ODE's, the local solution $Y$ on $[\delta , \delta + \epsilon_{0})$ as represented by the quadratic function $Y(r) = \alpha_1 (r-\delta) -\frac{\alpha_2}{2}(r-\delta )^2$ can be uniquely extended to a global solution $Z$ to equation \eqref{linearODE} on the whole interval $(\frac{\pi}{2a} , \frac{\pi}{a})$. In addition, since the coefficient functions $Q_1(r)$, and $Q_2(r)$ are real-analytic on $(\frac{\pi}{2a} , \frac{\pi}{a})$, such a global solution $Z$ must also be real analytic on $(\frac{\pi}{2a} , \frac{\pi}{a}) $.
So, the power series representation $Z(r) = \sum_{k=0}^{\infty}a_{k} (r - \frac{3\pi}{4a})^k$ of the real-analytic solution $Z$ about the point $\frac{3\pi}{4}$ should have its radius of convergence to be at least $\frac{\pi}{4a}$. That is, $Z$ can be represented by $\sum_{k=0}^{\infty}a_{k} (r - \frac{3\pi}{4a})^k$, which converges absolutely for all $r \in (\frac{\pi}{2a} , \frac{\pi}{a})$. Now, we consider the following two holomorphic functions, which are the complexifications of the real analytic functions $Y(r) = \alpha_1 (r-\delta) -\frac{\alpha_2}{2}(r-\delta )^2$ and $Z(r) = \sum_{k=0}^{\infty}a_{k} (r - \frac{3\pi}{4a})^k$ respectively.
\begin{equation}
\begin{split}
\textbf{Y}(w) &= \alpha_1 (w-\delta) -\frac{\alpha_2}{2}(w-\delta )^2 , \\
\textbf{Z} (w) & = \sum_{k=0}^{\infty}a_{k} \left(w - \frac{3\pi}{4a}\right)^k .
\end{split}
\end{equation}
Since the radius of convergence of $\sum_{k=0}^{\infty}a_{k} (r - \frac{3\pi}{4a})^k$ is at least $\frac{\pi}{4a}$, the holomorphic function $\textbf{Z} (w)  = \sum_{k=0}^{\infty}a_{k} (w - \frac{3\pi}{4a})^k$ is well-defined \emph{at least} on the open ball $\{w \in \mathbb{C} : |w - \frac{3\pi}{4a}| < \frac{\pi}{4a}\}$ in $\mathbb{C}$. Recall that the real analytic solution $Z$ to \eqref{linearODE} arises as the unique extension of the local solution $Y(r) = h(r-\delta )$ in that
\begin{equation}
Z|_{[\delta , \delta + \epsilon_{0})} =Y ,
\end{equation}
which means the same as saying that the two holomorphic functions $\textbf{Z}$ and $\textbf{Y}$ coincide on the line segment $\{r : \delta < r < \delta + \epsilon_{0} \}$, which by itself is included in the open disc $\{w \in \mathbb{C} : |w - \frac{3\pi}{4a}| < \frac{\pi}{4a}\}$. So, it follows from the \emph{identity theorem} of the complex function theory that $\textbf{Z}$ must be identical to the entire function $\textbf{Y}$ on $\mathbb{C}$, which basically says that the power series
$\sum_{k=0}^{\infty}a_{k} (w - \frac{3\pi}{4a})^k$ is identical to $ \alpha_1 (r-\delta) -\frac{\alpha_2}{2}(r-\delta )^2  $, for all $w \in \mathbb{C}$. So, we deduce that the global real-analytic solution $Z$ to \eqref{linearODE} must be \emph{identical to} the quadratic function $\alpha_1 (r-\delta) -\frac{\alpha_2}{2}(r-\delta )^2$ over the \emph{whole} interval $(\frac{\pi}{2a} , \frac{\pi}{a})$, which is the same as saying that we have the following identity for \emph{all} $r \in (\frac{\pi}{2a} , \frac{\pi}{a})$, with $Y(r)$ to be the quadratic function $Y(r) = \alpha_1 (r-\delta) -\frac{\alpha_2}{2}(r-\delta )^2$,
\begin{equation}\label{goodidentity}
Y''(r) + Q_{1}Y'(r) + Q_2(r)Y(r) = 0.
\end{equation}
To finish the argument, we will derive a contradiction against identity \eqref{goodidentity} through investigating the limiting behavior of $Y''(r) + Q_{1}(r)Y'(r) + Q_2(r)Y(r)$ as $r \rightarrow \frac{\pi}{a}^-$. Now, we will further split the discussion into two subcases subordinate to \textbf{Case Three}.
First, we consider the subcase when $\alpha_1 - \frac{\alpha_2}{2}(\frac{\pi}{a} -\delta ) $ is \emph{not} zero. In this subcase, we have the following relations, which follow from direct computations.
\begin{equation}
\begin{split}
\lim_{r\rightarrow \frac{\pi}{a}^-} \sin (ar) Y''(r)& = 0, \\
\lim_{r\rightarrow \frac{\pi}{a}^-} \sin (ar) Q_{1}(r) Y'(r) & = \frac{a}{2}\left[\alpha_1 - \alpha_2 \left(\frac{\pi}{a} -\delta \right)\right] \\
\lim_{r\rightarrow \frac{\pi}{a}^-} \sin (ar) Q_2(r) Y(r) & = \infty ,
\end{split}
\end{equation}
from which it follows at once that
\begin{equation}
\lim_{r\rightarrow \frac{\pi}{a}^-} \sin(ar) [Y''(r) + Q_{1}Y'(r) + Q_2(r)Y(r)] = \infty ,
\end{equation}
which is in direct conflict with identity \eqref{goodidentity}, which is supposed to hold on the \emph{whole interval} $(\frac{\pi}{2a} , \frac{\pi}{a})$. So, in this subcase, we can rule out the possibility of having a smooth function $P$ solving \eqref{Navier2} on the simply-connected open region $R_{\delta , \tau , \epsilon_{0}}$ as specified in \eqref{exteriorregion}.\\

\noindent
Next, we deal with the remaining subcase when $\alpha_1 = \frac{\alpha_2}{2}(\frac{\pi}{a} -\delta )$. For this case, it is immediate to see that
\begin{equation}
\begin{split}
\alpha_1 - \frac{\alpha_2}{2}(r-\delta ) &= \frac{\alpha_2}{2} \left(\frac{\pi}{a} - r\right), \\
\alpha_{1} - \alpha_2 (r-\delta ) &= \frac{\alpha_2}{2} \left( \frac{\pi}{a} -2r + \delta \right).
\end{split}
\end{equation}
Hence, it follows that
\begin{equation}\label{cosexpression}
\begin{split}
& \cos (ar) [Y''(r) + Q_1(r)Y'(r) + Q_2(r)Y(r)] \\
&= -2 \alpha_2 \cos (ar) + \frac{a \alpha_2}{2}\left(\frac{\pi}{a} - 2r + \delta \right) \sin (ar) + a^2 \alpha_2 (r-\delta )
\left(\frac{\pi}{a} -r \right)\cos (ar) \\
& + \frac{a \alpha_2}{2 \sin^2(ar)}\left\{ \left(2r - \frac{\pi}{a} - \delta \right) \sin (ar) -a (r-\delta )\left(r- \frac{\pi}{a}\right) \cos (ar)            \right \}
\end{split}
\end{equation}
However, we see that the following relations follow from direct computation.
\begin{equation}
\begin{split}
& \lim_{r\rightarrow \frac{\pi}{a}^-}  \left\{-2 \alpha_2 \cos (ar) + \frac{a \alpha_2}{2}\left(\frac{\pi}{a} - 2r + \delta \right) \sin (ar) + a^2 \alpha_2 (r-\delta )\left(\frac{\pi}{a} -r \right)\cos (ar) \right\} = 2 \alpha_2 \\
& \lim_{r\rightarrow \frac{\pi}{a}^-} \frac{a \alpha_2}{2 \sin^2(ar)}\left\{ \left(2r - \frac{\pi}{a} - \delta \right) \sin (ar) -a (r-\delta )\left(r- \frac{\pi}{a}\right) \cos (ar)         \right    \}
 = -\frac{\alpha_2}{2} ,
\end{split}
\end{equation}
from which we deduce that
\begin{equation}
\lim_{r\rightarrow \frac{\pi}{a}^-} \cos (ar) [Y''(r) + Q_1(r)Y'(r) + Q_2(r)Y(r)] = \frac{3 \alpha_2}{2} \in \mathbb{R}-\{0\},
\end{equation}
which again is in direct conflict with identity \eqref{goodidentity}, whose validity is supposed to be on the whole interval $(\frac{\pi}{2a} , \frac{\pi}{a})$. So, this contradiction rules out the possibility of having a smooth $P$ solving \eqref{Navier2} on the sector-shaped region $R_{\delta , \tau ,\epsilon_0 }$ as specified in \eqref{exteriorregion}. So, we have finished the proof for Theorem \ref{Noparallel}.
In the proof of Theorem \ref{Noparallel}, we have already used the following basic existence and uniqueness theorem for linear ODE with real analytic coefficients, which can be found in standard ODE textbooks.

\begin{thm}\label{ODEtheorem}
Consider the following linear equation about the unknown solution $Y(t)$:
\begin{equation}\label{ODE}
Y^{(n)}(t) + Q_{n-1}(t) Y^{(n-1)}(t)  + Q_{n-2}(t)  Y^{(n-2)}(t) + .... + Q_0(t) Y(t) = 0 ,
\end{equation}
where $Q_{j}$ are the prescribed real analytic coefficient functions defined on some open interval $(a, b)$. Here, the symbol $Y^{(j)}(t)$ stands for the $j$-order derivative of $Y$. Let $\delta \in (a,b)$ to be some selected based point in $(a,b)$. Then, for any prescribed real numbers $\beta_0$, $\beta_{1}$, $\beta_{2}$, ... $\beta_{n-1}$, there exists a unique real analytic solution $Y : (a, b) \rightarrow \mathbb{R}$ to the linear ODE \eqref{ODE} which at the same time satisfies the initial values $Y^{(j)}(\delta ) = \beta_j$, for every $0 \leq j \leq n-1$.
\end{thm}

\noindent
Now, with the help of Theorem \ref{ODEtheorem}, we can now give a very brief proof for Theorem \ref{existenceEasy} as follow.\\

\noindent{\bf Step 4 : The proof of Theorem \ref{existenceEasy}.}

\noindent
To begin, we again consider the following representation formula, which we obtain through the efforts spent in \textbf{Step 1} and \textbf{Step 2} of this Section, and which is valid for \emph{any} parallel laminar flow $u = -Y(r) \frac{a}{\sin (ar)}\frac{\partial}{\partial \theta }$, with $Y(r) = h(r-\delta )$ to be \textbf{any possible} smooth function defined on some open interval about the based point $\delta$.
\begin{equation}\label{generalexpSECOND}
\begin{split}
&d \{ \nu (-\triangle u^* -2 a^2 u^*) + \beta \cos (ar)* u^* + \overline{\nabla}_{u}u^*     \} \\
& = \nu \bigg\{Y'''(r) \frac{\sin (ar)}{a}  + 2 Y''(r) \cos (ar)  + a Y'(r) \left(\sin (ar) -\frac{1}{\sin (ar)}\right) \\
& + a^2 Y(r) \cos (ar)\left(2+ \frac{1}{\sin^2(ar)}\right)\bigg\} dr\wedge d\theta .
\end{split}
\end{equation}
As we have already mentioned in the introduction, to see whether equation \eqref{Navier2} will admit a solution in the form of $u = -Y(r) \frac{a}{\sin (ar)}\frac{\partial}{\partial \theta }$ on the simply-connected open region $\Omega_{\delta , \tau }$ as specified in \eqref{OMEGASPHERE}, with $Y(r) = h(r-\delta )$ to be some unknown function, it is enough to check whether the smooth $2$-form $d \big \{ \nu (-\triangle u^* -2 a^2 u^*) + \beta \cos (ar) *u^* + \overline{\nabla}_{u}u^*   \big   \}$ will vanish identically over the same simply-connected open region $\Omega_{\delta , \tau } \subset S^2(a^2)$. This is because the $d$-closed property of $\nu (-\triangle u^* -2 a^2 u^*) + \beta \cos (ar)* u^* + \overline{\nabla}_{u}u^*$ over any simply-connected open region in $S^2(a^2)-K$ will at once ensure the existence of a smooth pressure function $P$ which solves equation \eqref{Navier2} on the same simply-connected open region. In order words, $u = -Y(r) \frac{a}{\sin (ar)}\frac{\partial}{\partial \theta }$ will solve equation \eqref{Navier2} on the simply-connected open region $\Omega_{\delta , \tau }$ as specified in \eqref{OMEGASPHERE}  if and only if the unknown function $Y(r) = h(r-\delta )$ will satisfy the following $3$-order linear ODE on $[\delta , \frac{\pi}{a})$.
\begin{equation}\label{thirdorderODESECOND}
Y'''(r)  + 2 a Y''(r) \frac{\cos (ar)}{\sin (ar)} + a^2 Y'(r) \left(1 -\frac{1}{\sin^2 (ar)}\right) + a^3 Y(r) \frac{\cos (ar)}{\sin (ar)}\left(2+ \frac{1}{\sin^2(ar)}\right) = 0 .
\end{equation}
However, it is quite obvious that all the coefficient functions appearing in the above $3$-order linear ODE are all real-analytic functions on the open interval $(0,\frac{\pi}{a})$. As a result, Theorem \ref{ODEtheorem} ensures that there exits a \emph{unique} real analytic function $Y$ defined on the whole interval $(0,\frac{\pi}{a})$ which solves the $3$-order ODE \eqref{thirdorderODESECOND}, and which satisfies the prescribed initial values $Y(\delta ) = 0$, $Y'(\delta ) = \alpha_1$, and $Y''(\delta ) = - \alpha_2$, where $\alpha_1$, $\alpha_2$ are some arbitrary given positive numbers. For such a unique analytic solution $Y : (0,\frac{\pi}{a}) \rightarrow \mathbb{R}$ to equation \eqref{thirdorderODESECOND}, the associated parallel laminar flow $u = -Y(r) \frac{a}{\sin (ar)}\frac{\partial}{\partial \theta }$ will solve equation \eqref{Navier2} on the simply-connected open region $\Omega_{\delta , \tau }$ as specified in \eqref{OMEGASPHERE}. So, the proof of Theorem \ref{existenceEasy} is completed.

\section{Basic geometry of $\mathbb{H}^2(-a^2)$ : The visualization of $\mathbb{H}^2(-a^2)$ through the use of the hyperboloid model, and the geodesic normal polar coordinate on $\mathbb{H}^2(-a^2)$.}\label{HYPERBOLIDSECTION}

\noindent
 The objective of this section is to give a brief introduction to the geometry of the space form $\mathbb{H}^2(-a^2)$ of constant negative sectional curvature $-a^2 < 0$. Again, all the materials as presented in this section are standard facts well-known to differential geometers. The purpose of this section, however, is to spell out the basic notions and geometric language which we will use in those differential-geometric calculations in \textbf{Section \ref{hyperMAINSECTION}}. We will first give a description of $\mathbb{H}^2(-a^2)$ through the use of the so-called hyperboloid model.\\

\noindent
In the following presentation, we closely follow the standard construction of the hyperboloid model of $\mathbb{H}^2(-a^2)$ as given in pages 201 to 202 of \cite{Jost}.\\

\begin{defn}\label{hyperbolid}
\textbf{(The characterization of $\mathbb{H}^2(-a^2)$ by means of the hyperboloid model).} Consider the \emph{linear space}
$\mathbb{V}^3 = \{ (x_0, x_1 , x_2) : x_0, x_1, x_2 \in \mathbb{R}  \}$ which is equipped with the following quadratic form $<\cdot , \cdot > : \mathbb{V}^3\otimes \mathbb{V}^3 \rightarrow \mathbb{R}$.
\begin{equation}\label{quadraticform}
<x,y> =  -x_0y_0 + x_1y_1 + x_2y_2 ,
\end{equation}
with $x$, $y$ to be any two elements in the linear space $\mathbb{V}^3$.
\end{defn}

\noindent
Notice that $\mathbb{V}^3$ which is equipped with quadratic form \eqref{quadraticform} is \textbf{not} the same as the Euclidean space $\mathbb{R}^3$. Then, we define $\mathbb{H}^2(-a^2)$ as follow.
\begin{equation}
\mathbb{H}^2(-a^2) = \big \{ x \in \mathbb{V}^3 : <x,x> = \frac{-1}{a^2}, x_0 > 0 \big \} .
\end{equation}
Then, it is clear that $\mathbb{H}^2(-a^2)$ is represented as a branch of the hyperboloid $<x,x> = \frac{-1}{a^2}$. Hence, topologically, $\mathbb{H}^2(-a^2)$, as a differentiable manifold on its own, is diffeomorphic to $\mathbb{R}^2$ (Be careful, $\mathbb{H}^2(-a^2)$ as a \textbf{Riemannian manifold} is \textbf{not} the same as the Euclidean space $\mathbb{R}^2$). Next, we would like to construct the Riemannian metric $g(\cdot, \cdot )$ on the $2$-dimensional manifold $\mathbb{H}^2(-a^2)$ as follow. Here, for any point $p \in \mathbb{H}^2(-a^2)$, we consider the symmetric bilinear form $\big ( -dx_0\otimes dx_0 + dx_1\otimes dx_1 + dx_2\otimes dx_2 \big )\big |_{p}$ acting on the tangent space $T_p\mathbb{V}^3$ of $\mathbb{V}^3$ at $p$, which is described as follow.
\begin{equation}
\big ( -dx_0\otimes dx_0 + dx_1\otimes dx_1 + dx_2\otimes dx_2 \big )\big |_{p}(v,w) = -v_0w_0 + v_1w_1 + v_2w_2 ,
\end{equation}
where $v$, $w \in T_p\mathbb{V}^3$. Then, for the point $p \in \mathbb{H}^2(-a^2)$, we consider the positive definite inner product $g_{p}(\cdot, \cdot ) $ on the tangent space $T_{p}\mathbb{H}^2(-a^2)$, which is \emph{defined} to be the restriction of the symmetric bilinear form $\big ( -dx_0\otimes dx_0 + dx_1\otimes dx_1 + dx_2\otimes dx_2 \big )\big |_{p}$ onto the vector subspace $T_{p}\mathbb{H}^2(-a^2)$ of $T_p\mathbb{V}^3$. Then, this smoothly varying family of positive definite inner products $g_{p}(\cdot , \cdot ) : T_{p}\mathbb{H}^2(-a^2)\otimes T_{p}\mathbb{H}^2(-a^2) \rightarrow \mathbb{R}$, for $p \in \mathbb{H}^2(-a^2)$ \emph{constitutes} the Riemannian metric $g(\cdot , \cdot )$ on $\mathbb{H}^2(-a^2)$.
In order to further clarify the content of the above definition of $\mathbb{H}^2(-a^2)$, a few remarks are in order here.
\begin{remark}
The $2$-dimensional space form $\mathbb{H}^2(-a^2)$ with constant sectional curvature $-a^2 < 0$ is often called the $2$-dimensional hyperbolic space (or hyperbolic manifold) of constant sectional curvature $-a^2$. If we consider the group $O(2,1)$ which consists of all those linear maps on $\mathbb{V}^3$ which leave the quadratic form $<\cdot, \cdot >$ as specified in \eqref{quadraticform} invariant, and which map the $x_0$-axis onto itself, then this group $O(2,1)$ will leave $\mathbb{H}^2(-a^2)$ invariant. Actually, $O(2,1)$ is exactly the group of isometries which acts transitively on $\mathbb{H}^2(-a^2)$ (see discussions on page 202 of the textbook \cite{Jost}). So, it turns out that the geometric structure of the hyperbolic space $\mathbb{H}^2(-a^2)$ around any selected point $p \in \mathbb{H}^2(-a^2)$ looks \textbf{exactly the same}, regardless of where the selected point $p$ of reference is. This just says that, for any two points $p$, $q$ in $\mathbb{H}^2(-a^2)$, the geometric structure of $\mathbb{H}^2(-a^2)$ around $p$ is identical to the geometric structure of $\mathbb{H}^2(-a^2)$ around $q$, up to an isometry $T$ from $O(2,1)$ sending $p$ to $T(p)=q$. So, in the following discussion, we can just, \emph{without the loss of generality}, choose the preferred reference point $O$ in $\mathbb{H}^2(-a^2)$ to be just $O = (\frac{1}{a^2},0,0)$.
\end{remark}

\noindent
Here, for the clarity of our presentation, we just select our preferred based point $O$ in $\mathbb{H}^2(-a^2)$ to be $O = (\frac{1}{a^2},0,0)$. That is, $O$ is located at the vertex of the hyperboloid $\{x \in \mathbb{V}^3 : <x,x> = \frac{1}{a^2} , x_0 > 0 \}$. We stress again that such a choice of $O$ is a choice without the loss of generality, due to the homogeneous structure of the hyperbolic manifold $\mathbb{H}^2(-a^2)$. Now, we will give a concrete description of the exponential map $\exp_{O} : T_{O}\mathbb{H}^2(-a^2)\rightarrow \mathbb{H}^2(-a^2)$ which maps the tangent space $T_{O}\mathbb{H}^2(-a^2)$ of $\mathbb{H}^2(-a^2)$ at $O$ onto the manifold $\mathbb{H}^2(-a^2)$ itself. Recall that such an exponential map is defined abstractly in the following manner.

\begin{defn}\label{exphyper}
\textbf{The exponential map on the hyperbolic space $\mathbb{H}^2(-a^2)$ :} For any $v \in T_{O}\mathbb{H}^2(-a^2)$ with $\|v\| = 1$, we consider the uniquely determined (\emph{unit speed}) geodesic $c_{v}: [0,\infty )\rightarrow \mathbb{H}^2(-a^2)$ which satisfies $c_v(0) = O$, and
$\dot{c}(0) = v$. Then, it follows that for any $r > 0$, $\exp_{O}(rv)$ is \emph{defined} to be
\begin{equation}\label{exponentialmaphyper}
\exp_{O}(rv) = c_v(r) .
\end{equation}
\end{defn}

\noindent
It is a basic differential geometric fact that the exponential map $\exp_{O} :  T_{O}\mathbb{H}^2(-a^2)\rightarrow \mathbb{H}^2(-a^2) $ as specified in \textbf{Definition \ref{exphyper}} is a smooth bijective map of the tangent space $T_{O}\mathbb{H}^2(-a^2)$ onto $\mathbb{H}^2(-a^2)$, whose inverse map $\exp^{-1} : \mathbb{H}^2(-a^2) \rightarrow T_{O}\mathbb{H}^2(-a^2)$ is also smooth. That is, $exp_{O}$ is a diffeomorphism from $T_{O}\mathbb{H}^2(-a^2) = \mathbb{R}^2$ onto $\mathbb{H}^2(-a^2)$.\\

\noindent
Now, by means of the concrete hyperboloid model of $\mathbb{H}^2(-a^2)$ as given in \textbf{Definition \ref{hyperbolid}}, for the reference point
$O = (\frac{1}{a^2} , 0 , 0)$, the exponential map $\exp_{O} : T_{O}\mathbb{H}^2(-a^2)\rightarrow \mathbb{H}^2(-a^2)$ can be expressed concretely as follow. Notice that the tangent space $T_{O}\mathbb{H}^2(-a^2)$ can naturally be identified with the plane $\{ (0, v_1, v_2) : v_1 , v_2 \in \mathbb{R} \}$ in the linear space $\mathbb{V}^3$. So, for any vector $v = (0, v_1 ,v_2) \in T_{O}\mathbb{H}^2(-a^2)$ with $\|v\| = 1$, the geodesic
$c_v : [0,\infty ) \rightarrow \mathbb{H}^2(-a^2)$ satisfying $c_v(0) = O$, and $\dot{c}(0) = v$ is expressed concretely as follow.
\begin{equation}\label{concretegeohyper}
c_v(r) = \cosh (ar) (\frac{1}{a} , 0,0) + \frac{\sinh (ar)}{a} (0, v_1 , v_2).
\end{equation}
Hence, for any unit vector $v$ in $T_{O}\mathbb{H}^2(-a^2)$, and any $r > 0$, the term $\exp_O(rv)$ is given explicitly as follow.
\begin{equation}\label{concreteexphyper}
\exp_{O}(rv) = \cosh (ar) (\frac{1}{a} , 0,0) + \frac{\sinh (ar)}{a} (0, v_1 , v_2).
\end{equation}
We can now define the \emph{normal polar coordinate system} $(r , \theta )$ on $\mathbb{H}(-a^2)$ about the reference point $O$ as follow.

\begin{defn}\label{polarDefn}
\textbf{The normal polar coordinate system on $\mathbb{H}^2(-a^2)$ :} The normal polar coordinate system on $\mathbb{H}^2(-a^2)$ about the point $O$ is the smooth map $(r,\theta ) : \mathbb{H}^2(-a^2) \rightarrow (0, \infty )\times (0, 2\pi )$ defined by
\begin{equation}
(r,\theta) = (\overline{r} , \overline{\theta }) \circ \exp_{O}^{\Red{-1}} ,
\end{equation}
in which $(\overline{r} , \overline{\theta })$ is the standard polar coordinate system on the Euclidean $2$-space $\mathbb{R}^2$ (here, we identify the tangent space $T_{O}\mathbb{H}^2(-a^2)$ with $\mathbb{R}^2$).
\end{defn}

\noindent
By means of such a normal polar coordinate system $(r,\theta )$ on $\mathbb{H}^2(-a^2)$, we can define two natural vector fields
$\frac{\partial}{\partial r}$ and $\frac{\partial}{\partial \theta}$ on $\mathbb{H}^2(-a^2)$ as follow.

\begin{defn}\label{FRAMEHyperDef}
\textbf{Natural coordinate frame $\{ \frac{\partial}{\partial r} , \frac{\partial}{\partial \theta }  \}$ on $\mathbb{H}^2(-a^2)$ :} For each $p \in \mathbb{H}^2(-a^2)$, the vectors $\frac{\partial}{\partial r}\big |_{p}$ and $\frac{\partial}{\partial \theta } \big |_p$ in $T_p \mathbb{H}^2(-a^2)$ are defined as linear derivations acting on the space $C^{\infty}(\mathbb{H}^2(-a^2))$ of smooth functions on $\mathbb{H}^2(-a^2)$ through the following relations.
\begin{equation}\label{naturalvectorhyper}
\begin{split}
\frac{\partial}{\partial r}\bigg |_{p} f &= \frac{\partial}{\partial \overline{r}}\bigg |_{(r,\theta )(p)}[f \circ (r,\theta)^{-1}] =
     \frac{\partial}{\partial \overline{r}}\bigg |_{(r,\theta )(p)}[f \circ \exp_{O} \circ (\overline{r} , \overline{\theta})^{-1} ] \\
\frac{\partial}{\partial \theta }\bigg |_{p} f & =  \frac{\partial}{\partial \overline{\theta}}\bigg |_{(r,\theta )(p)}[f \circ (r,\theta)^{-1}] =
\frac{\partial}{\partial \overline{\theta}}\bigg |_{(r,\theta )(p)}[f \circ \exp_{O} \circ (\overline{r} , \overline{\theta})^{-1} ] ,
\end{split}
\end{equation}
where $f$ can be any smooth function on $\mathbb{H}^2(-a^2)$.
\end{defn}

\noindent
Now, thanks to the concrete expression for the exponential map $\exp_O$ as given in \eqref{concreteexphyper}, we can give concrete expressions for
the natural vector fields $\frac{\partial}{\partial r}$, and $\frac{\partial}{\partial \theta }$ on $\mathbb{H}^2(-a^2)$ as follow.
\begin{equation}\label{concreteexpression}
\begin{split}
\frac{\partial}{\partial r}\bigg |_{\exp_O(rv)} & = \frac{d}{dr}\big ( \exp_O(rv)\big ) = \sinh (ar) (1,0,0) + \cosh (ar)v , \\
\frac{\partial}{\partial \theta } \bigg |_{\exp_O(rv)} & = \frac{\sinh (ar)}{a} v^{\perp} ,
\end{split}
\end{equation}
where $v = (0, v_1 ,v_2) \in T_{O}\mathbb{H}^2(-a^2)$ is any unit vector in $T_{O}\mathbb{H}^2(-a^2)$, with $v^{\perp} = (0, -v_2 ,v_1)$,  and $r > 0$ is any positive number. Now if we consider the smooth vector field $e_2$ on $\mathbb{H}^2(-a^2)$ which is defined by
\begin{equation}
e_2 = \frac{a}{\sinh (ar)} \frac{\partial}{\partial \theta },
\end{equation}
then the restriction of such a smooth vector field $e_2$ \emph{along each geodesic $c_v$} must be a \emph{parallel vector field} along $c_v$, simply because the second expression in \eqref{concreteexpression} informs us that
\begin{equation}
e_2 \big |_{c_{v}(r)} = e_2 \big |_{\exp_O (rv)} = v^{\perp} .
\end{equation}
The above expression further informs us that we must have $\overline{\nabla}_{\dot{c_{v}}} e_2 = 0$, for any direction indicated by the unit vector $v = (0, v_1, v_2) \in T_{O}(\mathbb{H}^2(-a^2))$. Here, the symbol $\overline{\nabla}$ stands for the Levi-Civita connection acting on the space of smooth vector fields on $\mathbb{H}^2(-a^2)$, which is naturally induced by the intrinsic Riemannian geometry of $\mathbb{H}^2(-a^2)$.
So, it turns out that the vector fields $e_1 = \frac{\partial}{\partial r}$, $e_{2} = \frac{a}{\sinh (ar)}\frac{\partial}{\partial \theta }$  constitute a positively oriented orthonormal moving frame $\{e_1 , e_2\}$ of the tangent bundle $T\mathbb{H}^2(-a^2)$ of $\mathbb{H}^2(-a^2)$.\\

\noindent
Due to the fact that we will work with differential operators, such as the Hodge's Laplacian $(-\triangle ) = dd^* + d^*d$ or the exterior differential operator $d$ etc, which can only naturally operate on differential $1$-forms on $\mathbb{H}^2(-a^2)$, we will consider the two associated $1$-forms $e_1^* = g(e_1 ,\cdot )$, and $e_2^* = g(e_2 , \cdot )$ of the vector fields $e_1$, $e_2$ respectively, which together constitute the orthnormal co-frame of the cotangent bundle $T^*\mathbb{H}^2(-a^2)$. Indeed, the associated $1$-forms $e_1^*$, and $e_2^*$ on $\mathbb{H}^2(-a^2)$ are given by
\begin{equation}\label{coframehyper}
\begin{split}
e_1^* & = dr , \\
e_2^* & = \frac{\sinh (ar)}{a} d\theta .
\end{split}
\end{equation}
Then, the volume form $Vol_{\mathbb{H}^2(-a^2)}$ on the hyperbolic space $\mathbb{H}^2(-a^2)$ (see \textbf{Definition \ref{VolumeformDef}}) can be locally expressed by
\begin{equation}
Vol_{\mathbb{H}^2(-a^2)} = e_1^* \wedge e_2^* = \frac{\sinh (ar)}{a} dr \wedge d\theta .
\end{equation}
In dealing with the Hodge Laplacian $(-\triangle) = d d^* + d^*d$, and the co-adjoint operator $d^*$ of $d$, we will encounter the Hodge-star operator $* : C^{\infty}(T^*\mathbb{H}^2(-a^2)) \rightarrow C^{\infty}(T^*\mathbb{H}^2(-a^2))$ which sends $1$-forms into $1$-forms on $\mathbb{H}^2(-a^2)$, and also the Hodge-Star operator $* : C^{\infty}(\wedge^2T^*M) \rightarrow C^{\infty}(M)$ which sends $2$-forms into smooth functions on $\mathbb{H}^2(-a^2)$ (See \textbf{Definition \ref{HodgeDefinition}} for the precise definitions of these Hodge Star operators).\\

\noindent
Indeed, the Hodge-Star operator $* : C^{\infty}(T^*\mathbb{H}^2(-a^2)) \rightarrow C^{\infty}(T^*\mathbb{H}^2(-a^2))$ can be locally expressed in the following way, through the use of the orthonormal co-frame $\{e_1 , e_2\}$ as specified in \eqref{coframehyper}.
\begin{equation}\label{Hodgerotathyp}
\begin{split}
* e_1^* & = e_2^* ,\\
* e_2^* & = -e_1^* .
\end{split}
\end{equation}
Then, in accordance with the tensorial property $*(f \omega ) = f *( \omega )$, with $f$ to be a smooth function and $\omega$ to be a differential form on $\mathbb{H}^2(-a^2)$, of the Hodge-Star operator $*$, it is plain to see, from \eqref{Hodgerotathyp}, that
$* dr = \frac{\sinh (ar)}{a} d \theta $, and that $*d \theta = \frac{a}{\sinh (ar)} *e_2^* = -\frac{a}{\sinh (ar)}dr$. \\

\noindent
On the other hand, the Hodge-Star operator $*: C^{\infty}(\wedge^2T^*\mathbb{H}^2(-a^2)) \rightarrow C^{\infty}(\mathbb{H}^2(-a^2))$ sending smooth $2$-forms into smooth functions on $\mathbb{H}^2(-a^2)$ can be expressed by the following single relation, with $Vol_{\mathbb{H}^2(-a^2)}$ to be the standard volume form on $\mathbb{H}^2(-a^2)$.
\begin{equation}
* Vol_{\mathbb{H}^2(-a^2)} = 1 .
\end{equation}

\section{About stationary Navier-Stokes flows with circular-arc streamlines around an obstacle in $\mathbb{H}^2(-a^2)$ : The proof of Theorem \ref{existenceHyperbolic}  }\label{hyperMAINSECTION}

\noindent
To begin the argument, let $K$ to be a given compact region in $\mathbb{H}^2(-a^2)$ which is entirely contained in $\overline{B_O(\delta )}$, where
$B_O(\delta) = \{ p \in \mathbb{H}^2(-a^2) : d(p, O) < \delta  \}$ is the open geodesic ball centered at $O$, and with radius $\delta > 0$ in the hyperbolic manifold $\mathbb{H}^2(-a^2)$. Let $(r , \theta )$ to be the normal polar coordinate system on $\mathbb{H}^2(-a^2)$ about the reference point $O \in \mathbb{H}^2(-a^2)$ which is specified in \textbf{Definition \ref{polarDefn}}. Suppose further that $\partial K$ contains a circular-arc portion $C_{\delta , \tau } = \{ p \in \mathbb{H}^2(-a^2) : r(p) = d(p , O ) = \delta ,  0 < \theta (p) < \tau   \}$, with some angle $\tau \in (0, 2\pi )$. Let $\frac{\partial}{\partial r}$, $\frac{\partial}{\partial \theta }$ to be the two natural vector fields on $\mathbb{H}^2(-a^2)$ induced by the normal polar coordinate system $(r , \theta )$ (see \textbf{Definition \ref{FRAMEHyperDef}}). Under such setting, we now consider a velocity field of the following form, with $h \in C^{\infty} ( [\delta , \delta + \epsilon_{0} ))$ to be some smooth function defined on an interval $[\delta , \delta + \epsilon_{0} )$ of length $\epsilon_{0} > 0$.
\begin{equation}\label{circulararchyperbolic}
u = -h(r- \delta ) e_{2} = - h(r-\delta ) \frac{a}{ \sinh (ar )} \frac{\partial}{\partial \theta }
\end{equation}
Recall that $e_1 = \frac{\partial}{\partial r}$, $e_2 = \frac{a}{ \sinh (ar )} \frac{\partial}{\partial \theta }$ together constitute a positively oriented orthonormal moving frame $\{e_1 , e_2\}$ on $\mathbb{H}^2(-a^2)$, whose orthonormal coframe $\{ e_1^* , e_2^* \}$ is constituted by
the differential $1$-forms $e_1^* =dr$, and $e_2^* = \frac{\sinh (ar)}{a} d\theta$. So, the associated $1$-form $u^*$ of the velocity field $u$ in
\eqref{circulararchyperbolic} is just given by
\begin{equation}\label{associated1fromhyper}
u^* = -h(r- \delta ) e_{2}^* = - h(r-\delta ) \frac{\sinh (ar)}{a} d\theta .
\end{equation}
Notice that under this setting, both the velocity field $u$ as specified in \eqref{circulararchyperbolic} and its associated $1$-form $u^*$ are defined on the sector-shaped open region $R_{\delta , \tau , \epsilon } = \{ p \in \mathbb{H}^2(-a^2) : \delta < d(p,O) < \delta + \epsilon_0 ,  0 < \theta (p) < \tau \}$ (The same open region as the one specified in \eqref{RegionHYPER}) whose boundary shares the same circular-arc boundary portion $C_{\delta , \tau }$ with $\partial K$.

In accordance with expression \eqref{associated1fromhyper} for $u^*$, we will compute $(-\triangle ) u^*$, and $\overline{\nabla}_{u}u$ step by step, just as what we did in dealing with the spherical case.\\

\noindent{\bf Step 1 : Checking the divergence free property of the velocity field $u$ as specified in \eqref{circulararchyperbolic}.   }
For $u$ as given in \eqref{circulararchyperbolic}, in order to verify the divergence property $d^* u^* = 0$ on the sector-shaped open region $R_{\delta , \tau , \epsilon }$ of $\mathbb{H}^2(-a^2)$ as specified in \eqref{RegionHYPER}, we just carry out the following straightforward computation, which is \emph{formally} identical to those computations being done in \eqref{divfreesphere}.
\begin{equation}
d^* u^* = -*d* \big \{ - h(r-\delta ) e_{2}^* \big \} = -*d \{ h(r-\delta ) dr\} = 0 .
\end{equation}
\noindent{\bf Step $2$ : The computation of $(-\triangle ) u^*$, for $u$ to be given in \eqref{circulararchyperbolic}.}
Since we have $d^* u = 0$ on the open region $R_{\delta , \tau , \epsilon }$ of $\mathbb{H}^2(-a^2)$ as specified in \eqref{RegionHYPER}, it follows that $(-\triangle )u^* = (dd^* + d^* d)u^* = d^* du^*$. So, as in the spherical case, we first compute $du^*$, for the velocity field $u$ given in \eqref{circulararchyperbolic}, as follow.
\begin{equation}\label{vorticityDischyp}
\begin{split}
du^* & = -\frac{\partial}{\partial r} \bigg \{ \frac{\sinh (ar)}{a} h(r-\delta )  \bigg \} dr \wedge d\theta \\
& = - \frac{a}{\sinh (ar)} \frac{\partial}{\partial r} \bigg \{ \frac{\sinh (ar)}{a} h(r-\delta )  \bigg \} Vol_{\mathbb{H}^2(-a^2)} \\
& = - \bigg \{  h'(r-\delta ) +   \frac{a \cosh (ar)}{\sinh (ar)} h(r-\delta )     \bigg \} Vol_{\mathbb{H}^2(-a^2)} ,
\end{split}
\end{equation}
where we recall that $Vol_{\mathbb{H}^2(-a^2)} = e_1^* \wedge e_2^* = \frac{\sinh (ar)}{a} dr \wedge d\theta$ is the volume form on $\mathbb{H}^2(-a^2)$. For $u$ as given in \eqref{circulararchyperbolic}, we now compute $(-\triangle )u^* = d^* du^*$ as follow.
\begin{equation}\label{Laplacianhyper}
\begin{split}
(-\triangle )u^* & = *d* \bigg \{ h'(r-\delta ) +   \frac{a \cosh (ar)}{\sinh (ar)} h(r-\delta )  \bigg \}  Vol_{\mathbb{H}^2(-a^2)} \\
& = * d \bigg \{ h'(r-\delta ) +   \frac{a \cosh (ar)}{\sinh (ar)} h(r-\delta )  \bigg \} \\
& = \bigg \{ h''(r- \delta ) +  \frac{a \cosh (ar)}{\sinh (ar)} h'(r-\delta ) +   \frac{\partial}{\partial r}\bigg (\frac{a \cosh (ar)}{\sinh (ar)} \bigg ) h(r-\delta )   \bigg \} * dr \\
& = \bigg \{ \frac{\sinh (ar)}{a} h''(r-\delta ) + \cosh (ar) h'(r-\delta ) - \frac{a}{\sinh (ar)} h(r-\delta )     \bigg \} d\theta ,
\end{split}
\end{equation}
in which the second equal sign holds due to the fact that $*  Vol_{\mathbb{H}^2(-a^2)} = 1$, and the last equal sign holds since
$* dr = e_2^* = \frac{\sinh (ar)}{a} d\theta$.
To prepare for the proof of Theorem \ref{existenceHyperbolic}, we will need the expression of $d \big \{ (-\triangle )u^* -2Ric (u^*)  \big \}$, for $u$ as specified in \eqref{circulararchyperbolic}. Since $Ric (X^*) = -a^2 X^*$ always holds for any smooth vector field $X$ on $\mathbb{H}^2(-a^2)$, it follows from \eqref{Laplacianhyper} and a direct computation that we have the following expression of $d \big \{ (-\triangle )u^* -2Ric (u^*)  \big \}$, for $u$ as specified in \eqref{circulararchyperbolic}.
\begin{equation}\label{VorticityhyperDisc}
\begin{split}
d \big \{ (-\triangle )u^* -2Ric (u^*)  \big \} & = \bigg \{ \frac{\sinh (ar)}{a}h'''(r-\delta ) + 2 \cosh (ar) h''(r-\delta ) \\
& -a\bigg ( \sinh (ar) + \frac{1}{\sinh (ar)}   \bigg ) h'(r-\delta ) \\
& + a^2\cosh (ar) \bigg ( \frac{1}{\sinh^2(ar)} -2 \bigg ) h(r-\delta ) \bigg \} dr \wedge d \theta .
\end{split}
\end{equation}

\noindent{\bf Step 3 :  The computation of the nonlinear convection term $\overline{\nabla}_uu$ for $u$ as given in \eqref{circulararchyperbolic} . }

\noindent
To compute $\overline{\nabla}_uu$ for $u$ as given in \eqref{circulararchyperbolic}, we express $\overline{\nabla}_uu$ as follow.
\begin{equation}\label{convectionhyper}
\overline{\nabla}_uu = A \frac{\partial}{\partial r} + B \frac{\partial}{\partial \theta } ,
\end{equation}
with $A$, $B$ to be the two component functions on $\mathbb{H}^2(-a^2)$. To compute $B$, we take the inner product with $\frac{\partial}{\partial \theta }$ on both sides of \eqref{convectionhyper}, and carry out the following computation.
\begin{equation}
\begin{split}
B \bigg (\frac{\sinh^2 (ar)}{a^2}\bigg ) & = g( \overline{\nabla}_uu , \frac{\partial}{\partial \theta }  ) \\
& = -\bigg ( \frac{\sinh (ar)}{a} \bigg ) \frac{1}{h(r-\delta )} g( \overline{\nabla}_uu , u  ) \\
& = -\bigg ( \frac{\sinh (ar)}{a} \bigg ) \frac{1}{2 h(r-\delta )} u \big (  |u|^2  \big ) \\
& = \frac{1}{2} \frac{\partial}{\partial \theta } \big ( (h(r-\delta ))^2 \big ) = 0 .
\end{split}
\end{equation}
In the above computation, the third equal sign follows from property \textbf{(1)} of the Levi-Civita connection $\overline{\nabla}$ as stated in \textbf{Definition \ref{LeviCivitadefinition}}. The symbol $u \big (  |u|^2  \big )$ stands for the derivative of the function $|u|^2$ along the direction of the vector field $u$.\\

\noindent
Next, we compute the component $A$ which appears in \eqref{convectionhyper}, by taking inner product with $\frac{\partial}{\partial r}$ on both sides of \eqref{convectionhyper} and  we get

\begin{equation}
\begin{split}
A & = g( \overline{\nabla}_uu , \frac{\partial}{\partial r} ) \\
& = -h(r-\delta ) \frac{a}{\sinh(ar)} g ( \overline{\nabla}_{\frac{\partial}{\partial \theta } }u ,  \frac{\partial}{\partial r}    ) \\
& = h(r-\delta ) \frac{a}{\sinh(ar)} g ( u , \overline{\nabla}_{\frac{\partial}{\partial \theta } } \frac{\partial}{\partial r} ) \\
& = h(r-\delta ) \frac{a}{\sinh(ar)} g ( u , \overline{\nabla}_{\frac{\partial}{\partial r}}\frac{\partial}{\partial \theta }   ) .
\end{split}
\end{equation}
In the above computation, the second equal sign follows directly from the tensorial property (property \textbf{(4)} in \textbf{Definition \ref{LeviCivitadefinition}}) of $\overline{\nabla}$. The third equal sign follows from a direct application of property \textbf{(1)} in \textbf{Definition \ref{LeviCivitadefinition}} of $\overline{\nabla}$. While the last equal sign follows from the torsion free property (property \textbf{(2)} of \textbf{Definition \ref{LeviCivitadefinition}} ) of $\overline{\nabla}$.\\

\noindent
Here, recall that $e_2 = \frac{a}{\sinh(ar)} \frac{\partial}{\partial \theta }$, once being restricted on each geodesic ray starting from the base point $O \in \mathbb{H}^2(-a^2)$ of the normal polar coordinate system $(r,\theta )$, is \emph{parallel} along that geodesic ray. This simply means that we have $\overline{\nabla}_{\frac{\partial}{\partial r}} e_2 = 0$. Hence, we can carry out the following computation in accordance with the Leibniz's rule of the connection $\overline{\nabla}$ (property \textbf{(3)} of \textbf{Definition \ref{LeviCivitadefinition}} ).
\begin{equation}
\overline{\nabla}_{\frac{\partial}{\partial r}} \frac{\partial}{\partial \theta } = \frac{\partial}{\partial r} \bigg ( \frac{\sinh (ar)}{a}\bigg ) e_2
= \cosh (ar) e_2 .
\end{equation}
Hence, it follows that the component function $A$ which appears in \eqref{convectionhyper} is given by
\begin{equation}
A = h(r-\delta ) \frac{a}{\sinh(ar)} g(u , \cosh (ar) e_2  ) = -a \bigg ( \frac{\cosh (ar)}{\sinh (ar)} \bigg ) \big ( h(r-\delta )\big )^2.
\end{equation}
So, it follows that for $u$ as given in \eqref{circulararchyperbolic}, the term $\overline{\nabla}_uu$ is given by
\begin{equation}
\overline{\nabla}_uu =  -a \bigg ( \frac{\cosh (ar)}{\sinh (ar)} \bigg ) \big ( h(r-\delta )\big )^2 \frac{\partial}{\partial r} ,
\end{equation}
whose associated $1$-form $[\overline{\nabla}_uu]^*$ is given by
\begin{equation}\label{nonlinearconvectionhyper}
[\overline{\nabla}_uu]^* =  -a \bigg ( \frac{\cosh (ar)}{\sinh (ar)} \bigg ) \big ( h(r-\delta )\big )^2 dr .
\end{equation}
So, by taking the operator $d$ on both sides of \eqref{nonlinearconvectionhyper}, it follows that the following relation holds on the sector-shaped region $R_{\delta , \tau \epsilon_0}$ of $\mathbb{H}^2(-a^2)$ as specified in \eqref{RegionHYPER}, for $u$ to be given  by \eqref{circulararchyperbolic}.
\begin{equation}\label{dconzerohyper}
d [\overline{\nabla}_uu]^* = 0 .
\end{equation}

\noindent{\bf Step 4 : The proof of \textbf{Assertion I} in Theorem \ref{existenceHyperbolic}}

\noindent
Here, we will give a simple proof of \textbf{Assertion I} in Theorem \ref{existenceHyperbolic}, which states that for any quadratic profile $h(\lambda ) = \alpha_1 \lambda - \frac{\alpha_2}{2} \lambda^2$, with prescribed constants $\alpha_1 > 0$, and $\alpha_2 > 0$, the velocity field $u$ as specified in \eqref{circulararchyperbolic} does not satisfies equation \eqref{NavierStokeshyperbolic} on the sector shaped region $R_{\delta , \tau , \epsilon_0 }$ of $\mathbb{H}^2(-a^2)$ as specified in \eqref{RegionHYPER} , regardless of how small $\epsilon_0 > 0$ is.
Now, assume towards contradiction that for a certain choice of constants $\alpha_1 >0$, $\alpha_2 > 0$, the velocity field $u$ as given in \eqref{circulararchyperbolic} does satisfy equation \eqref{NavierStokeshyperbolic} on the sector-shaped region $R_{\delta , \tau , \epsilon_0 }$ of $\mathbb{H}^2(-a^2)$ as given in \eqref{RegionHYPER}, for some $\epsilon_0 > 0$. Then, for such a $u$ as given in \eqref{circulararchyperbolic}, we take the operator $d$ on both sides of the main equation in \eqref{NavierStokeshyperbolic}, and deduce the following vorticity equation from \eqref{vorticityDischyp} and \eqref{dconzerohyper}.

\begin{equation}\label{VorticityequationDischyp}
\begin{split}
0 & = d \bigg \{ \nu \big ( (-\triangle )u^* - Ric(u^*) \big ) + [\overline{\nabla}_uu]^* + dP \bigg \} \\
& =  \nu d\bigg \{ \big ( (-\triangle )u^* - Ric(u^*) \big ) \bigg \} \\
& = \nu \bigg \{ \frac{\sinh (ar)}{a}h'''(r-\delta ) + 2 \cosh (ar) h''(r-\delta ) \\
& -a\bigg ( \sinh (ar) + \frac{1}{\sinh (ar)}   \bigg ) h'(r-\delta ) \\
& + a^2\cosh (ar) \bigg ( \frac{1}{\sinh^2(ar)} -2 \bigg ) h(r-\delta ) \bigg \} dr \wedge d \theta .
\end{split}
\end{equation}
For the quadratic profile $h(\lambda ) = \alpha_1 \lambda -\frac{\alpha_2}{2} \lambda^2$, the vorticity equation \eqref{VorticityequationDischyp} reduces to the following form
\begin{equation}\label{Vortsimplehyper}
0 = G_{\alpha_1 , \alpha_2 , \delta} (r) dr \wedge d \theta ,
\end{equation}
where the function $G_{\alpha_1 , \alpha_2 , \delta} (r)$ is given by the following expression.
\begin{equation}\label{TrivialhyperDisc}
\begin{split}
G_{\alpha_1 , \alpha_2 , \delta} (r) & = -2\alpha_2 \cosh (ar) -a \bigg ( \sinh (ar) + \frac{1}{\sinh(ar)} \bigg ) (\alpha_1 - \alpha_2 \lambda ) \\
&  + a^2 \cosh (ar) \bigg ( \frac{1}{\sinh^2(ar)} -2 \bigg ) \big ( \alpha_1 \lambda -\frac{\alpha_2}{2} \lambda^2 \big ) .
\end{split}
\end{equation}
So, if $u = - h(r-\delta ) \frac{a}{\sinh (ar)} \frac{\partial}{\partial \theta }$ does satisfies equation \eqref{NavierStokeshyperbolic} on the sector-shaped region $R_{\delta , \tau , \epsilon_0 }$ of $\mathbb{H}^2(-a^2)$ as specified in \eqref{RegionHYPER}, then, it must follow, in accordance with \eqref{Vortsimplehyper}, that the function $G_{\alpha_1 , \alpha_2 , \delta} (r)$ must totally vanish on the interval $(\delta , \delta + \epsilon_0 )$. However, we observe that
\begin{equation}\label{Gnegative}
G_{\alpha_1 , \alpha_2 , \delta} (\delta ) = -2\alpha_2 \cosh (a \delta ) - a \alpha_1 \bigg ( \sinh(a\delta ) + \frac{1}{\sinh (a\delta )} \bigg ) < 0 .
\end{equation}
Since $G_{\alpha_1 , \alpha_2 , \delta} (\delta )$ is continuous on $[\delta ,\infty )$, \eqref{Gnegative} immediately implies that $G_{\alpha_1 , \alpha_2 , \delta}< 0 $ holds on some interval $[\delta , \delta + \epsilon_1)$, for some $0 < \epsilon_1 < \epsilon_0$, which directly violates \eqref{Vortsimplehyper}. A contradiction has been achieved, and we are done in proving \textbf{Assertion I} of Theorem.\\

\noindent{\bf Step 5 : The proof of Assertion II in Theorem \ref{existenceHyperbolic}}

\noindent
Here, for any prescribed positive constants $\alpha_1 > 0$, and $\alpha_2>0$, we consider the velocity field $u = - Y(r) \frac{a}{\sinh (ar)} \frac{\partial}{\partial \theta }$ on the region $\Omega_{\delta , \tau }$ as given in \eqref{OMEGAhyper} of Theorem \ref{existenceHyperbolic}, with $Y \in C^{\infty } ([\delta ,\infty ))$ to be a smooth function on $[\delta ,\infty )$ satisfying $Y(\delta ) =0$, $Y'(\delta ) = \alpha_1$, and $Y''(\delta ) = -\alpha_2$. Since the sector-shaped region $\Omega_{\delta , \tau }$ as given in \eqref{OMEGAhyper} is simply connected in $\mathbb{H}^2(-a^2)$, we know that such
$u = - Y(r) \frac{a}{\sinh (ar)} \frac{\partial}{\partial \theta }$ will satisfy equation \eqref{NavierStokeshyperbolic} with some globally defined smooth pressure $P$ on $\Omega_{\delta , \tau }$ \emph{if and only if} the vorticity equation \eqref{VorticityequationDischyp} holds on the simply-connected region $\Omega_{\delta , \tau }$ of $\mathbb{H}^2(-a^2)$. However, saying that the vorticity equation \eqref{VorticityequationDischyp} holds on the simply connected region $\Omega_{\delta , \tau }$ as specified in \eqref{OMEGAhyper} of Theorem \ref{existenceHyperbolic} is equivalent to saying that the smooth function $Y \in C^{\infty } ([\delta ,\infty ))$ is a solution to the following third order O.D.E. on $[\delta ,\infty )$ with initial values $Y (\delta ) = 0$, $Y'(\delta )= \alpha_1$, and $Y''(\delta ) =-\alpha_2$.

\begin{equation}\label{ODEhyperbolicDisc}
\begin{split}
0 & = \frac{\sinh (ar)}{a}Y'''(r) + 2 \cosh (ar) Y''(r) -a\bigg ( \sinh (ar) + \frac{1}{\sinh (ar)}   \bigg ) Y'(r) \\
& + a^2\cosh (ar) \bigg ( \frac{1}{\sinh^2(ar)} -2 \bigg ) Y(r).
\end{split}
\end{equation}
However, in accordance with the basic existence theorem in the theory of linear O.D.E.(i.e. Theorem \ref{ODEtheorem}) , we know that there exists a unique smooth solution $Y \in C^{\infty } ([\delta ,\infty ))$ to \eqref{ODEhyperbolicDisc} satisfying initial values $Y (\delta ) = 0$, $Y'(\delta )= \alpha_1$, and $Y''(\delta ) =-\alpha_2$. Moreover, since the coefficient functions involved in \eqref{ODEhyperbolicDisc} are all real analytic on $[\delta , \infty )$, it follows that such a unique solution $Y \in C^{\infty } ([\delta ,\infty ))$ must also be \textbf{real-analytic} on $[\delta , \infty )$. So, according to these observations, we can now deduce that, for any prescribed positive constants $\alpha_1 > 0$, and $\alpha_2 >0$, there exists a unique smooth function $Y \in C^{\infty } ([\delta ,\infty ))$ with $Y (\delta ) = 0$, $Y'(\delta )= \alpha_1$, and $Y''(\delta ) =-\alpha_2$ such that $u = - Y(r) \frac{a}{\sinh (ar)} \frac{\partial}{\partial \theta }$ is a solution to equation \eqref{NavierStokeshyperbolic} on the simply connected region $\Omega_{\delta , \tau }$ as specified in \eqref{OMEGAhyper}. Moreover, such a unique smooth function $Y \in C^{\infty } ([\delta ,\infty ))$ is further known to be real-analytic on $[\delta , \infty )$. So, we are done in proving \textbf{Assertion II} of Theorem \ref{existenceHyperbolic}.

\section{The ''Cartesian coordinate system'' on the hyperbolic space $\mathbb{H}^2(-a^2)$}\label{Cartesianhyperbolicsection}

\noindent
The purpose of this section is to introduce a natural coordinate system $\Phi : \mathbb{R}^2 \rightarrow \mathbb{H}^2(-a^2)$ on the hyperbolic space $\mathbb{H}^2(-a^2)$, with $\mathbb{R}^2 = \{(\tau , s) : \tau , s \in (-\infty , \infty )\}$ to be the parameter space being used to parameterize the manifold $\mathbb{H}^2(-a^2)$. Since this natural coordinate system $\Phi(\tau , s )$ which we proceed to construct here is the closest possible analog to the standard cartesian coordinate system of the Euclidean space $\mathbb{R}^2$, we will just call $\Phi(\tau , s )$ which we construct below to be a ''Cartesian coordinate system'' introduced on $\mathbb{H}^2(-a^2)$. Again, such a ''Cartesian coordinate system'' $\Phi (\tau , s)$ on $\mathbb{H}^2(-a^2)$ which we are going to describe is also another piece of standard knowledge in Riemannian Geometry. Indeed, our discussions in this section can be viewed as a special case of the well-known procedure of constructing Jacobi-fields along geodesics on a general Riemannian manifold, and we refer the interested readers to pages 185-193 of the textbook by Jost \cite{Jost}.\\

\noindent
To begin the construction, let $O \in \mathbb{H}^2(-a^2)$ to be any selected point in the hyperbolic space $\mathbb{H}^2(-a^2)$ of constant negative sectional curvature $-a^2$. Recall that $\mathbb{H}^2(-a^2)$ (just as $\mathbb{R}^2$, or $S^2(a^2)$ ) is a symmetric space in that the \emph{its geometric structure looks exactly the same at any selected reference point, up to isometries preserving the Riemannian metric on $\mathbb{H}^2(-a^2)$}. So, we can just select any reference point $O \in \mathbb{H}^2(-a^2)$, which will play the role of the origin of the Cartesian coordinate system as introduced below.\\
With such a reference point $O$ in $\mathbb{H}^2(-a^2)$ to be chosen and fixed, consider a geodesic $c : (-\infty , \infty ) \rightarrow  \mathbb{H}^2(-a^2)$ which passes through $O$ in that $c(0)=O$, and which travels towards the East-direction as the parameter $\tau \in (-\infty , \infty )$ increases. Recall that a geodesic on a Riemannian manifold is really a straight line with respect to the Riemannian structure of that manifold.   Such a geodesic $c(\tau )$ which we choose will play the role of the \emph{$x$-axis} for our Cartesian coordinate system on the hyperbolic space $\mathbb{H}^2(-a^2)$. In order to specify the appropriate $y$-axis, we will regard $\mathbb{H}^2(-a^2)$ to be an oriented manifold with the orientation compatible with the anti-clockwise rotation. Then, choose $w \in T_{O}\mathbb{H}^2(-a^2)$ to be the unit vector which, together with $\dot{c}(0) = \frac{d}{d\tau }c|_{\tau =0} \in T_{O}\mathbb{H}^2(-a^2)$, constitute a positively oriented \emph{orthonormal} basis $\{\dot{c}(0) , w \}$ of $T_{O}\mathbb{H}^2(-a^2)$ compatible with the anti-clockwise rotation on $T_{O}\mathbb{H}^2(-a^2)$. Then, we just consider the geodesic $\gamma :(-\infty , \infty ) \rightarrow \mathbb{H}^2(-a^2)$ which satisfies the properties that $\gamma (0) = O $ and that $\dot{\gamma}(0) = \frac{d}{ds}\gamma|_{s= 0}=w$. So, $\gamma$ will be a straight line (i.e. geodesic) which passes through $O$ and which travels towards the North-direction as the parameter $s \in (-\infty , \infty )$ increases. So, the geodesic $\gamma$, which intersects the geodesic $c$ at the reference point $O$ in an orthogonal manner, will play the proper role of the $y$-axis of our Cartesian coordinate system on $\mathbb{H}^2(-a^2)$. For a technical purpose, consider $V(s)$ to be the \emph{parallel vector field} along the geodesic $\gamma$ which satisfies $V(0) = \dot{c}(0)$ in the sense as specified in the following Definition

\begin{defn}
Let $\gamma : (a, b) \rightarrow M$ to be a geodesic on a $N$-dimensional Riemannian manifold $M$. A parallel vector field $V(s)$ along $\gamma$ is a smooth map $s \in (a,b) \to V(s)\in T_{\gamma (s)}M $ for which the property $\overline{\nabla}_{\dot{\gamma}}V = 0$ holds on $(a,b)$. Here,   $\overline{\nabla}$ is the Levi-Civita connection (covariant derivative) acting on the space of smooth vector fields on $M$.
\end{defn}

Now, with such a parallel vector field $V(s)$ along the geodesic $\gamma$ with $V(0) = \dot{c}(0)$, we consider, for each real value $s\in \mathbb{R}$, the geodesic $c_{s}: (-\infty , \infty ) \rightarrow \mathbb{H}^2(-a^2) $ which passes through the point $\gamma (s)$ in that $c_{s}(0) = \gamma (s)$, and which satisfies the condition $\dot{c_{s}}(0) = \frac{d}{d\tau }c_{s}|_{\tau =0} = V(s) \in T_{\gamma(s)}\mathbb{H}^2(-a^2)$. In accordance with the definition of the exponential map about a reference point on a Riemannian manifold, we can express the geodesic $c_{s}$ in terms of the exponential map $\exp_{\gamma (s)} : T_{\gamma (s)} \mathbb{H}^2(-a^2) \rightarrow \mathbb{H}^2(-a^2)$ as follow, where $\tau \in (-\infty , \infty )$ is the parameter of the geodesic $c_{s}$,
\begin{equation}
c_{s}(\tau ) = \exp_{\gamma (s)}(\tau V(s)) .
\end{equation}
We can now define the smooth bijective map $\Phi : \mathbb{R}^2 \rightarrow \mathbb{H}^2(-a^2)$ in accordance with the following rule, where $\tau ,s \in \mathbb{R}$ are arbitrary real parameters.
\begin{equation}\label{coordinatesystemhyper}
\Phi (\tau , s) = \exp_{\gamma (s)}(\tau V(s)) = c_{s}(\tau ) .
\end{equation}
Indeed, the smooth map $\Phi$, which maps the parameter space $\mathbb{R}^2$ bijectively onto $\mathbb{H}^2(-a^2)$ with smooth inverse $\Phi^{-1}$, is exactly the coordinate system which we introduced on the hyperbolic manifold $\mathbb{H}^2(-a^2)$. Next, by means of this natural coordinate system $\Phi( \tau ,s)$, we will define two natural vector fields $\frac{\partial}{\partial \tau }$, and $\frac{\partial}{\partial s}$ \textbf{on the hyperbolic space} $\mathbb{H}^2(-a^2)$ as follow.

\begin{defn}
$($Natural vector fields $\frac{\partial}{\partial \tau }$, and $\frac{\partial}{\partial s}$ on $\mathbb{H}^2(-a^2)$ via the coordinate system $\Phi(\tau ,s)$ as given in \eqref{coordinatesystemhyper}$)$.\\
For any point $p \in \mathbb{H}^2(-a^2)$, the vectors $\frac{\partial}{\partial \tau }|_{p}$, $\frac{\partial}{\partial s}|_{p}$ in the tangent space $T_{p}\mathbb{H}^2(-a^2)$ of $\mathbb{H}^2(-a^2)$ at $p$ are characterized (as linear derivations acting on smooth functions) by the following rules.
\begin{equation}\label{naturalrelationhyp}
\begin{split}
\frac{\partial}{\partial \tau } \bigg |_{p} f &= \frac{\partial}{\partial \tau }(f \circ \Phi ) \bigg |_{\Phi^{-1}(p)} , \\
\frac{\partial}{\partial s} \bigg |_{p} f & = \frac{\partial}{\partial s}(f \circ \Phi ) \bigg |_{\Phi^{-1}(p)} ,
\end{split}
\end{equation}
where $f \in C^{\infty}(\mathbb{H}^{2}(-a^2))$ is any smooth function on $\mathbb{H}^{2}(-a^2)$. Here, we remark that the same symbol $\frac{\partial}{\partial \tau}$ means totally \emph{different} things on the two sides of relation \eqref{naturalrelationhyp}. The symbol $\frac{\partial}{\partial \tau } |_{p}$ on the left stands for the vector in $T_{p}\mathbb{H}^2(-a^2)$ which is to be defined through the right hand side. While, the symbol $\frac{\partial}{\partial \tau }$ appearing on the righthand-side is just the ordinary partial derivative acting on the Euclidean space $\mathbb{R}^2$ at the point $\Phi^{-1}(p) \in \mathbb{R}^2$. The same remark also applies to the symbol $\frac{\partial}{\partial s}$ appearing in the second line of \eqref{naturalrelationhyp}.
\end{defn}

\noindent
In accordance with the above rigorous definition for the vector fields $\frac{\partial}{\partial \tau }$, and $\frac{\partial}{\partial s}$ on  $\mathbb{H}^2(-a^2)$, for each pair of parameters $(\tau ,s) \in \mathbb{R}^2$, we can \emph{think of} the two vectors $\frac{\partial}{\partial \tau }|_{\Phi (\tau ,s)}$ and $\frac{\partial}{\partial s}|_{\Phi (\tau ,s)}$ in $T_{\Phi (\tau ,s )} \mathbb{H}^2(-a^2)$ in the following intuitive manner.
\begin{equation}\label{intuitivehyp}
\begin{split}
\frac{\partial}{\partial \tau } \bigg |_{\Phi (\tau ,s)} & = \partial_{\tau}\{ \exp_{\gamma (s)}(\tau V(s)) \} = \partial_{\tau } \{ c_{s}(\tau )  \} , \\
\frac{\partial}{\partial s} \bigg |_{\Phi (\tau ,s)} & = \partial_{s}\{ \exp_{\gamma (s)}(\tau V(s)) \}
\end{split}
\end{equation}
As a result of relation \eqref{intuitivehyp}, it follows that for any point $p \in \mathbb{H}^{2}(-a^2)$, which is parameterized by the pair $(\tau ,s) \in \mathbb{R}^2$ $($i.e. $(\tau ,s)= \Phi^{-1}(p)$ $)$, we will have the following relation
\begin{equation}
\frac{\partial}{\partial \tau } \bigg |_{p} = \dot{c_{s}}(\tau ) .
\end{equation}
Notice that $\dot{c_{s}}(\tau )$ itself is a tangential parallel vector field along the geodesic $c_{s}$. So, it turns out that $\frac{\partial}{\partial \tau }$ must be of unit length. That is, we have $\|\frac{\partial}{\partial \tau } \| = 1$ holds everywhere on $\mathbb{H}^{2}(-a^2)$.\\

\noindent
Next, in accordance with basic knowledge in Riemannian geometry, the vector field $\frac{\partial}{\partial s}$, when being restricted to each geodesic $c_{s}$, is the \emph{uniquely determined} Jacobi field $W^{(s)}(\tau )$
along $c_{s}$ satisfying the initial conditions $W^{(s)}(0) = \dot{\gamma}(s)$, and $\overline{\nabla}_{\dot{c_{s}}}W^{(s)}|_{\tau = 0} = 0 $, in the sense of the following definition.

\begin{defn}\label{Jacobifield}
On a $N$-dimensional Riemannian manifold $M$, let $c : (a,b) \rightarrow M$ to be a geodesic. A vector field $X(\tau )$ along the geodesic $c(\tau )$ is called a Jacobi field along $c$ if $V$ satisfies the following Jacobi-field equation on $\tau \in (-\infty , \infty )$.
\begin{equation}
\overline{\nabla}_{\dot{c}}\overline{\nabla}_{\dot{c}} X + R(X,\dot{c})\dot{c} = 0 ,
\end{equation}
in which the symbol $\overline{\nabla}$ is again the Levi-Civita connection acting on the space of all smooth vector fields on $M$, and the symbol $R(\cdot ,\cdot )$ stands for the Riemannian curvature tensor which is defined in the following relation.
\begin{equation}
R(X,Y)Z = \overline{\nabla}_{X} \overline{\nabla}_{Y} Z - \overline{\nabla}_{Y} \overline{\nabla}_{X} Z - \overline{\nabla}_{[X,Y]}Z ,
\end{equation}
with $X$, $Y$, $Z$ to be smooth vector fields on $M$, and $[X, Y]$ is the smooth vector field given by $[X,Y] = XY-YX$.
\end{defn}

\noindent
The above definition of Jacobi-fields along geodesic on a given manifold may looks strange to those readers who are not familiar with Riemannian geometry. Giving a detailed discussion of the precise geometric meaning of the concept of Jacobi fields on a Riemannian manifold is out the scope of this paper. Here, we just simply mention that, on an intuitive level, the magnitude of a non-tangential Jacobi field along a geodesic on a Riemannian manifold encodes the growth rate of the spatial structure of the Riemannian manifold in the far range. But, fortunately, the very special symmetric structure of $\mathbb{H}^2(-a^2)$ ensures that the Riemannian curvature tensor $R(\cdot, \cdot )$ on $\mathbb{H}^2(-a^2)$ satisfies the following simple relation.
\begin{equation}\label{simplecurvaturetensor}
R(X, \dot{c}) \dot{c} = -a^2 X ,
\end{equation}
where $c :(-\infty , \infty ) \rightarrow \mathbb{H}^2(-a^2)$ is a geodesic on $\mathbb{H}^2(-a^2)$, and $X$ is a smooth vector field along the geodesic $c$. So, in the case of the hyperbolic manifold $\mathbb{H}^2(-a^2)$, we can just take relation \eqref{simplecurvaturetensor} for granted and hence the Jacobi-field equation as in Definition \eqref{Jacobifield} will reduce down to the following one in the case of a smooth vector field $X$ defined along a geodesic $c$ on $\mathbb{H}^2(-a^2)$.
\begin{equation}\label{Jacobifieldshyper}
\overline{\nabla}_{\dot{c}}\overline{\nabla}_{\dot{c}} X - a^2 X = 0.
\end{equation}
With the Jacobi-field equation as in \eqref{Jacobifieldshyper}, we can give a geometric description for the vector field $\frac{\partial}{\partial s}$  as follow. Here, For each fixed $s \in (-\infty , \infty )$, let $e_2^{(s)}(\tau )$ to be the parallel vector field along the geodesic $c_{s}$ which satisfies the initial condition $e_2^{(s)}(0) = \frac{d}{ds}\gamma|_{s} = \dot{\gamma}(s) $. Recall that, from our construction, we have
\begin{equation}
\dot{c_{s}}(0) = \frac{d}{d\tau }c_{s}\bigg |_{\tau =0 } = V(s) ,
\end{equation}
Where $V(s)$ is the parallel unit vector field along the geodesic $\gamma$ which satisfies $V(0) = \frac{d}{d\tau }c|_{\tau = 0}= \dot{c}(0)$. Since $V(s)$, as a parallel vector field along a geodesic, always preserves its inscribed angle with $\dot{\gamma}(s)$, $\dot{c_{s}}(0) = V(s)$ must be  \emph{orthonormal} to $\dot{\gamma}(s) = e_2^{(s)}(0)$ in the tangent space $T_{\gamma (s)}\mathbb{H}^2(-a^2)$. That is, we have $e_2^{(s)}(0)\perp \dot{c_{s}}(0)$. Then, $e_2^{(s)}(\tau )$, as a parallel vector field along $c_{s}$, must preserves the its inscribed right-angle with $\dot{c_{s}}(\tau )$. That is, the parallel vector field $e_2^{(s)}(\tau )$ along $c_{s}$ is everywhere orthogonal to the geodesic $c_s$ itself. Now, since the coordinate system $\Phi : \mathbb{R}^2 \rightarrow \mathbb{H}^2(-a^2)$, as specified in \eqref{coordinatesystemhyper}, maps $\mathbb{R}^2$ bijectively onto $\mathbb{H}^2(-a^2)$, we know that each point $p \in \mathbb{H}^2(-a^2)$ has to be passed through by exactly one geodesic $c_{s}$. So, we can define a global smooth vector field $e_{2}$ on $\mathbb{H}^2(-a^2)$ by the following relation, where $\tau$, $s$ are arbitrary parameters.
\begin{equation}\label{e2hyper}
e_2(\Phi (\tau , s ) ) = e_2^{(s)}(\tau ) .
\end{equation}
Again, the vector field $e_2$ as defined above is everywhere orthonormal to $\frac{\partial}{\partial \tau }$ (just recall that $\frac{\partial}{\partial \tau }|_{\Phi(\tau ,s )} = \dot{c_s}(\tau )$ ). Now, we claim that the vector field $\frac{\partial}{\partial s}$ is related to $e_2$ through the following relation.
\begin{equation}\label{coshgrowth}
\frac{\partial}{\partial s} \bigg |_{\Phi(\tau , s)} = \cosh (ar) e_2(\Phi (\tau ,s ))
\end{equation}
To justify the above relation, we just have to recall that, for each $s \in (-\infty , \infty )$, the restriction of $\frac{\partial}{\partial s}$ on the geodesic $c_{s}$ is known to be the unique Jacobi field $W(\tau )$ along $c_{s}$ which satisfies the following initial conditions
\begin{equation}\label{initialvalues}
\begin{split}
W(0) & = \dot{\gamma}(s) ,\\
(\overline{\nabla}_{\dot{c}} W) \bigg |_{\tau = 0} & = 0.
\end{split}
\end{equation}
So, we just have to show that the vector field $W^{(s)}(\tau )$ along $c_s$ as \emph{defined} by $W^{(s)} (\tau ) = \cosh (ar) e_2(\Phi (\tau ,s ))$ is also the Jacobi field satisfying the two initial conditions specified in \eqref{initialvalues}. Once this is done, the uniqueness property of Jacobi field will immediately give the validity of relation \eqref{coshgrowth}. Now, observe that we must have $\overline{\nabla}_{\dot{c_s}} e_2^{(s)} = 0$ holds everywhere on $c_s$, simply because $e_2^{(s)}$ is a parallel vector field along $c_s$. So, by means of the product rule
$\overline{\nabla}_{X}(f Y) = (Xf) Y + f \overline{\nabla}_{X}Y$, which is one of the characteristic properties of any covariant derivative, it follows that
\begin{equation}
\overline{\nabla}_{\dot{c_s}}W^{(s)} = a \sinh (ar) e_2^{(s)} ,
\end{equation}
from which we immediately deduce that $W^{(s)}$ satisfies the second initial condition as specified in \eqref{initialvalues}. Of course $W^{s}$ clearly also satisfies the first condition in \eqref{initialvalues}. Now, by taking one more covariant derivative $\overline{\nabla}_{\dot{c}(s)}$ on both sides of the above equation and using the fact that $\overline{\nabla}_{\dot{c_s}} e_2^{(s)} = 0$, we immediately obtain

\begin{equation}
\overline{\nabla}_{\dot{c_s}} \overline{\nabla}_{\dot{c_s}}W^{(s)} = a^2 \cosh (ar) e_2^{(s)} = a^2 W^{(s)},
\end{equation}
from which we see immediately that $W^{(s)}$ clearly satisfies the Jacobi-field equation \eqref{Jacobifieldshyper}. As a result, $W^{(s)}$ is really a Jacobi field along $c_s$ which satisfies the same initial conditions as the Jacobi field $\frac{\partial}{\partial s}$ along $c_s$. Hence, by uniqueness, relation \eqref{coshgrowth} is true.

\section{The proof of Theorem \ref{parallelflowhyperbolicThm}: about parallel laminar flow along a straight edge of an obstacle in $\mathbb{H}^2(-a^2)$}\label{ProofofHypstraightSec}

\noindent
The goal of this section is to study parallel laminar flows along a geodesic, which represents the boundary (straight edge) of an obstacle, in the $2$-dimensional space form $\mathbb{H}^2(-a^2)$ with constant sectional curvature $-a^2$.\\

\noindent
Now, let us consider the ''Cartesian coordinate system'' $\Phi : \mathbb{R}^2 \rightarrow \mathbb{H}^2(-a^2)$ as given in \eqref{coordinatesystemhyper}, which we have just constructed in the previous section. For each point $p$, $\tau (p)$ and $s(p)$ stand for the first and second components of $\Phi^{-1}(p)$ in $\mathbb{R}^2$ respectively. That is, we have $\Phi^{-1}(p) = (\tau (p) , s(p))$. Then, we will consider the following solid region
\begin{equation}\label{K}
K = \{ p \in \mathbb{H}^2(-a^2) :  \tau (p) \leqslant 0 \} ,
\end{equation}
which will represent a solid obstacle with the straight edge $\partial K$ along which we study parallel laminar flow under the ''no-slip'' condition. According to the setting as given in \textbf{Section \ref{Cartesianhyperbolicsection}}, the straight edge $\partial K$ is exactly the geodesic $\gamma$, which represents the ''$y$-axis'' of the ''Cartesian coordinate system'' $\Phi (\tau , s)$ on $\mathbb{H}^2(-a^2)$. Now, we consider the following vector field as defined on $\mathbb{H}^2(-a^2) - K = \{ p \in \mathbb{H}^2(-a^2) \tau (p) > 0 \}$.
\begin{equation}\label{parallelhypone}
u (p) = - h (\tau (p)) e_{2}(p) = - h(\tau (p)) \frac{1}{\cosh (a\tau (p))} \frac{\partial}{\partial s}\bigg |_{p} .
\end{equation}
From now on, the two real-valued \emph{smooth functions} $p \to \tau (p)$, and $p \to s(p)$ on $\mathbb{H}^2(-a^2)$ will simply be denoted by $\tau$ and $s$ respectively. Here, we remind our readers that, from now on, the symbols $\tau$, and $s$ stand for the first- and second-component functions of the map $\Phi^{-1} : \mathbb{H}^2(-a^2) \rightarrow \mathbb{R}^2$. (So, our readers should not confuse them with the use of the same symbols ''$\tau $'', and ''$s$'' in $\mathbb{R}$ which represents parameters of the geodesics $c_{s}$, and $\gamma$ as in the previous section). With this convention for our notations, we can simply just write expression \eqref{parallelhypone} in the following ''short-hand'' form.
\begin{equation}\label{parallelhyptwo}
u = h(\tau ) e_2 = -h(\tau ) \frac{1}{\cosh (a\tau )} \frac{\partial}{\partial s} .
\end{equation}
Again, let $e_{2}$ to be the globally defined smooth unit vector field on $\mathbb{H}^2(-a^2)$ which we construct in \textbf{Section \ref{Cartesianhyperbolicsection}}. Recall, that, in terms of the notations as specified in the previous section, $e_{2}$ is everywhere orthonormal to the unit vector field $\frac{\partial}{\partial \tau } = \dot{c_s}$. Now, let us denote the vector field $\frac{\partial}{\partial \tau}$ on $\mathbb{H}^2(-a^2)$ by the symbol $e_1$. On the other hand, we recall, form \eqref{coshgrowth} of the previous section, that we have
$e_{2} = \frac{1}{\cosh (a\tau )} \frac{\partial}{\partial s}$. Then, in accordance with our construction in the previous section, we know that
\begin{itemize}
\item $e_1 = \frac{\partial}{\partial \tau}$ and $e_2 = \frac{1}{\cosh (a\tau )} \frac{\partial}{\partial s}$ constitute a positively oriented orthonormal moving frame on $\mathbb{H}^2(-a^2)$.
\end{itemize}
Then, it follows that the associated $1$-forms $e_1^* = d\tau$, and $e_2^* = \cosh (a\tau ) ds$ will constitute an orthonormal coframe with respect to the induced Riemannian metric on the cotangent bundle $T^*\mathbb{H}^2(-a^2)$ of $\mathbb{H}^2(-a^2)$. Then, the volume form on $\mathbb{H}^2(-a^2)$ can be expressed as
\begin{equation}\label{volumehyptwo}
Vol_{\mathbb{H}^2(-a^2)} = e_1^*\wedge e_2^* = \cosh (a\tau ) d\tau \wedge ds .
\end{equation}
Then, the Hodge-Star operator $* : C^{\infty}(\wedge^2T^* \mathbb{H}^2(-a^2)) \rightarrow C^{\infty}(\mathbb{H}^2(-a^2))$ sending $2$-forms back to the space of smooth functions is characterized by the \emph{defining} relation $* Vol_{\mathbb{H}^2(-a^2)} = 1 $, and the tensorial property $*(f\omega ) = f *\omega$, with $f$ to be a smooth function and $\omega$ to be a differential form. We also recall that the Hodge-star operator
$* C^{\infty}(T^* \mathbb{H}^2(-a^2)) \rightarrow C^{\infty}(T^* \mathbb{H}^2(-a^2))$ sending $1$-forms back to the space of $1$-forms can now be represented by the following rules.
\begin{equation}\label{rotationhyper}
\begin{split}
*(d\tau ) = *e_1^* & = e_2^* = \cosh (a\tau ) ds \\
*e_2^* = -e_1^* &= -d\tau .
\end{split}
\end{equation}
Now, we can proceed to compute $(-\triangle )u^*$, where $u^* = g(u, \cdot )$ is the associated $1$-form of the vector field $u$ as specified in \eqref{parallelhyptwo}.

\noindent{\bf Step $1$: Checking the divergence free property $d^*u^* = 0$ for $u$ as given in \eqref{parallelhyptwo}.}
Recall that the operator $d^* = -*d*$ which sends $1$-forms on $\mathbb{H}^2(-a^2)$ to the space of smooth functions on $\mathbb{H}^2(-a^2)$ is just the operator $- div$ acting on the space of smooth vector fields on $\mathbb{H}^2(-a^2)$. So, the desired divergence free property for $u$ as given in \eqref{parallelhyptwo} is expressed in the form of $d^* u^* =0$, which can easily be obtained through the following computation.
\begin{equation}
\begin{split}
d^*u^* & = *d[h(\tau ) * e_2^*] \\
& = - * d [h(\tau ) d\tau ] \\
& = - * \left\{ \frac{\partial}{\partial s}(h(\tau )) ds \wedge d\tau   \right\} \\
& = 0.
\end{split}
\end{equation}

\noindent{\bf Step $2$: The computation of $(-\triangle ) u^*$ for the velocity field $u$ given by \eqref{parallelhyptwo}.}

\noindent
Recall that the Hodge Laplacian $(-\triangle )$, which sends the space of smooth $1$-forms into itself, is given by $(-\triangle ) = dd^* + d^* d$. Since the velocity field $u$ as given in \eqref{parallelhyptwo} satisfies the divergence free property $d^* u^* = 0$ on $\mathbb{H}^2(-a^2) - K$, with $K$ to be the solid obstacle with a straight edge boundary as specified in \eqref{K}, it follows that we have the following relation.
\begin{equation}
(-\triangle )u^* = d^* du .
\end{equation}
Now, we first carry out the following computation for the $2$-form $du^*$, which represents the vorticity of $u$ on $\mathbb{H}^2(-a^2)-K$.
\begin{equation}
\begin{split}
du^* & = d \{ - h(\tau ) \cosh (a\tau ) ds \} \\
& = - \frac{\partial}{\partial \tau } [h(\tau ) \cosh (a\tau )] d\tau \wedge ds \\
& = -\frac{1}{\cosh (a\tau )} \frac{\partial}{\partial \tau } [h(\tau ) \cosh (a\tau )] Vol_{\mathbb{H}^2(-a^2)} ,
\end{split}
\end{equation}
where $Vol_{\mathbb{H}^2(-a^2)}$ is the volume-form on $\mathbb{H}^2(-a^2)$ as described in \eqref{volumehyptwo}. Since $* Vol_{\mathbb{H}^2(-a^2)} = 1$, it follows that
\begin{equation}\label{Viscosityhyper}
\begin{split}
(-\triangle )u^* & = d^* du^* \\
& = -*d* \bigg \{ -\frac{1}{\cosh (a\tau )} \frac{\partial}{\partial \tau } [h(\tau ) \cosh (a\tau )] Vol_{\mathbb{H}^2(-a^2)}  \bigg \} \\
& = * d \bigg \{ \frac{1}{\cosh (a\tau )} \frac{\partial}{\partial \tau } [h(\tau ) \cosh (a\tau )] \bigg \} \\
& = * d \bigg \{  h'(\tau ) + a h(\tau ) \cdot \frac{\sinh (a\tau )}{\cosh (a\tau )}   \big \} \\
& = * \frac{\partial}{\partial \tau } \bigg \{  h'(\tau ) + a h(\tau ) \cdot \frac{\sinh (a\tau )}{\cosh (a\tau )}       \bigg \} d \tau \\
& = \frac{\partial}{\partial \tau } \bigg \{  h'(\tau ) + a h(\tau ) \cdot \frac{\sinh (a\tau )}{\cosh (a\tau )}       \bigg \} * d\tau \\
& = \frac{\partial}{\partial \tau } \bigg \{  h'(\tau ) + a h(\tau ) \cdot \frac{\sinh (a\tau )}{\cosh (a\tau )}       \bigg \} \cosh (a\tau ) ds \\
& = \bigg \{  h''(\tau ) \cosh (a\tau ) + a \sinh (a \tau ) h'(\tau ) + \frac{a^2}{\cosh (a\tau )} h(\tau )           \bigg \} ds .
\end{split}
\end{equation}
In the above computation, the symbol $h'(\tau )$ stands for $h'(\tau ) = \frac{\partial h}{\partial \tau }$. In the same way, $h''(\tau )$ and $h'''(\tau)$ mean $h''(\tau ) = \big (\frac{\partial}{\partial \tau }\big )^2 h$ and $h'''(\tau ) = \big ( \frac{\partial}{\partial \tau }\big )^3 h$ respectively. In the sixth equality of the above calculation, we have used the basic tensorial property $* (f \omega ) = f * (\omega )$ of the Hodge-Star operator $*$, with $f$ to be a smooth function and $\omega$ to be a smooth $1$-form.

\noindent{\bf Step 3: the computation of the nonlinear convection term $(\overline{\nabla}_{u}u)^*$, for $u$ to be given in \eqref{parallelhyptwo}.}

\noindent
In the computation for the nonlinear convection term $(\overline{\nabla}_{u}u)^*$, it is convenient for us to work directly with the computation of  $\overline{\nabla}_{u}u$ at the level of smooth vector fields. Recall that the natural vector fields $\frac{\partial}{\partial \tau }$ and
$\frac{\partial}{\partial s }$ are everywhere orthogonal to each other on $\mathbb{H}^2(-a^2)$. So, we now express the vector field $\overline{\nabla}_{u}u$ in terms of the linear combination of $\frac{\partial}{\partial \tau }$ and
$\frac{\partial}{\partial s }$ as follow.
\begin{equation}\label{LeviCivitahyper}
\overline{\nabla}_{u}u = A \frac{\partial}{\partial \tau } + B \frac{\partial}{\partial s },
\end{equation}
where $A$, and $B$ are some smooth functions on $\mathbb{H}^2(-a^2)$ which we will figure out in a minute. First, by taking the inner product with $\frac{\partial}{\partial s }$ on both sides of \eqref{LeviCivitahyper}, we deduce that
\begin{equation}
\begin{split}
\cosh^2(a\tau ) B & = g(\overline{\nabla}_{u}u , \frac{\partial}{\partial s } ) \\
& =  - \frac{\cosh (a\tau )}{h(\tau )} g(\overline{\nabla}_{u}u , u  )\\
& = - \frac{\cosh (a\tau )}{h(\tau )} \cdot \frac{1}{2} u \big ( |u|^2 \big ) \\
& = \frac{1}{2}  \frac{\partial}{\partial s } \big [ (h(\tau ))^2 \big ] = 0,
\end{split}
\end{equation}
from which  we immediately get $B = 0$ on $\mathbb{H}^2(-a^2) - K$. We remark that in the above computation, the symbol $u \big ( |u|^2 \big )$ stands for \emph{the derivative of the function $|u|^2$ along the direction of the vector field $u$}.  On the other hand, by taking inner product with $\frac{\partial}{\partial \tau }$ on both sides of \eqref{LeviCivitahyper}, we can compute the smooth function $A$ as follow.
\begin{equation}
A  = g( \overline{\nabla}_{u}u , \frac{\partial}{\partial \tau } )
 = \frac{-h(\tau )}{\cosh (a\tau )} g(  \overline{\nabla}_{\frac{\partial}{\partial s }}u , \frac{\partial}{\partial \tau }   )
\end{equation}
However, since $0= g(u , \frac{\partial}{\partial \tau } )$ holds on $\mathbb{H}^2(-a^2) - K$, it follows that we have
\begin{equation}
0 = \frac{\partial}{\partial s } g(u , \frac{\partial}{\partial \tau } ) = g(\overline{\nabla}_{\frac{\partial}{\partial s }} u , \frac{\partial}{\partial \tau } ) + g(u , \overline{\nabla}_{\frac{\partial}{\partial s }}  \frac{\partial}{\partial \tau } ) .
\end{equation}
Hence, we can resume the calculation for $A$ as follow.
\begin{equation}\label{C1hyper}
A  = \frac{-h(\tau )}{\cosh (a\tau )} g(  \overline{\nabla}_{\frac{\partial}{\partial s }}u , \frac{\partial}{\partial \tau }   )
= \frac{h(\tau )}{\cosh (a\tau )}  g(u , \overline{\nabla}_{\frac{\partial}{\partial s }}  \frac{\partial}{\partial \tau } )
= \frac{h(\tau )}{\cosh (a\tau )}  g(u , \overline{\nabla}_{\frac{\partial}{\partial \tau }}  \frac{\partial}{\partial s } ) ,
\end{equation}
where in the last equality sign, we have used the property $\overline{\nabla}_{\frac{\partial}{\partial s }}  \frac{\partial}{\partial \tau } = \overline{\nabla}_{\frac{\partial}{\partial \tau }}  \frac{\partial}{\partial s }$, which follows from the torsion free property of the Levi-Civita connection. Also, according with the Leibniz's rule satisfied by the Levi-Civita connection, we can carry out the following computation for $\overline{\nabla}_{\frac{\partial}{\partial \tau }}  \frac{\partial}{\partial s }$, with $\overline{\nabla}_{\frac{\partial}{\partial \tau }}e_{2} = 0$ (which is true since $e_2$ is parallel along each geodesic $c_s$) being taken into account in the calculation.
\begin{equation}\label{C2hyper}
\overline{\nabla}_{\frac{\partial}{\partial \tau }}  \frac{\partial}{\partial s } =
\overline{\nabla}_{\frac{\partial}{\partial \tau }}  \big [  \cosh (a \tau ) e_2   \big ] =
a\sinh (a\tau ) e_2 =
\frac{a \sinh(a\tau )}{\cosh(a\tau )} \frac{\partial}{\partial s } .
\end{equation}
Hence, it follows from \eqref{C1hyper}, and \eqref{C2hyper} that
\begin{equation}
A  = \frac{h(\tau )}{\cosh (a\tau )} g\big ( - h(\tau ) e_2 , a \sinh(a\tau ) e_2 \big ) \\
 = \frac{-a \sinh(a\tau )}{\cosh (a\tau )} (h(\tau ))^2 .
\end{equation}
As a result, we finally conclude that
\begin{equation}
\overline{\nabla}_u u = \frac{-a \sinh(a\tau )}{\cosh (a\tau )} (h(\tau ))^2 \frac{\partial}{\partial \tau } ,
\end{equation}
with associated $1$-form $\big ( \overline{\nabla}_u u \big )^*$ to be given by
\begin{equation}\label{nonlinearhyper}
\big ( \overline{\nabla}_u u \big )^* = \frac{-a \sinh(a\tau )}{\cosh (a\tau )} (h(\tau ))^2 d\tau  .
\end{equation}
Now, for our forthcoming applications of these calculations in the proof of Theorem \ref{parallelflowhyperbolicThm}, we will need concrete expressions of $d \big ( (-\triangle )u^* -2 Ric (u^*)   \big )$ and $d \big (  \overline{\nabla}_uu   \big )^*$. Indeed, by taking into account of the basic fact in Riemannian Geometry that $Ric (X^*) = -a^2 X^*$ always holds for any smooth vector fields $X$ on $\mathbb{H}^2(-a^2)$, we can use expression \eqref{Viscosityhyper} to derive the following expression for $d \big ( (-\triangle )u^* -2 Ric (u^*)   \big )$.
\begin{equation}\label{vorticitytermhyper}
\begin{split}
&d \big \{(-\triangle )u^* -2Ric (u^*) \big \} \\
&= d\bigg \{  h''(\tau ) \cosh (a\tau ) + a \sinh (a \tau ) h'(\tau ) + a^2 \bigg ( \frac{1}{\cosh (a\tau )} - 2 \cosh (a\tau )\bigg ) h(\tau )  \bigg \} \wedge ds \\
&= \bigg \{ h'''(\tau ) \cosh (a\tau ) + 2a \sinh (a\tau ) h''(\tau ) - \frac{a^2\sinh^2(a\tau)}{\cosh(a\tau)}h'(\tau ) \\
& - a^3\sinh(a\tau ) \bigg (   2 + \frac{1}{\cosh^2(a\tau )} \bigg ) h(\tau ) \bigg \} d\tau \wedge ds .
\end{split}
\end{equation}
Also, by taking the operator $d$ on both sides of \eqref{nonlinearhyper}, we immediately obtain the following relation,
\begin{equation}\label{TRIVIALhyper}
d \big ( \overline{\nabla}_u u \big )^* = \frac{\partial}{\partial s}\bigg \{ \frac{-a \sinh(a\tau )}{\cosh (a\tau )} (h(\tau ))^2 \bigg \} ds \wedge d\tau = 0.
\end{equation}

\noindent{\bf Step 4 : The proof of \textbf{Assertion I} in Theorem \ref{parallelflowhyperbolicThm}}\\

\noindent
With the preparations in the previous steps of this section, we are now ready to give a proof for \textbf{Assertion I} in Theorem \ref{parallelflowhyperbolicThm} as follow.\\
To begin the argument, assume towards contradiction that the parallel laminar flow $u = - h(\tau ) \frac{1}{\cosh (a\tau )} \frac{\partial}{\partial s}$, with the quadratic profile $h(\tau ) = \alpha_1 \tau -\frac{\alpha_2}{2}\tau^2$, does satisfy equation \eqref{NavierStokeshyperbolic} on the simply connected region $\Omega_{\tau_0} = \Phi \big ( (0,\tau_0 )\times \mathbb{R} \big )$ of $\mathbb{H}^2(-a^2)$, for certain positive constants $\tau_0 > 0$, $\alpha_1 > 0$, and $\alpha_2 > 0$. That is, the associated one form $u^* = -h(\tau ) \cosh (a\tau ) ds$ will satisfies the following Stationary Navier-Stokes equation on $\Omega_{\tau_0} = \Phi \big ( (0,\tau_0 )\times \mathbb{R} \big )$.
\begin{equation}\label{NShyperbolictwo}
\nu ((-\triangle)u^* -2 Ric (u^*) ) +\overline{\nabla}_{u}u^* + dP = 0 .
\end{equation}
Now, by taking the exterior differential operator $d$ on both sides of \eqref{NShyperbolictwo}, we deduce from \eqref{vorticitytermhyper} and \eqref{TRIVIALhyper} that the following vorticity equation would hold on $\Omega_{\tau_0} = \Phi \big ( (0,\tau_0 ) \times \mathbb{R} \big )$.
\begin{equation}\label{vorticityequationhyperplane}
\begin{split}
0 & = d \big \{ (-\triangle)u^* -2 Ric (u^*)  \big \} \\
& = \bigg \{ h'''(\tau ) \cosh (a\tau ) + 2a\sinh (a\tau ) h''(\tau ) \\
&- \frac{a^2\sinh^2(a\tau )}{\cosh (a\tau )} h'(\tau )
- a^3 \sinh (a\tau ) \bigg ( 2+ \frac{1}{\cosh^2(a\tau ) }\bigg ) h(\tau )            \bigg \} d\tau \wedge ds .
\end{split}
\end{equation}
For convenience, we will consider the following smooth function
\begin{equation}\label{GoodFunctionhyperplane}
\begin{split}
F(\tau ) & =  \bigg \{ 2a\sinh (a\tau ) h''(\tau )
 - \frac{a^2\sinh^2(a\tau )}{\cosh (a\tau )} h'(\tau )
- a^3 \sinh (a\tau ) \bigg ( 2+ \frac{1}{\cosh^2(a\tau ) }\bigg ) h(\tau )            \bigg \} .
\end{split}
\end{equation}
Since for $h(\tau ) = \alpha_1 \tau -\frac{\alpha_2}{2}\tau^2$, we have $h'''(\tau ) = 0$ for all $\tau \in \mathbb{R}$, it follows that, for the velocity field $u = - h(\tau ) \frac{1}{\cosh (a\tau )} \frac{\partial}{\partial s}$, with the quadratic profile $h(\tau ) = \alpha_1 \tau -\frac{\alpha_2}{2}\tau^2$, equation \eqref{vorticityequationhyperplane} on $\Omega_{\tau_0} = \Phi \big ( (0,\tau_0 )\times \mathbb{R} \big ) $ should be equivalent to the following equation, which would hold on $\Omega_{\tau_0} = \Phi \big ( (0,\tau_0 )\times \mathbb{R} \big )$.
\begin{equation}\label{VorticityeqContradiction}
0 = F(\tau ) d\tau \wedge ds ,
\end{equation}
where $F(\tau )$ is the real-analytic function as given in \eqref{GoodFunctionhyperplane}. Since the differential $2$-form $d\tau \wedge ds$ is everywhere non-vanishing on $\mathbb{H}^2(-a^2)$, the validity of equation \eqref{VorticityeqContradiction} on $\Omega_{\tau_0} = \Phi \big ( (0,\tau_0 )\times \mathbb{R} \big )$ is equivalent to saying that
\begin{itemize}
\item The real-analytic function $F(\tau )$ as given in \eqref{GoodFunctionhyperplane} should totally vanish over $\Omega_{\tau_0} = \Phi \big ( (0,\tau_0 ) \times \mathbb{R} \big )$. That is, we should have $F (\tau ) = 0$, for all $\tau \in [0, \tau_0 )$.
\end{itemize}
So, to finish the proof for \textbf{Assertion I} of Theorem \ref{parallelflowhyperbolicThm}, we just need to arrive at a contradiction against the \emph{everywhere vanishing property of $F(\tau )$ over the interval $[0, \tau_{0})$}. To achieve this, we simply just compute the term $\frac{\partial F}{\partial \tau} \big |_{\tau = 0} = F'(0)$ as follow. Indeed, for the quadratic profile $h(\tau ) = \alpha_1 \tau -\frac{\alpha_2}{2}\tau^2$, we can carry out the following straightforward computations.
\begin{equation}\label{calculuscomphyper}
\begin{split}
\frac{\partial}{\partial \tau }\bigg \{ 2a\sinh (a\tau )h''(\tau )   \bigg \} \bigg |_{\tau = 0} & = 2a^2h''(0) = -2a^2\alpha_2 , \\
\frac{\partial}{\partial \tau }\bigg \{ - \frac{a^2\sinh^2(a\tau )}{\cosh (a\tau )} h'(\tau )  \bigg \}  \bigg |_{\tau = 0} & = 0 ,\\
\frac{\partial}{\partial \tau }\bigg \{ - a^3 \sinh (a\tau ) \bigg ( 2+ \frac{1}{\cosh^2(a\tau ) }\bigg ) h(\tau )   \bigg \} \bigg |_{\tau = 0} & = 0 ,
\end{split}
\end{equation}
where in the above computation, we have taken the information $h(0)=0$, $h''(0) = -\alpha_2$, into account. So, in light of \eqref{calculuscomphyper}, we can easily deduce from \eqref{GoodFunctionhyperplane} that we must have
\begin{equation}
\frac{\partial F}{\partial \tau} \big |_{\tau = 0} = F'(0) = -2a^2 \alpha_2 < 0 ,
\end{equation}
since $\alpha_2 > 0$ by our hypothesis. The above relation clearly indicates that the function $F(\tau )$ must be a strictly decreasing function in some small open interval $(-\epsilon_0 , \epsilon_0 )$ around the point $\tau =0$. This fact, together with $F(0) = 0$, will imply that the function $F(\tau )$ must be strictly negative, for all $\tau \in (0, \epsilon_0 )$. This clearly violates the everywhere vanishing property: $F(\tau ) = 0$, for all $\tau \in [0, \tau_0 )$. So, we have derive a contradiction against the validity of the vorticity equation \eqref{VorticityeqContradiction}. Hence, we conclude that the velocity field $u = - h(\tau ) \frac{1}{\cosh (a\tau )} \frac{\partial}{\partial s}$, with the quadratic profile $h(\tau ) = \alpha_1 \tau -\frac{\alpha_2}{2}\tau^2$ is not a solution to equation \eqref{NavierStokeshyperbolic} on $\Omega_{\tau_0} = \Phi \big ((0,\tau_0 )\times \mathbb{R} \big )$, \emph{no matter which positive constants $\tau_0$, $\alpha_1 > 0$, $\alpha_2 > 0$ we take}.\\

\noindent{\bf Step $5$: The proof of \textbf{Assertion II} in Theorem \ref{parallelflowhyperbolicThm}}\\

\noindent
To begin the proof of \textbf{Assertion II} in Theorem \ref{parallelflowhyperbolicThm}, recall that all the computations in \textbf{Section \ref{ProofofHypstraightSec}} up to \eqref{vorticitytermhyper} are valid for a velocity field $u = -Y(\tau ) \frac{1}{\cosh (a\tau )}\frac{\partial}{\partial s}$, with $Y (\tau )$ to be any smooth function on $[0,\infty )$. Now, we are interested in the question of whether a velocity field of the form $u = -Y(\tau ) \frac{1}{\cosh (a\tau )}\frac{\partial}{\partial s}$, with $Y :[0,\infty )\rightarrow \mathbb{R}$ to be a smooth function, is a solution to equation \eqref{NavierStokeshyperbolic} on the whole exterior region $\mathbb{H}^2(-a^2)-K = \{ \Phi (\tau , s) \in \mathbb{H}^2(-a^2) : \tau > 0 , s \in \mathbb{R}  \}$, with the obstacle $K = \{\Phi (\tau ,s ) \in \mathbb{H}^2(-a^2) : \tau \leqslant 0 , s \in \mathbb{R}    \}$. Since $ \mathbb{H}^2(-a^2)-K = \Phi \big ( (0,\infty )\times \mathbb{R} \big )$ is clearly simply-connected, there exists a globally defined smooth function $P$ on $\mathbb{H}^2(-a^2)-K$ which satisfies equation \ref{NavierStokeshyperbolic} on $\mathbb{H}^2(-a^2)-K$ \emph{if and only if} $u = -Y(\tau ) \frac{1}{\cosh (a\tau )}\frac{\partial}{\partial s}$ satisfies the following vorticity equation on $\mathbb{H}^2(-a^2)-K$.
\begin{equation}\label{vorticityhyperLASTeq}
d \big \{ \nu \big ( (-\triangle ) u^* -2 Ric (u^*)     \big ) + \overline{\nabla}_{u}u^*   \big \} = 0 ,
\end{equation}
which, in accordance with \eqref{vorticitytermhyper} and \eqref{TRIVIALhyper}, is equivalent to the following third-order O.D.E. with real-analytic coefficients in the variable $\tau \in \mathbb{R}$.
\begin{equation}\label{3rdODEhyperplane}
\begin{split}
0 & = Y'''(\tau ) \cosh (a\tau ) + 2a\sinh (a\tau ) Y''(\tau ) \\
& - \frac{a^2\sinh^2(a\tau )}{\cosh (a\tau )} Y'(\tau )
- a^3 \sinh (a\tau ) \bigg ( 2+ \frac{1}{\cosh^2(a\tau ) }\bigg ) Y(\tau ).
\end{split}
\end{equation}
In other words, $u = -Y(\tau ) \frac{1}{\cosh (a\tau )}\frac{\partial}{\partial s}$ will satisfy equation \eqref{NavierStokeshyperbolic} on the simply-connected region $\mathbb{H}^2(-a^2)-K$ with some globally defined pressure $P$ on $\mathbb{H}^2(-a^2)-K$ \emph{if and only if} the smooth function $Y \in C^{\infty }\big ( [0,\infty ) \big )$ is a solution to the third-order ODE \eqref{3rdODEhyperplane}. However, in accordance with the basic existence theorem in the theory of linear O.D.E. $($i.e. Theorem \ref{ODEtheorem}$)$, we deduce that for any prescribed positive constants $\alpha_1 > 0$, and $\alpha_2 > 0$, there exists a unique smooth solution $Y \in C^{\infty}\big ( [0,\infty ) \big )$ to \eqref{3rdODEhyperplane}, which satisfies the initial values $Y(0) =0$, $Y'(0) = \alpha_1$, and $Y''(0) = -\alpha_2$. Since the coefficients in the third-order O.D.E. \eqref{3rdODEhyperplane} are all real analytic on $\mathbb{R}$, it turns out that such a unique solution $Y \in C^{\infty}\big ( [0,\infty ) \big )$ to \eqref{3rdODEhyperplane} must itself also be real-analytic on $[0,\infty )$. With such a real-analytic solution $Y$ to \eqref{3rdODEhyperplane} to be available to us, we can now conclude that: for any given $\alpha_1 > 0$, and $\alpha_2 > 0$, the real-analytic solution $Y :[0,\infty ) \rightarrow \mathbb{R}$  to \eqref{3rdODEhyperplane}, satisfying the initial values $Y(0) = 0$, $Y'(0) = \alpha_1$, $Y''(0) = - \alpha_2$, is the one and only one smooth function on $[0,\infty )$ for which the velocity field $u = -Y(\tau ) \frac{1}{\cosh (a\tau )}\frac{\partial}{\partial s}$ will satisfy  equation \eqref{NavierStokeshyperbolic} $($with some globally defined pressure function $P$$)$ on the simply-connected open region $\mathbb{H}^2(-a^2)-K = \Phi \big( (0,\infty ) \times \mathbb{R} \big )$. So, we are done in proving \textbf{Assertion II} of Theorem \ref{parallelflowhyperbolicThm}.






\bibliography{ChiHinMagda}
\bibliographystyle{plain}

\end{document}